\documentclass[11pt]{amsart}

\usepackage{amscd,amsmath,amssymb}
\usepackage[mathscr]{eucal}

\makeatletter \@addtoreset{equation}{section}\makeatother

\newtheorem{theorem}{Theorem}[section]
\newtheorem{lemma}[theorem]{Lemma}
\newtheorem{proposition}[theorem]{Proposition}
\newtheorem{corollary}[theorem]{Corollary}

\newtheorem{remark}[theorem]{Remark}

\newcommand{\Z}{{\mathbb{Z}}}
\newcommand{\ZZ}{{\mathbb{Z}/2}}
\newcommand{\C}{{\mathbb{C}}}
\newcommand{\hH}{{\widehat{H}}}
\newcommand{\R}{{\mathbb{R}}}
\newcommand{\lm}{{\mu}}
\newcommand{\ld}{{\delta}}
\newcommand{\bd}{{b(\delta)}}
\newcommand{\bm}{{b(\mu)}}
\newcommand{\bw}{{b^e(w)}}
\newcommand{\bdw}{{b^e(D(w))}}
\newcommand{\bwne}{{b(w)}}
\newcommand{\bdwne}{{b(D(w))}}
\newcommand{\bed}{{b^e(\delta)}}
\newcommand{\bem}{{b^e(\mu)}}
\newcommand{\bmv}{{b^v(\mu)}}
\newcommand{\bmh}{{b^h(\mu)}}
\newcommand{\h}{{h}}
\newcommand{\End}{{\mathrm{End}}}
\newcommand{\rC}{{\mathsf{C}}}
\newcommand{\rHP}{{\mathsf{HP}}}
\newcommand{\rrC}{{\mathsf{C}^\Pi}}
\newcommand{\rrCe}{{\mathsf{C}^{\Pi,\,e}}}
\newcommand{\str}{{\mathsf{str}}}

\newcommand{\rHH}{{\mathsf{HH}}}

\newcommand{\ID}{{b^e(a)}}
\newcommand{\JD}{{b^e(a')}}
\newcommand{\IdD}{{b^e(da)}}
\newcommand{\dD}{{\mathrm{ad}(a)}}
\newcommand{\lD}{{l(a)}}
\newcommand{\rD}{{r(a)}}
\newcommand{\ldw}{{l(D(w))}}
\newcommand{\rdw}{{r(D(w))}}
\newcommand{\lw}{{l(w)}}
\newcommand{\rw}{{r(w)}}

\newcommand{\cG}{{\mathcal{G}}}
\newcommand{\cD}{{\mathscr{D}}}
\newcommand{\Q}{{\mathbb{Q}}}
\newcommand{\cH}{{\mathcal{H}}}

\newcommand{\DD}{{\mathrm{End}_\C V}}

\newcommand{\cA}{{\mathscr{A}}}
\newcommand{\bK}{{\boldsymbol{k}}}

\newcommand{\Udl}{{\mathscr{U}(\delta)}}
\newcommand{\Vdl}{{\mathscr{V}(\delta)}}
\newcommand{\Uml}{{\mathscr{U}({\mu})}}
\newcommand{\Vml}{{\mathscr{V}({\mu})}}

\newcommand{\Ud}{{\mathscr{U}^\circ(\delta)}}
\newcommand{\Vd}{{\mathscr{V}^\circ(\delta)}}

\newcommand{\nUd}{{\mathscr{U}^{un}(\delta)}}
\newcommand{\nVd}{{\mathscr{V}^{un}(\delta)}}
\newcommand{\nWd}{{\mathscr{W}^{un}(\delta)}}

\newcommand{\Uwc}{{\mathscr{U}(w)}}
\newcommand{\Vwc}{{\mathscr{V}(w)}}
\newcommand{\Umc}{{\mathscr{U}({\mu})}}
\newcommand{\Vmc}{{\mathscr{V}({\mu})}}

\title[NC Hodge structures: matching categorical and geometric examples]{Non-commutative Hodge structures: Towards matching categorical and geometric examples}
\author{D. Shklyarov}
\address{Lehrstuhl f\"ur Analysis und Geometrie, Universit\"at Augsburg, Institut f\"ur Mathematik,
86135 Augsburg, Germany}
\email{dmytro.shklyarov@math.uni-augsburg.de}
\thanks{This research was supported by the ERC Starting Independent Researcher Grant StG No. 204757-TQFT (K.~Wendland PI)}
\begin{document}
\begin{abstract}
The subject of the present work is the de Rham part of non-commutative Hodge structures on the periodic cyclic homology of differential graded algebras and categories. We discuss explicit formulas for the corresponding connection on the periodic cyclic homology viewed as a bundle over the punctured formal disk. Our main result says that for the category of matrix factorizations of a polynomial the formulas reproduce, up to a certain shift, a well-known connection on the associated twisted de Rham cohomology which plays a central role in the geometric approach to the Hodge theory of isolated singularities.
\end{abstract}
\maketitle
\tableofcontents
\baselineskip 1.5pc

\section{Introduction}

The study of various generalizations of the classical Hodge theory has evolved into a vast branch of mathematics; a comprehensive overview of relevant notions and references can be found in \cite{Sab}. The present work is in the framework of the Hodge theory of categories, as outlined in \cite{Bar1,Bar2,KKP,Kon,KS}. We refer the reader to these sources for an introduction to the circle of ideas surrounding the subject, as well as motivation and references. Our main goal is to present a piece of evidence (see Theorem \ref{intro2} below) in favor of the idea that the categorical Hodge theory should include some known geometric examples of generalized Hodge structures as special cases.

The periodic cyclic homology of a differential graded (dg) category is well-known to be a direct generalization of the de Rham cohomology of a space, and it is natural to expect it to carry a Hodge-like structure. This expectation has been converted into a precise conjecture in \cite[Section 2.2.2]{KKP}, at least in the case of the so-called homologically smooth and proper dg categories which are to be thought of as analogs of smooth and proper varieties. According to the conjecture, the cyclic homology of such a category can be endowed with a {\it non-commutative Hodge structure}. Despite the terminology, non-commutative Hodge structures are defined without referring to non-commutative mathematics; they are a subject of their own, with a number of geometric applications \cite{HS,Iri,KKP}.

The definition of a non-commutative Hodge structure involves two sets of data called the ``de Rham data'' and the ``Betti data'' in \cite{KKP}. Roughly, the former generalizes the Hodge filtration of a classical Hodge structure (cf. section \ref{NCHodge}) while the latter is an analog of the $\Q$-structure. The main obstacle to proving the aforementioned conjecture seems to have been the fact that the periodic cyclic homology of, say, a $\C$-linear dg category does not carry an obvious $\Q$-structure in general. More on this can be found in \cite[Section 2.2.6]{KKP} and \cite[Section 8]{Ka}. The present work concerns the de Rham part of the sought-for non-commutative Hodge structure on the cyclic homology whose origin has been understood for some time now.

The de Rham data of a non-commutative Hodge structure can be defined as a pair $(\cH, \nabla)$ where $\cH$ is a $\ZZ$-graded free $\C\{u\}$-module of finite rank and $\nabla$ is a meromorphic connection on $\cH$ with a pole of order at most 2 at the origin. There is a formal counterpart of such data in which the ring $\C\{u\}$ of convergent series is replaced with the ring $\C[[u]]$ of formal series. We will be considering only the formal analog of the de Rham data due to the very nature of our main example, the cyclic homology.

In what follows, the term ``bundle on the formal (resp. punctured formal)  disk'' is used as a synonym of ``free $\C[[u]]$-module of finite rank (resp. finite-dimensional $\C((u))$-vector space)''. Accordingly, we will speak of ``bundles with connection (on the formal disk)'' instead of ``the de Rham data of  non-commutative Hodge structures''.

A general way of producing such bundles with connections was outlined in \cite[Section 4.2.4]{KKP}. Let us reproduce this idea here using three examples.

Let $X$ be a compact K\"{a}hler manifold and $\cA(X)$ the $\ZZ$-graded space of complex $C^\infty$-forms on $X$, with the $\ZZ$-grading given by the parity of the degree of differential forms. Consider the $\ZZ$-graded $\C[[u]]$-linear complex
$
(\cA(X)[[u]], \bar{\partial}+u{\partial})
$
where $\cA(X)[[u]]$ stands for the space of formal power series with coefficients in $\cA(X)$. By the classical Hodge theory, its cohomology $H^*(\cA(X)[[u]], \bar{\partial}+u{\partial})$ is a vector bundle on the formal disk of rank $\mathrm{dim}_\C\,H_{DR}(X)$. It carries a meromorphic connection defined as follows.
First, define $\nabla^X: \cA(X)((u))\to\cA(X)((u))$ by
$$
\nabla^X=\frac{d}{du}+\frac{\Gamma'}{2u}, \quad \Gamma'\in\End(\cA(X)), \quad \Gamma'|_{\cA^{p,q}(X)}=q-p
$$
and then observe that
\begin{equation}\label{1}
[\nabla^X,\bar{\partial}+u{\partial}]=\frac{1}{2u}(\bar{\partial}+u{\partial})
\end{equation}
Therefore, $\nabla^X$ induces a connection on the cohomology $H^*(\cA(X)((u)), \bar{\partial}+u{\partial})$, the latter being the restriction of the vector bundle $H^*(\cA(X)[[u]], \bar{\partial}+u{\partial})$ to the punctured formal disk.

For our second example, let us take a polynomial $w=w(y_1,\ldots,y_k)$ on $Y=\C^k$  such that
\begin{equation}\label{crit}
w(0,\ldots,0)=\partial_1w(0,\ldots,0)=\ldots=\partial_kw(0,\ldots,0)=0
\end{equation}
We will assume that the origin is the only critical point of $w$. To produce a bundle with connection associated with $w$ we will mimic the above construction for K\"{a}hler manifolds. Namely, let $\Omega(Y)$ be the $\ZZ$-graded space of (holomorphic) differential forms with polynomial coefficients on $Y$. Consider the so-called {\it twisted de Rham complex}
$
(\Omega(Y)[[u]], -dw+ud)
$
where $d$ is the (holomorphic) de Rham differential and $dw$ is the operator of wedge multiplication with $dw$
(as before, we view it as a $\ZZ$-graded $\C[[u]]$-linear complex). It is a classical fact that
$$
H^*(\Omega(Y)[[u]], -dw +ud)=\begin{cases} \Omega^k(Y)[[u]]/(-dw+ud)\Omega^{k-1}(Y)[[u]] & *=k\,\,\text{mod}\,\,2\\ 0 & \text{otherwise}\end{cases}
$$
and that the cohomology is a free $\C[[u]]$-module, i.e. a vector bundle on the formal disk. It carries a connection by the same argument as before:
the operator $\nabla^w: \Omega(Y)((u))\to\Omega(Y)((u))$ given by
\begin{equation}\label{nablaw}
\nabla^w=\frac{d}{du}+\frac{w}{u^2}+ \frac{\Gamma}{u}, \quad \Gamma\in\End(\Omega(Y)), \quad \Gamma|_{\Omega^{p}(Y)}=-\frac{p}{2}
\end{equation}
satisfies (cf. (\ref{1}))
\begin{equation}\label{2}
[\nabla^w,-dw +ud]=\frac{1}{2u}(-dw +ud)
\end{equation}
thereby giving rise to a connection on the cohomology (cf. \cite[Lemma 3.10]{KKP}).

The latter example is closely related to the subject of the Gauss-Manin systems and the Brieskorn lattices (see section \ref{applicat}) and through that to the study of generalized Hodge structures and Frobenius manifolds associated with isolated singularities (see \cite{HSe,ST} for review and references).

Our last example is algebraic. Let $A$ be a differential $\ZZ$-graded algebra. The previous two examples suggest that we should replace the complexes
$$
(\cA(X)[[u]], \bar{\partial}+u{\partial}), \quad (\Omega(Y)[[u]], -dw+ud)
$$
with the complex
$$
(\rC(A)[[u]], b+uB)
$$
(here $(\rC(A), b)$ is the Hochschild chain complex of $A$ and $B$ is the Connes differential) and try to repeat the construction. This does not quite work since, in general, the cohomology of the above complex has lots of torsion $\C[[u]]$-submodules meaning it is not a vector bundle. There are two possible options at this point: either we work with those dg algebras for which the torsion submodules do not occur (conjecturally, this is so for any homologically smooth and proper dg algebra \cite{KS}), or we replace the cohomology of $(\rC(A)[[u]], b+uB)$ by its image in the cohomology of $(\rC(A)((u)), b+uB)$, the periodic cyclic homology of $A$. Choosing the first option is analogous to working with compact K\"{a}hler manifolds only. Then, as we discussed previously, there is a chance that the resulting bundles with connections can be promoted to full-fledged non-commutative Hodge structures. We will keep considering arbitrary dg algebras or, more precisely, dg algebras with finite-dimensional periodic cyclic homology; the latter condition guarantees that we are still dealing with {\it finite rank} vector bundles on the formal disk.

The periodic cyclic homology, viewed as a bundle on the punctured formal disk, should carry a connection by a general argument (\cite[Section 11.5]{KS},\cite[Section 2.2.5]{KKP}) involving the so-called non-commutative Gauss-Manin connection \cite{Get} (along the parameter of a one-parameter deformation of $A$). In order to get a connection at the level of the periodic cyclic complex one should, perhaps, repeat the same argument but for a refined version of the non-commutative Gauss-Manin connection obtained in \cite{Ts} (see also \cite{DTT}). We hope to return to this idea on another occasion.
At the moment, borrowing an idea from \cite[Section 11.5]{KS}, we can write out ${\it a}$ connection $\nabla^A$ on the periodic cyclic complex of $A$ satisfying the same property as in our geometric examples above:
\begin{equation}\label{3}
[{\nabla^A}, b+uB]=\frac{1}{2u}(b+uB)
\end{equation}
An explicit formula for $\nabla^A$ can be found in section \ref{canuconn} (Proposition \ref{conncycl}). The origin of the formula is easy to explain. The reader familiar with the non-commutative differential calculus will notice that $\nabla^A$ is built by combining a grading operator on $\rC(A)$ with the non-commutative Cartan homotopy operator (cf. \cite[Section 2]{Get}) corresponding to the differential $d$ of $A$. Then the relation (\ref{3}) for this connection is a slight variation on the non-commutative Cartan homotopy formula \cite[(2.1)]{Get} (or, rather, a very special case of the latter which goes back to \cite{Ri}).

Since $\nabla^A$ is written in terms of some basic operations on the Hochschild complex, it can be easily generalized to define a connection $\nabla^{\mathscr{D}}$ on the cyclic complex of an arbitrary differential $\ZZ$-graded {\it category} $\mathscr{D}$.

Notice that when the $\ZZ$-grading on $A$ can be lifted to a $\Z$-grading, there is another obvious connection on the cyclic complex satisfying (\ref{3}), namely
$$
\nabla^A_{gr}=\frac{d}{du}+\frac{\Gamma'}{2u}
$$
where $\Gamma'$ is the corresponding $\Z$-grading operator on $\rC(A)$. It turns out that ${\nabla^A}$ and $\nabla^A_{gr}$ coincide at the cohomological level (see\footnote{If we knew that $\nabla^A$ was ${\it the}$ connection that the aforementioned argument produced then this statement would also follow from the argument.} section \ref{zgraded}). This observation shows that the de Rham part of non-commutative Hodge structures on the periodic cyclic homology is not very interesting in the case of $\Z$-graded categories (in other words, the connections are not essential for developing the Hodge theory in the $\Z$-graded setting; for more on this, see section \ref{NCHodge}). On the other hand, this same observation implies that our first geometric example is a special case of the algebraic one, at least, when $X$ is a smooth projective variety. Let us present the argument in a very sketchy manner since this is known, and is not the main point of the paper anyway.

The claim is that there is an isomorphism of bundles on the formal disk
$$
H^*(\cA(X)[[u]], \bar{\partial}+u{\partial})\cong H^*(\rC(\mathrm{par}_{\mathrm{dg}}X)[[u]], b+uB)
$$
(here $\mathrm{par}_{\mathrm{dg}}X$ stands for the dg category of perfect complexes\footnote{See \cite[Section 5.3]{Kel} and references therein for a discussion of the cyclic homology of this dg category.} on $X$) which induces an isomorphism of connections\footnote{Let us use the same notation for connections at the level of complexes and at the level of cohomology.}
$$
(H^*(\cA(X)((u)), \bar{\partial}+u{\partial}), \nabla^X)\cong (H^*(\rC(\mathrm{par}_{\mathrm{dg}}X)((u)), b+uB), \nabla^{\mathrm{par}_{\mathrm{dg}}X})
$$
To see this, note first that the category $\mathrm{par}_{\mathrm{dg}}X$ is $\Z$-graded and therefore we can replace $\nabla^{\mathrm{par}_{\mathrm{dg}}X}$ with $\nabla_{gr}^{\mathrm{par}_{\mathrm{dg}}X}$. Furthermore,
$$
H^*(\cA(X)[[u]], \bar{\partial}+u{\partial})\cong H^*(\cA(X), \bar{\partial})[[u]]
$$
$$
H^*(\rC(\mathrm{par}_{\mathrm{dg}}X)[[u]], b+uB)\cong H^*(\rC(\mathrm{par}_{\mathrm{dg}}X), b)[[u]]
$$
since the left-hand sides are known to be free $\C[[u]]$-modules. The right-hand sides can be endowed with meromorphic connections given by the same formulas that define $\nabla^X$ and $\nabla_{gr}^{\mathrm{par}_{\mathrm{dg}}X}$ and the isomorphisms above can be easily chosen so that they preserve the connections. To conclude the argument, we observe that the two $\Gamma'$s - in the definitions of $\nabla^X$ and $\nabla_{gr}^{\mathrm{par}_{\mathrm{dg}}X}$, respectively - match under a well-known isomorphism
$$
H^*(\cA(X),\bar{\partial})\cong H^*(\rC(\mathrm{par}_{\mathrm{dg}}X), b)
$$

The goal of the present paper is to explain why our {\it second} geometric example is a special case of the algebraic one. Namely, we will prove

\begin{theorem}\label{intro2} There is an isomorphism of bundles on the formal disk
$$
H^*(\Omega(Y)[[u]], -dw +ud)\cong H^*(\rC(\mathrm{MF}(w))[[u]], b+uB)
$$
which induces an isomorphism of connections
$$
\left(H^*(\Omega(Y)((u)), -dw +ud), \nabla^w\right)\cong \left(H^*(\rC(\mathrm{MF}(w))((u)), b+uB), \nabla^{\mathrm{MF}(w)}\right)
$$
Here $\mathrm{MF}(w)$ stands for the dg category of matrix factorizations of $w$.
\end{theorem}

Let us omit the definition of $\mathrm{MF}(w)$ since we will not need it. We will use as a black box a result of \cite{Dyck} which says that this dg category is (quasi-)equivalent to the category of dg modules over a certain differential $\ZZ$-graded algebra, $A_w$. This result allows us to replace the cyclic complex of the category $\mathrm{MF}(w)$ with the much smaller cyclic complex of $A_w$. The dg algebra is actually quite simple, so let us reproduce its definition here.

As a $\ZZ$-graded algebra, it is the tensor product $\C[Y]\otimes\mathrm{End}_\C V$ where $V$ is the $\ZZ$-graded space of polynomials in $k$ odd variables $\theta_1,\ldots,\theta_k$. The only part of $A_w$ that depends on $w$ is the differential. The differential depends, in fact, on more than just $w$ itself \footnote{Thus, our notation $A_w$ is a bit misleading.}: one also needs to pick a decomposition of $w$ of the form
$$
w=y_1w_1+\ldots+y_kw_k
$$
Then the differential on $A_w$ is simply the commutator with $D(w):=\sum_i\left(y_i\frac{\partial}{\partial \theta_i}+w_i\theta_i\right)$. Clearly, $D(w)^2=w$ and so $d_w=[D(w),-]$ squares to 0.

The second key ingredient used in the proof of Theorem \ref{intro2} is an explicit quasi-isomorphism from the Hochschild complex of $\mathrm{MF}(w)$ to the complex $(\Omega(Y), -dw)$ constructed recently in \cite{Seg} by combining results of {\it loc.cit.} with some facts obtained earlier in \cite{CalTu,Dyck}. We will need only the composition of this quasi-isomorphism with the embedding of the Hochschild complexes
$$
(\rC(A_w),b)\hookrightarrow (\rC(\mathrm{MF}(w)), b)
$$
The composition is still a quasi-isomorphism, by the aforementioned result of \cite{Dyck}.

It is worthwhile to emphasize that the results from \cite{Dyck,Seg}
we need (and, as a consequence, Theorem \ref{intro2} itself) hold
true if we replace the pair $(\C[Y],w)$ with more general pairs
$(R,w)$ having isolated critical loci (cf. definitions in
\cite[Section 3]{Dyck} and a discussion and citations at the end of
\cite[Section 3]{Seg}). For example, we could have started with the
pair $(\C[Y]_{(0)},w)$ where $Y=\C^k$, $\C[Y]_{(0)}$ is the local
ring at the origin, and $w$ satisfies (\ref{crit}) and has an
isolated singularity at the origin. For the sake of simplicity we
will stick to the polynomial algebra and will not attempt to present
the results in their ``most natural'' generality.

To conclude, let us point out that Theorem \ref{intro2} has some implications for the classical singularity theory (not related to non-commutative Hodge structures). These will be discussed briefly in section \ref{applicat}.

\bigskip
{\bf Organization of the paper.} The main body of the paper comprises four sections.

In section \ref{mixcompl} we develop some simple convenient formalism of mixed complexes with connections. It can be viewed as a step toward defining the {\it derived category} of mixed complexes with connections but we do not pursue this idea due to lack of motivation.

In section \ref{cyclcompl} we recall the definition of the mixed cyclic complexes of dg algebras, both unital and non-unital, and then move on to describing connections on these complexes. We present several equivalent sets of formulas for the connections and discuss separately the case of differential $\Z$-graded algebras. The section is concluded with an explanation of why the language of connections is suitable for developing the Hodge theory in the categorical framework.

In section \ref{locsing} we explain the proof of Theorem \ref{intro2}.

Section \ref{applicat} is devoted to some applications of Theorem \ref{intro2}.

In order to keep the main body of the text as short as possible, we have collected all the proofs in several appendices which occupy most of the paper.

\bigskip
{\bf Conventions.} All complexes and dg algebras in this paper are $\ZZ$-graded, unless we explicitly say otherwise. That is, our complexes are pairs $(\rC,d)$, where $\rC$ is a $\ZZ$-graded vector space $\rC=\rC_\mathrm{even}\oplus\rC_\mathrm{odd}$ and $d$ is an odd operator such that $d^2=0$, and our dg algebras are such complexes equipped with a compatible multiplication.

In sections \ref{mc}--\ref{zgraded} we are working over a ground field $\bK$ whose characteristic is not equal to 2. Starting from section \ref{NCHodge} the ground field is $\C$. We will be also considering complexes over the field of formal Laurent series but it will be clear from the context (or an explicit comment) which of the fields is important at a given moment.

\medskip
{\bf Acknowledgements.} I would like to thank Y.~Soibelman for explaining some basics of the categorical Hodge theory to me several years ago. I am also grateful to M.~Herbst, C.~Hertling, D.~Murfet, A.~Polishchuk, H.~Ruddat, E.~Scheidegger, C.~Sevenheck, and K.~Wendland for inspiring discussions on matrix factorizations and Hodge theory. Special thanks are due to C.~Sevenheck for a number of interesting comments on the subject matter of this paper, to the anonymous referee for very valuable suggestions on how the exposition could be improved, and to both of them for spotting several typos.

\bigskip

\section{Formalism of mixed complexes and $u$-connections}\label{mixcompl}

\subsection{Mixed complexes}\label{mc}

A {\it mixed complex} is a triple $(\rC, b, B)$ where $(\rC, b)$ is a ($\bK$-linear) complex and $B$ is an odd operator on $\rC$ that anti-commutes with $b$ and squares to 0. In the $\ZZ$ setting, the mixed complexes are merely double complexes but the use of the term ``mixed'' suggests the analogy with the mixed complexes in the conventional $\Z$-graded setting \cite{Lod}, and also emphasizes the fact that $b$ and $B$ play different roles.

Recall that we are using $u$ to denote a formal (even) variable, playing the role of a coordinate on the formal disk,  and $V((u))$ to denote the $\bK((u))$-vector space of formal Laurent series with coefficients in $V$ where $V$ is any $\bK$-vector space. Given a mixed complex $(\rC, b, B)$, we will call the $\bK((u))$-linear complex
$$
(\rC((u)), b+uB)
$$
the {\it $u$-totalization} of $(\rC, b, B)$.

As it should be clear from the Introduction, the mixed complexes and their $u$-totalizations are used as means to construct vector bundles over the formal disk. With this in mind, we will restrict ourselves to mixed complexes of {\it finite type}, i.e. those complexes whose $u$-totalizations have finite-dimensional cohomology over $\bK((u))$.

A morphism of mixed complexes is defined as a morphism of the underlying ordinary complexes that commutes with the $B$-operators. We will use also the following weaker notion \footnote{This is a version of the $S$-morphisms \cite{Lod}.}:  a {\it $u$-morphism} of mixed complexes is a $\bK((u))$-linear morphism of the corresponding $u$-totalizations (as complexes). Any $u$-morphism $f(u): (\rC, b, B)\to (\rC', b', B')$ can be written in the form
\begin{eqnarray*}
f(u)=\sum_{i=k}^\infty u^i f_i
\end{eqnarray*}
for some $k\in\Z$ where the coefficients $f_i: \rC\to\rC'$ are even operators.

\medskip
The reader, we hope, understands how to extend this terminology to other notions used in the study of complexes, such as quasi-isomorphisms, homotopy equivalences, etc. For instance, two $u$-morphisms $f(u)$ and $g(u)$, with the same domain and range, will be called {\it $u$-homotopic} if there is an odd $\bK((u))$-linear operator $H(u): \rC((u))\to\rC'((u))$ (a {\it $u$-homotopy}) such that
$$
f(u)-g(u)=H(u)(b+uB)+(b'+uB')H(u)
$$
Again, each $u$-homotopy has the form
\begin{eqnarray*}
H(u)=\sum_{i=k}^\infty u^i H_i
\end{eqnarray*}
where $k\in\Z$ and $H_i: \rC\to\rC'$ are odd operators. Accordingly, two mixed complexes $(\rC, b, B)$ and $(\rC', b', B')$ will be called {\it $u$-homotopy equivalent} if there are $u$-morphisms
\begin{eqnarray}\label{u-hom}
\iota(u): (\rC, b, B)\rightleftarrows (\rC', b', B') :p(u)
\end{eqnarray}
such that $\iota(u)p(u)$ and $p(u)\iota(u)$ are $u$-homotopic to $id_{\rC'}$ and $id_{\rC}$, respectively.

It is a standard fact that for complexes of vector spaces there is not much difference between homotopy equivalences and quasi-isomorphisms. This fact applies to our setting since we are working with complexes of vector spaces over $\bK((u))$:
\begin{proposition}\label{easylemma}
Any $u$-quasi-isomorphism $f(u): (\rC, b, B)\to (\rC', b', B')$ can be completed to a $u$-homotopy equivalence $f(u): (\rC, b, B)\rightleftarrows(\rC', b', B'): g(u)$.
\end{proposition}

Given a mixed complex $(\rC, b, B)$, the image of the cohomology of $(\rC[[u]], b+uB)$ in that of $(\rC((u)), b+uB)$ under the canonical ($\bK[[u]]$-linear) embedding of the complexes will be called the {\it canonical $u$-lattice}. The canonical $u$-lattice plays the role of the bundle on the formal disk associated with the mixed complex. We will say that a $u$-morphism of mixed complexes is {\it regular} if the image of the canonical $u$-lattice under the induced operator on the cohomology of the $u$-totalizations belongs to the canonical $u$-lattice. In particular, a $u$-homotopy equivalence will be called regular if both $u$-morphisms involved are regular.

\medskip

\subsection{$u$-Connections}\label{uconnections}

The following definition has been inspired by the examples we discussed in the Introduction (see (\ref{1}), (\ref{2}), (\ref{3})):
A {\it $u$-connection} on a mixed complex $(\rC, b, B)$ is a $\bK$-linear operator
$$
\nabla: \rC((u))\to \rC((u))
$$ satisfying the properties:
$$[\nabla,u]=id_{\rC}$$ (meaning $\nabla=\frac{d}{du}+\cA(u)$
where $\cA(u)$ is a $\bK((u))$-linear operator from $\rC((u))$ to itself)
and
\begin{eqnarray}\label{u-conn}
[\nabla,b+uB]=\frac{1}{2u}(b+uB)\quad (\,\,\Leftrightarrow\,\,[\cA(u),b+uB]=\frac{1}{2u}(b-uB)\,\,)
\end{eqnarray}

The equality (\ref{u-conn}) implies that $\nabla$ descends to the cohomology of the $u$-totalization, thereby giving rise to a connection on the corresponding bundle on the punctured formal disk.

\medskip

\begin{remark}{\rm Notice that since $\nabla$ does not commute with the differential, a $u$-connection on a mixed complex is not a $\ZZ$-graded complex of $D$-modules on the punctured formal disk, in the conventional sense. However, with any mixed complex with a $u$-connection one can associate a $\Z$-graded (unbounded) complex of $D$-modules as follows. Let $(\rC, b, B, \nabla)$ be as above and $D^n=D^n(\rC, b, B, \nabla)$ be the sequence of $D$-modules given by
\begin{eqnarray}\label{compofdmod}
D^n=
\begin{cases}
\left(\rC_\mathrm{even}((u)), \nabla|_{\rC_\mathrm{even}((u))}\right)\langle \frac{n}2\rangle, & n\,\,\,\text{even}\\
\left(\rC_\mathrm{odd}((u)), \nabla|_{\rC_\mathrm{odd}((u))}\right)\langle \frac{n}2\rangle, & n\,\,\,\text{odd}
\end{cases}
\end{eqnarray}
where $(V((u)), \nabla)\langle \frac{n}2\rangle:=(V((u)), \nabla-\frac{n}2\cdot\frac1{u})$, a Tate twist of $(V((u)), \nabla)$ (cf. \cite[Section 2.1.7]{KKP}). Then by (\ref{u-conn}) the differential $b+uB: D^n\to D^{n+1}$ is a morphism of $D$-modules for all $n$.
We will return to this construction in section \ref{NCHodge}; apart from that section, it will not be used in what follows. }
\end{remark}

\medskip

To define morphisms of mixed complexes with $u$-connections, we need
\medskip
\begin{lemma}\label{diffconn} Let $(\rC, b, B, \nabla=\frac{d}{du}+\cA(u))$ and $(\rC', b', B',\nabla'=\frac{d}{du}+\cA'(u))$ be mixed complexes with $u$-connections. Then for any $u$-morphism $f(u)$ of the mixed complexes
\begin{eqnarray}\label{ucum}
\frac{df(u)}{du}+\cA'(u)f(u)-f(u)\cA(u)
\end{eqnarray}
is a $u$-morphism.
\end{lemma}
\noindent{\bf Proof} is given in section \ref{proofdiffconn}. \hfill $\Box$

\medskip
We will say that  $f(u)$ defines a {\it morphism of mixed complexes with $u$-connections}, or simply a {\it morphism of $u$-connections}, if (\ref{ucum}) is $u$-homotopic to 0.

\medskip

As the next claim shows, being a morphism of $u$-connections is stable under passing to $u$-homotopic $u$-morphisms.

\begin{proposition}\label{hommaps}
Let $(\rC, b, B, \nabla)$, $(\rC', b', B',\nabla')$ be mixed complexes with $u$-connections and $f(u), g(u): (\rC, b, B)\to(\rC', b', B')$ two $u$-morphisms that are $u$-homotopic to each other. If $f(u)$ is a morphism of $u$-connections then so is $g(u)$.
\end{proposition}
\noindent{\bf Proof} is given in section \ref{proofhommaps}. \hfill $\Box$

\medskip

The following statement follows immediately from the definitions:
\begin{proposition}\label{composition}
The composition of morphisms of $u$-connections is a morphism of $u$-connections.
\end{proposition}

\medskip

We will say that two mixed complexes with $u$-connections $(\rC, b, B, \nabla)$ and $(\rC', b', B',\nabla')$ (or simply two $u$-connections $\nabla$ and $\nabla'$) are {\it homotopy gauge equivalent} if there is a $u$-homotopy equivalence (\ref{u-hom}) such that $p(u)$ and $\iota(u)$ are morphisms of $u$-connections.

\medskip
By Proposition \ref{composition}, homotopy gauge equivalence is an equivalence relation.

\medskip

\begin{proposition}\label{qishom}
If $f(u): (\rC, b, B, \nabla)\to(\rC', b', B',\nabla')$ is a $u$-quasi-isomorphism and a morphism of $u$-connections then
$\nabla$ and $\nabla'$ are homotopy gauge equivalent.
\end{proposition}
\noindent{\bf Proof} is given in section \ref{proofqishom}. \hfill $\Box$

\medskip

Homotopy gauge equivalence can be refined as follows: We may require it to be furnished by regular $u$-homotopy equivalences, in the sense of the previous section. This is similar to the notion of holomorphic equivalence in the classical theory of connections.
In all the examples we will consider in the remaining sections, the $u$-homotopy equivalences will be regular, and so in all the statements we will formulate homotopy gauge equivalence can be understood in the above stronger sense.

A very special instance of homotopy gauge equivalence is the case when one $u$-connection is obtained from another one, living on the same mixed complex, by adding a $u$-endomorphism that is $u$-homotopic to 0 (in this case the two $u$-connections are homotopy gauge equivalent with respect to the identity endomorphism of the mixed complex). Since there is not any actual ``gauging'' in this cases, we will simply say that two such $u$-connections are {\it equal up to a $u$-homotopy}.

\medskip
We will conclude this section with the following result which allows one to transfer $u$-connections using $u$-homotopy equivalences.

\begin{proposition}\label{gauge}
Given a $u$-homotopy equivalence (\ref{u-hom}) and a $u$-connection $\nabla'=\frac{d}{du}+\cA'(u)$ on $(\rC', b', B')$, there exists a $u$-connection $\nabla$ on $(\rC, b, B)$ that is homotopy gauge equivalent to $\nabla'$ with respect to $p(u)$ and $\iota(u)$. Explicitly, one can set
\begin{eqnarray}\label{gaugeconn}
\nabla=\frac{d}{du}+\cA(u),\quad \cA(u):=p(u)\frac{d\iota(u)}{du}+p(u)\cA'(u)\iota(u)+\frac{1}{2u}H(u)(b-uB)
\end{eqnarray}
where $H(u)$ is any $u$-homotopy such that $p(u)\iota(u)=id_{\rC} + (b+uB)H(u)+H(u)(b+uB)$.
\end{proposition}
\noindent{\bf Proof} is given in section \ref{gaugeequiv}.
\hfill $\Box$

\medskip

\section{$u$-Connections on the cyclic complexes of dg algebras}\label{cyclcompl}

\subsection{The cyclic complexes of unital dg algebras}\label{hoch}

Let $A$ be a (not necessarily unital) dg algebra and $sA$ stand for $A$ with the reversed $\ZZ$-grading. Given $a\in A$, the corresponding element of
$sA$ will be denoted by $sa$. The parity of $a$ will be denoted by $|a|$; thus, $|sa|=|a|-1$. Let
\begin{eqnarray*}
\rC (A)=A\otimes
T(sA)=\bigoplus\limits_{n=0}^\infty A\otimes sA^{\otimes n}
\end{eqnarray*}
equipped with the induced $\ZZ$-grading. We will write the elements of
$A\otimes sA^{\otimes n}$ as $$a_0[a_1|a_2|\ldots |a_n]$$ (i.e. $a_0[a_1|a_2|\ldots
|a_n]=a_0\otimes sa_1\otimes sa_2\otimes\ldots\otimes sa_n$), or simply $a_0$, if $n=0$.

\medskip

Throughout the paper, we will use the following convention: for an operator $T$ from $\rC (A)$ to anywhere we will write $T_{n+1}$ for the restriction of $T$ onto $A\otimes sA^{\otimes n}$.

\medskip

$\rC (A)$ is the underlying $\ZZ$-graded space of the Hochschild chain complex of $A$. Let us recall the definition of the Hochschild differential,  $b$.

Let $\tau$ denote the cyclic permutation on $\rC (A)$:
\begin{eqnarray*}
\tau_{n+1}(a_0[a_1|a_2|\ldots |a_n])=(-1)^{|sa_0|\sum_{i=1}^n|sa_i|}a_1[a_2|\ldots |a_n|a_0]
\end{eqnarray*}
One can easily see that $\tau_{n+1}^{n+1}=1$. Set
\begin{eqnarray*}
\ld^{(0)}_{n+1}(a_0[a_1|a_2|\ldots |a_n])=da_0[a_1|a_2|\ldots |a_n]
\end{eqnarray*}
\begin{eqnarray*}\lm^{(0)}_{n+1}(a_0[a_1|a_2|\ldots |a_n])=
\begin{cases}
0 & n=0\\
(-1)^{|a_0|}a_0a_1[a_2|\ldots |a_n] & n\geq1\\
\end{cases}
\end{eqnarray*}
and
\begin{eqnarray*}
\ld^{(i)}_{n+1}:=\tau_{n+1}^{-i}\,\ld^{(0)}_{n+1}\, \tau_{n+1}^{i}, \quad \lm^{(i)}_{n+1}:=\tau_{n}^{-i}\,\lm^{(0)}_{n+1}\, \tau_{n+1}^{i},\quad i=1,
\ldots,n
\end{eqnarray*}
Explicitly,
\begin{eqnarray}\label{Di}
\ld^{(i)}_{n+1}(a_0[a_1|a_2|\ldots |a_n])=(-1)^{\sum_{k=0}^{i-1}|sa_k|}a_0[a_1|\ldots|da_i|\ldots |a_n]
\end{eqnarray}
\begin{eqnarray}\label{Mi}
\lm^{(i)}_{n+1}(a_0[a_1|a_2|\ldots |a_n])=\begin{cases} (-1)^{\sum_{k=0}^{i}|sa_k|+1}a_0[a_1|\ldots|a_ia_{i+1}|\ldots |a_n] & i<n \\
-(-1)^{|sa_n|\left(|a_0|+\sum_{k=1}^{n-1}|sa_k|\right)}a_na_0[a_1|\ldots |a_{n-1}] & i=n
\end{cases}
\end{eqnarray}
Then the Hochschild differential $b$ is defined as follows\footnote{This definition of the Hochschild differential is easily seen to be equivalent to the one given in \cite[Section 1]{Get}.}
\begin{eqnarray}\label{newb0b1}
b=\bd+\bm, \quad \bd_{n+1}=\sum_{i=0}^{n}\ld^{(i)}_{n+1}, \quad
\bm_{n+1}=\sum_{i=0}^{n}\lm^{(i)}_{n+1}
\end{eqnarray}
That $b$ squares to 0 can be deduced form the formulas we write out in section \ref{auxlemmas}. The cohomology of $(\rC (A), b)$ is called the Hochschild homology of $A$.

\medskip

Let us recall now the definition of the mixed complex that computes the {\it cyclic homology} of (unital) dg algebras. The underlying ordinary complex is $(\rC (A), b)$, and $B$ is defined as follows.

Assume $A$ is unital. Set
\begin{eqnarray*}
N_{n+1}=\sum_{i=0}^{n}\tau_{n+1}^{i},\quad \h_{n+1}(a_0[a_1|a_2|\ldots |a_n])=1[a_0|a_1|a_2|\ldots |a_n]
\end{eqnarray*}
Then
\begin{eqnarray*}
B_{n+1}=(1-\tau^{-1}_{n+2})\h_{n+1}N_{n+1}
\end{eqnarray*}
Clearly, $B$ is an odd operator that squares to 0. Again, the formulas from section \ref{auxlemmas} imply that $B$ anti-commutes with both $\bd$ and $\bm$. The cohomology of the $u$-totalization of this mixed complex is called the periodic cyclic homology of $A$.

\medskip

\subsection{The cyclic complexes of non-unital dg algebras}\label{hoch1}

Notice that when $A$ is non-unital, the definition of $B$ from the previous section does not make sense. To fix the problem, one replaces the Hochschild chain complex with a quasi-isomorphic one, which contains a replacement for the unit. The aim of the present section is to recall the definition of this new complex.

Let $A$ be a dg algebra, with or without unit. Denote by $A^+$ the dg algebra obtained from $A$ by adjoining a unit:
\begin{eqnarray*}
A^e=A\oplus\bK e, \quad |e|=0\\
de=0,\quad ae=ea=a \quad \forall a\in A
\end{eqnarray*}
Consider the following $\ZZ$-graded vector space:
\begin{eqnarray*}
\rC^e (A)=A\oplus\bigoplus\limits_{n=1}^\infty (A^e\otimes sA^{\otimes n})
\end{eqnarray*}
The operators $\ld^{(i)}_{n+1}$ and $\lm^{(i)}_{n+1}$ on $\rC (A)$ extend to $\rC^e (A)$ via the formulas (\ref{Di}) and (\ref{Mi}), and we get a  complex $(\rC^e (A), b^e=\bed+\bem)$ where the differentials $\bed$, $\bem$ are defined as before. It is nothing but the so-called normalization of the usual Hochschild complex of $A^e$ \cite{Lod}.

\medskip

The vector space $\rC^e (A)$ admits the following decomposition:
\begin{eqnarray}\label{decom}
\rC^e (A)=\rC (A)\oplus\rC^+ (A)
\end{eqnarray}
where
\begin{eqnarray*}
\rC^+ (A)=\bigoplus\limits_{n=1}^\infty(\bK e\otimes sA^{\otimes n})
\end{eqnarray*}
Note, however, that $\rC^+ (A)$ is not stable under $\bem$, so it is not a subcomplex.

With (\ref{decom}) in mind, we can represent the elements of $\rC^e (A)$ by column-vectors with two components, the first component being an element of $\rC (A)$ and the second one an element of $\rC^+ (A)$. Accordingly, operators on $\rC^e (A)$ can be represented by $2\times2$ matrices. For example, the differential $\bed$ preserves the decomposition (\ref{decom}) and therefore can be represented by the matrix
$$
\begin{pmatrix}\bd & 0 \\0 & \bd \end{pmatrix}
$$
The operator $\lm^{(0)}$, on the other hand, preserves $\rC (A)$ but maps $\rC^+ (A)$ to $\rC (A)$; it is then represented by the matrix
$$
\begin{pmatrix}\lm^{(0)} & \lm^{(0)} \\0 & 0 \end{pmatrix}
$$
Note that, formally, the two copies of $\bd$ in the first matrix (or $\lm^{(0)}$ in the second matrix) represent {\it different} operators, since they have different domains/ranges.

In this notation, the differential $\bem$ is easily seen to correspond to the following upper-triangular matrix
\begin{eqnarray*}
\begin{pmatrix}\bm &  \bmv \\0 & \bmh \end{pmatrix},
\end{eqnarray*}
where\footnote{We extend our previous convention regarding subscripts as follows: For an operator $T$ from $\rC^e (A)$ to anywhere we will write $T_{n+1}$ for the restriction of $T$ onto $(A\otimes sA^{\otimes n})\oplus(\bK e\otimes sA^{\otimes n})$.}
\begin{eqnarray*}
\bmv_{n+1}=\lm^{(0)}_{n+1}+\lm^{(n)}_{n+1}: \bK e\otimes sA^{\otimes n}\to A\otimes sA^{\otimes (n-1)}
\end{eqnarray*}
\begin{eqnarray*}
\bmh_{n+1}=\sum_{i=1}^{n-1}\lm^{(i)}_{n+1}: \bK e\otimes sA^{\otimes n}\to \bK e\otimes sA^{\otimes (n-1)}
\end{eqnarray*}

\medskip

Let us describe the analog of $B$ for the new complex $(\rC^e (A), b^e)$. Define
$
h^e: \rC (A)\to \rC^+ (A)
$
by
\begin{eqnarray*}
h^e_{n+1}(a_0[a_1|a_2|\ldots |a_n])=e[a_0|a_1|a_2|\ldots |a_n]
\end{eqnarray*}
Then, in our matrix notation, $B^e: \rC^e (A)\to \rC^e (A)$ is given by
\begin{eqnarray*}
B^e_{n+1}=\begin{pmatrix}0 & 0 \\h^e_{n+1}N_{n+1} & 0 \end{pmatrix}
\end{eqnarray*}
$B^e$ is an odd operator that squares to 0 and anti-commutes with $b^e$:
$$\quad \bed B^e=-B^e\bed, \quad \bem B^e=-B^e\bem$$

\medskip

We have now two a priori different definitions of the cyclic homology for unital dg algebras, namely, via $(\rC(A), b, B)$ and $(\rC^e (A), b^e, B^e)$. Let us explain why the two definitions are equivalent.

Let $A$ be a unital dg algebra. Consider the following $u$-morphisms
\begin{eqnarray}\label{iota}
\iota(u): (\rC(A), b, B)\to (\rC^e (A), b^e, B^e), \quad \iota(u)=\begin{pmatrix}id_{\rC (A)} \\ 0  \end{pmatrix}+u\begin{pmatrix}0 \\ h^ehN  \end{pmatrix}
\end{eqnarray}
\begin{eqnarray}\label{p}
p(u): (\rC^e (A), b^e, B^e)\to (\rC(A), b, B), \quad p(u)=\begin{pmatrix}id_{\rC (A)} & (1-\tau^{-1})h\lm^{(0)}  \end{pmatrix}
\end{eqnarray}
(the operator $h\lm^{(0)}$ in the last line maps $e[a_0|a_1|a_2|\ldots |a_n]$ to $1[a_0|a_1|a_2|\ldots |a_n]$).
\medskip

\begin{proposition}\label{quasiisocycl}
$\iota(u)$ and $p(u)$ establish a $u$-homotopy equivalence between $(\rC(A), b, B)$ and $(\rC^e(A), b^e, B^e)$.
\end{proposition}
\noindent{\bf Proof} is given in section \ref{homotopy}.
\hfill $\Box$

\medskip

The last result in this section relates the cyclic complexes of a dg algebra $A$ and its {\it opposite} dg algebra $A^\circ$. Recall that $A^\circ$ is just $A$ with a new product:
$$
a'\otimes a''\mapsto(-1)^{|a'||a''|}a''a'
$$

Consider the following isomorphism of $\ZZ$-graded vector spaces:
\begin{eqnarray*}
\Phi=\Phi_A: \rC^e(A)\to\rC^e(A^\circ),\quad a_0[a_1|a_2|\ldots |a_n]\mapsto(-1)^{n+\sum_{1\leq i<j\leq
n}|sa_i||sa_j|}a_0[a_n|a_{n-1}|\ldots |a_1]
\end{eqnarray*}

\medskip

\begin{proposition}\label{opp} $\Phi$ is an isomorphism from the mixed complex $(\rC^e(A), b^e, B^e)$ to the mixed complex $(\rC^e(A^\circ), b^e, -B^e)$.
\end{proposition}
\noindent{\bf Proof} is based on the following properties of $\Phi$, which are easy to verify using the definitions of the maps involved:
\begin{eqnarray}\label{Phi}
\Phi\tau=\tau^{-1}\Phi,\quad \Phi\ld^{(0)}=\ld^{(0)}\Phi, \quad \Phi\lm^{(0)}=\lm^{(0)}\tau^{-1}\Phi, \quad \Phi h^e=-h^e\tau\Phi
\end{eqnarray}
It follows from these properties that
\begin{eqnarray*}\label{PhibB}
\Phi N=N\Phi,\quad \Phi_{n+1}\ld^{(i)}_{n+1}=\ld^{(n+1-i)}_{n+1}\Phi_{n+1}, \quad \Phi_{n}\lm^{(i)}_{n+1}=\lm^{(n-i)}_{n+1}\Phi_{n+1}
\end{eqnarray*}
which implies the statement.
\hfill $\Box$

\medskip

To conclude our discussion of the cyclic complexes, let us point out a variation on the definition of $(\rC(A), b, B)$ and $(\rC^e(A), b^e, B^e)$ which we will use in section \ref{locsing}. Namely, in the definition of $\rC(A)$ (or $\rC^e(A)$) one can replace the direct sums by the direct products:
\begin{eqnarray*}
\rrC(A)=\prod_n\left(A\otimes sA^{\otimes n}\right)_{even}\bigoplus \prod_n\left(A\otimes sA^{\otimes n}\right)_{odd}
\end{eqnarray*}
(and a similar version for $\rC^e(A)$ which we will denote by $\rrCe(A)$). Then all the above definitions (e.g. those of the differentials $b$, $B$) and conclusions (e.g. Proposition \ref{quasiisocycl}) are easily checked to hold true in the new setting. We refer the reader to \cite{PP} for a detailed study of these versions of the Hochschild complexes.

\medskip

\subsection{Canonical $u$-connections on the cyclic complexes}\label{canuconn}

Let us omit a detailed discussion on what the adjective ``canonical'' in the title refers to. In plain words, by a canonical $u$-connection one should understand a $u$-connection given by a {\it universal} formula, i.e. a formula that works simultaneously for all dg algebras. Such a $u$-connection can be extended to the cyclic complexes of dg categories thereby giving rise to a functor from the category of (small $\bK$-linear) dg categories to that of mixed complexes with connections which has various reasonable properties (i.e. Morita equivalent dg categories correspond to homotopically gauge equivalent $u$-connections, etc.).

Observe that any canonical $u$-connection has a ``twin'' obtained by means of the morphism $\Phi$ from the previous section (see Proposition \ref{opp}). Namely,
suppose
$$
A\rightsquigarrow \nabla^A=\frac{d}{du}+\cA_A(u)
$$
is a canonical $u$-connection. The following statement is an immediate consequence of Proposition \ref{opp}.

\medskip

\begin{proposition}
\begin{eqnarray}
A\rightsquigarrow (\nabla^{\circ})^A:=\frac{d}{du}-\Phi^{-1}_A\cdot \cA_{A^\circ}(-u)\cdot\Phi_A
\end{eqnarray}
is also a canonical $u$-connection.
\end{proposition}

Let us call this $u$-connection {\it dual} to the original one. The rest of this section is devoted to one example of a canonical $u$-connection. As we will show, it is self-dual (up to a $u$-homotopy).

\medskip

Let $A$ be a (not necessarily unital) dg algebra. Consider the following operators on $\rC^e(A)$:
$$
\Gamma_{n+1}=-\frac{n}2\cdot id_{n+1}
$$
\begin{eqnarray*}
\Udl_{n+1}=\frac12\begin{pmatrix}-\lm_{n+1}^{(0)}\ld_{n+1}^{(1)} & -\lm_{n+1}^{(0)}\ld_{n+1}^{(1)} \\0 & 0 \end{pmatrix}, \quad
\Vdl_{n+1}=\frac12\begin{pmatrix}0 & 0 \\ -h^e_{n+1}\sum_{i=1}^{n}\sum_{j=i+1}^{n+1}\tau_{n+1}^{j}\ld_{n+1}^{(i)} & 0 \end{pmatrix}
\end{eqnarray*}
\begin{eqnarray*}
\Ud_{n+1}=\frac12\begin{pmatrix}\lm_{n+1}^{(n)}\ld_{n+1}^{(n)} & \lm_{n+1}^{(n)}\ld_{n+1}^{(n)} \\0 & 0 \end{pmatrix}, \quad
\Vd_{n+1}=\frac12\begin{pmatrix}0 & 0 \\ h^e_{n+1}\sum_{i=1}^{n}\sum_{j=1}^{i}\tau_{n+1}^{j}\ld_{n+1}^{(i)} & 0 \end{pmatrix}
\end{eqnarray*}

\medskip

\begin{proposition}\label{conncycl}
\begin{eqnarray}\label{ext}
{\nabla}:=\frac{d}{du}+\frac{\Udl}{u^2}+\frac{\Vdl+\Gamma}{u}
\end{eqnarray}
is a $u$-connection on $(\rC^e(A), b^e, B^e)$; its dual is given by
\begin{eqnarray}\label{extR}
{\nabla^{\circ}}:=\frac{d}{du}+\frac{\Ud}{u^2}+\frac{\Vd+\Gamma}{u}
\end{eqnarray}
\end{proposition}
\noindent{\bf Proof} is given in section \ref{proofUV}. \hfill $\Box$

\medskip

\begin{proposition}\label{nabla-nabla} $\nabla$ and $\nabla^\circ$ are equal up to a $u$-homotopy.
\end{proposition}
\noindent{\bf Proof} is given in section \ref{proofnabla-nabla}.
 \hfill $\Box$

\medskip

Let us point out two corollaries of Proposition \ref{conncycl}. To formulate the first corollary, we need two more operators on $\rC^e(A)$:

\begin{eqnarray*}
\Uml_{n+1}=\frac12\begin{pmatrix}-\lm_{n}^{(0)}\lm_{n+1}^{(1)} & -\lm_{n}^{(0)}\lm_{n+1}^{(1)} \\0 & 0 \end{pmatrix},
\quad \Vml_{n+1}=\frac12\begin{pmatrix}id_{n+1} & 0 \\ -h^e_n\sum_{i=1}^{n-1}\sum_{j=i+1}^{n}\tau_{n}^{j}\lm_{n+1}^{(i)} & 0 \end{pmatrix}
\end{eqnarray*}

\medskip

\begin{corollary}\label{alterconncycl}
\begin{eqnarray}\label{alterext}
\widetilde{\nabla}=\frac{d}{du}+\frac{\Udl+\Uml}{u^2}+\frac{\Vdl+\Vml}{u}
\end{eqnarray}
is a $u$-connection on $(\rC^e(A), b^e, B^e)$ equal to ${\nabla}$ up to a $u$-homotopy.
\end{corollary}
\noindent{\bf Proof} is given in section \ref{proofalterconncycl}.
 \hfill $\Box$

\medskip

The second corollary concerns a canonical $u$-connection on $\rC(A)$ in the case when $A$ has unit. Consider the following operators on $\rC(A)$:
\begin{eqnarray*}
\nUd_{n+1}=-\frac12\lm_{n+1}^{(0)}\ld_{n+1}^{(1)}, \quad
\nVd_{n+1}=-\frac12(1-\tau^{-1}_{n+2})h_{n+1}\sum_{i=1}^{n}\sum_{j=i+1}^{n+1}\tau_{n+1}^{j}\ld_{n+1}^{(i)},
\end{eqnarray*}
and
\begin{eqnarray*}
 \nWd=\frac{1}{2}(1-\tau^{-1})hhhN\bd
\end{eqnarray*}

\medskip

\begin{corollary}\label{nonext}
\begin{eqnarray*}
\nabla^{un}:=\frac{d}{du}+\frac{\nUd}{u^2}+\frac{\nVd+\Gamma}{u}+\nWd
\end{eqnarray*}
is a $u$-connection on $(\rC(A), b, B)$, homotopy gauge equivalent to ${\nabla}$.
\end{corollary}
\noindent{\bf Proof} is given in section \ref{proofnonext}.
 \hfill $\Box$

\medskip

In section \ref{locsing}, we will apply the above results to the dg algebra $A_w$ defined in the Introduction. In that example, it will be important for us to work with the version of the cyclic complex defined at the end of section \ref{hoch1}, namely, ${\rrCe}(A)$. Let us therefore emphasize the obvious fact that the above $u$-connections extend to this complex.

\medskip

\subsection{On the case of differential $\Z$-graded algebras}\label{zgraded}
Let us assume that $A$ is $\Z$-graded or, more precisely, that its $\ZZ$-grading can be lifted to a $\Z$-grading.

Let us denote the $\Z$-degree of an elements $a\in A$ by $\deg(a)$. Clearly, the $\ZZ$-grading on the unitalization $A^e$ of $A$ also lifts to a $\Z$-grading, if we set $\deg(e)=0$. Therefore, the space $\rC^e(A)$ has a natural $\Z$-grading determined by the following operator
$$
\Gamma'(a_0[a_1|\ldots|a_n])=(\sum_{i=0}^n \deg(a_i)-n)a_0[a_1|\ldots|a_n]
$$
With respect to this grading $b^e$ and $B^e$ have degrees 1 and $-1$, respectively, i.e.
$$
[\Gamma', b^e]=b^e,\quad [\Gamma', B^e]=-B^e
$$
This shows that
$$
\nabla_{gr}=\frac{d}{du}+\frac{\Gamma'}{2u}
$$
is a $u$-connection.

\medskip

\begin{proposition}\label{zgr}
The $u$-connections $\nabla$ and $\nabla_{gr}$  are equal up to a $u$-homotopy.
\end{proposition}
\noindent{\bf Proof} is given in section \ref{proofzgr}.
 \hfill $\Box$

An important Hodge theoretic implication of this result will be explained in the next section.

\medskip

\subsection{Non-commutative Hodge filtrations}\label{NCHodge}
One of the goals of the present work is to contribute to the idea that the canonical $u$-connections form the correct framework for developing the Hodge theory of dg algebras and dg categories (at least, over the field of complex numbers). However, in the case of differential $\Z$-graded categories the Hodge theory can be formulated in an alternative way, with no mention of connections (see, for instance, \cite{Ka}). Proposition \ref{zgr} implies that the two approaches are equivalent, and the aim of this section is to explain this in more detail.

Let us start by recalling a description of the periodic cyclic homology of a differential $\Z$-graded category as a 2-periodic sequence of filtered vector spaces \cite[Sections 2,3]{Ka} where the filtrations are thought of as generalizing the classical Hodge filtration. The ground field in this section is $\C$.

Let $\cD$ be a differential $\Z$-graded category. Then the Hochschild complex of $\cD$ carries a $\Z$-grading\footnote{Let us not care about the various versions ($\rC,\rC^e$ etc.) of the Hochschild complex in this section.} (cf. section \ref{zgraded}), and one can consider a refined version of the $u$-totalization of $(\rC(\cD), b, B)$, namely,
$
(\rC_\bullet(\cD)((u))^{gr}, b+uB)
$
where $u$ is a formal variable {\it of degree 2} and ``${\rm gr}$'' means that we look at the subspace in $\rC_\bullet(\cD)((u))$ spanned by homogeneous series. The cohomology groups $\rHP_\bullet(\cD)$ of the resulting $\Z$-graded complex are known as the periodic cyclic homology groups of $\cD$. Clearly, all the groups $\rHP_{2k}(\cD)$ resp. $\rHP_{2k+1}(\cD)$ are isomorphic to each other.

In the example $\cD=\mathrm{par}_{\mathrm{dg}}X$ (see the Introduction) where $X$ is a smooth projective complex variety we have \cite{Ke,We}:
\begin{equation}\label{shkl1}
\rHP_{2k}(\cD)\simeq \bigoplus_p H^{2p}_{DR}(X),\quad \rHP_{2k+1}(\cD)\simeq \bigoplus_p H^{2p+1}_{DR}(X)
\end{equation}

For simplicity, we will assume from now on that the periodic cyclic homology groups of our categories are finite-dimensional.

The spaces $\rHP_{\bullet}(\cD)$ come equipped with canonical decreasing filtrations which generalize the  Hodge filtration. These filtrations are defined via an analog of the Hodge-to-de Rham spectral sequence  (cf. \cite{Ka}), in complete agreement with the classical picture. Namely, the complex
$
(\rC_\bullet(\cD)((u))^{gr}, b+uB)
$
is filtered by the sub-complexes
$
(u^n\rC_\bullet(\cD)[[u]]^{gr}, b+uB),
$
and the induced filtration on $\rHP_\bullet(\cD)$ is what we need. It follows from the definition that
\begin{equation}\label{2-period}
(\rHP_\bullet(\cD), F^\bullet\rHP_\bullet(\cD))\simeq(\rHP_{\bullet+2}(\cD), F^{\bullet+1}\rHP_{\bullet+2}(\cD))
\end{equation}
This is the aforementioned 2-periodicity property.

To see the connection to the classical case, note that the above filtration on the cyclic complex gives rise to a spectral sequence with
$
E_1=\rHH_\bullet(\cD)((u))^{gr}
$
where $\rHH_\bullet$ is the Hochschild homology. Given that
\[
\rHH_k(\cD)\simeq\bigoplus_{q-p=k} H^{q}(X, \Omega^{p}(X))
\]
when $\cD=\mathrm{par}_{\mathrm{dg}}X$ (cf. \cite{Ke,We}), this spectral sequence is a close relative of the classical Hodge-to-de Rham spectral sequence. In fact, in the above example, the non-commutative Hodge filtration on, say, ${\rHP}_{0}(\mathrm{par}_{\mathrm{dg}}X)$ transforms under the isomorphism (\ref{shkl1}) into
\[F^n\left(\oplus_p H^{2p}_{DR}(X)\right)=\oplus_{p-q\geq 2n}H^{p,q}(X)\]
Note that on each $H^{2p}_{DR}(X)$ it differs from the usual Hodge filtration by a shift.

Let us explain now how the canonical $u$-connection $\nabla^\cD$ gives rise, under certain conditions, to such a 2-periodic sequence of filtered vector spaces even when $\cD$ is not necessarily $\Z$-graded.

We will begin by recalling briefly some basic facts regarding {\it regular} singular connections in one variable (see \cite{S} and \cite[Section 2.2]{Sa} for a detailed treatment and a brief reminder, respectively). Given a finite-dimensional vector space $\widehat{\cG}$ over $\mathbb{C}((u))$ equipped with a regular singular connection $\nabla$ there is a canonical decreasing filtration of $\widehat{\cG}$ by free $\C[[u]]$-submodules ${V}^{\alpha}\widehat{\cG}$ determined by certain simple properties which we will not recall here (for instance, $u\nabla-\alpha$ is nilpotent on $\mathrm{Gr}_V^\alpha(\widehat{\cG}))$. In the most general case, it is indexed by $\C$ (equipped with a proper total order) but we will assume for simplicity that $\alpha\in\R$. Furthermore, there is a functor $\psi$ from the category of triples $(\widehat{\cG}, \nabla, \widehat{\cG}^{\,0})$, where $\widehat{\cG}$ and $\nabla$ are as above and  $\widehat{\cG}^{\,0}\subset\widehat{\cG}$ is a lattice, to the category of finite-dimensional filtered vector spaces (over $\mathbb{C}$) defined as follows:
\begin{eqnarray}\label{shkl5}
\psi(\widehat{\cG}, \nabla, \widehat{\cG}^{\,0})=(\mathsf{H}, F^\bullet \mathsf{H}),\nonumber\\
\mathsf{H}:=\bigoplus_{\alpha\in(-1,0]}\mathrm{Gr}_V^\alpha\widehat{\cG},\quad F^n\mathsf{H}:=\bigoplus_{\alpha\in(-1,0]}\mathrm{Gr}_V^\alpha(u^n\widehat{\cG}^{\,0})
\end{eqnarray}

Let $\cD$ be a dg category with the additional property that the connection $\nabla^\cD$ induces a regular singular connection on the cohomology of $(\rC(\cD)((u)), b+uB)$. Then, using (\ref{compofdmod}), we obtain a $\Z$-graded complex of $D$-modules $(D^\bullet(\rC(\cD), b, B, \nabla^\cD), b+uB)$ whose $n$th cohomology will be denoted by $(\widehat{\cG}_n(\cD), \nabla_n^\cD)$; that is,
\[
\widehat{\cG}_n(\cD)=H^{n\,\text{mod}\,2}(\rC(\cD)((u)), b+uB),\quad \nabla_n^\cD=\nabla^\cD-\frac{n}{2}\cdot\frac1{u}
\]
Let us also denote by $\widehat{\cG}^{\,0}_n(\cD)\subset\widehat{\cG}_n(\cD)$ the image of $H^{n\,\text{mod}\,2}(\rC(\cD)[[u]], b+uB)$ under the canonical map. Applying (\ref{shkl5}) to the resulting triples, we obtain a sequence of filtered vector spaces
\begin{eqnarray}\label{shkl7}
(\rHP_n(\cD), F^\bullet\rHP_n(\cD)):=\psi(\widehat{\cG}_n(\cD), \nabla_n^\cD,\widehat{\cG}^{\,0}_n(\cD))
\end{eqnarray}
One can check easily that these filtered vector spaces satisfy (\ref{2-period}).

When $\cD$ is $\Z$-graded, the regularity assumption is satisfied by Proposition \ref{zgr}. Moreover, it is easy to see that the above construction, when applied to $\nabla^\cD_{gr}$, reproduces the periodic cyclic homology groups with their non-commutative Hodge filtrations.

The previous discussion explains why the canonical $u$-connection should be viewed as a true generalization of the Hodge filtration.
Note, however, that in the $\Z$-graded setting the canonical $u$-connection does not give us any essential information beyond the filtered periodic cyclic homology. In this sense the $u$-connections are not needed for developing the Hodge theory in the $\Z$-graded case.

As we will see in section \ref{applicat}, where the example $\cD=\mathsf{MF}(w)$ is treated in some more detail, the canonical $u$-connections, in general, carry more information than the filtrations.

\medskip

\section{Mixed complexes and $u$-connections associated with polynomials}\label{locsing}
We remind the reader that from now on the ground field is $\C$. In this section, $Y$, $w$, $A_w$ etc. have the same meaning as in the Introduction, so we will not repeat their definitions. We will be also using all the other definitions and notation related to that part of the Introduction.

\medskip

\subsection{Mixed complexes associated with $w$}

As we mentioned in the Introduction, there is a quasi-isomorphism, constructed in \cite{Seg}, from the Hochschild complex of the algebra $A_w$ to the complex $(\Omega(Y), -dw)$. The purpose of this section is to show that the quasi-isomorphism can be promoted to a $u$-quasi-isomorphism of the corresponding mixed complexes. This requires working with the extended versions, $({\rC^e}(A_w), b^e, B^e)$ and $({\rrCe}(A_w), b^e, B^e)$, of the mixed cyclic complexes.
Apart from that, the only addition to the original construction of \cite{Seg} we need is the following well-known fact:

\medskip

\begin{lemma}\label{standard}
Let $f: (\rC,b,B)\to(\rC',b',B')$ be a morphism of mixed complexes that induces a quasi-isomorphism from $(\rC,b)$ to $(\rC',b')$. If $H^*(\rC,b)$ and $H^*(\rC',b')$ are finite-dimensional and are purely even or purely odd (i.e. sit in one degree only) then $f$ is a $u$-quasi-isomorphism.
\end{lemma}
\noindent{\bf Sketch of the Proof.} The conditions of the Lemma imply, by a simple inductive argument, that $H^*(\rC[[u]],b+uB)$ and $H^*(\rC'[[u]],b'+uB')$ are {\it free} $\bK[[u]]$-modules of finite rank. Obviously, $f$ lifts to a $\bK[[u]]$-linear morphism $(\rC[[u]],b+uB)\to(\rC'[[u]],b'+uB')$. The induced map $\bar{f}$ on the cohomology is a morphism of free $\bK[[u]]$-modules whose reduction modulo $u$ is an isomorphism, that is, $\bar{f}$ is itself an isomorphism. To conclude the proof, it remains to use the exactness of the functor $\bK((u))\otimes_{\bK[[u]]}-$.
 \hfill $\Box$

\medskip

To describe the $u$-morphism from $({\rC^e}(A_w), b^e, B^e)$ to $(\Omega(Y)\,,\, -dw\,,\, d)$ we need some additional ingredients.

First, let us introduce the following two operators on ${\rC^e}(A_w)$: the odd operator $\bw$ determined by
$$
\bw_{n+1}:=\sum_{i=1}^{n+1} w^{(i)}_{n+1},\quad  w^{(i)}_{n+1}(a_0[a_1|\ldots |a_n]):=(-1)^{\sum_{j=0}^{i-1}|sa_j|}a_0[a_1|\ldots|a_{i-1}|w|a_{i}| \ldots|a_n]
$$
and the even operator $\bdw$ determined by
$$
\bdw_{n+1}:=\sum_{i=1}^{n+1}D(w)^{(i)}_{n+1},\quad D(w)^{(i)}_{n+1}(a_0[a_1|\ldots |a_n]):=a_0[a_1|\ldots|a_{i-1}|D(w)|a_{i}| \ldots|a_n]
$$
Observe that the operators extend to ${\rrCe}(A_w)$. Also note that $\bw$ preserves the image of the embedding ${\rrCe}(\C[Y])\to\rrCe(A_w)$ induced by the obvious embedding $\C[Y]\to A_w$.

\medskip

The second ingredient we need is a map $\str: {\rrCe}(A_w)\to {\rrCe}(\C[Y])$ defined as follows\footnote{It is, of course, a super-analog of the standard trace map from \cite[(1.2.1)]{Lod}.}. Recall that $A_w=\C[Y]\otimes\DD$. Pick a basis $\{v_1,\ldots, v_{2^{n}}\}$ in $V$ consisting of even and odd vectors and denote by $E_{ij}$, $i,j=1,\ldots2^{n}$, the operators on $V$ given by the elementary $2^{n}\times 2^{n}$ matrices in the above basis. Then
\begin{eqnarray*}
\str: (\phi_0\otimes E_{i_0i_1})[\phi_1\otimes E_{i_1i_2}|\ldots |\phi_n\otimes E_{i_ni_{n+1}}]\mapsto \begin{cases}(-1)^\ast\phi_0[\phi_1|\ldots |\phi_n] & i_0=i_{n+1}\\ 0 & \text{otherwise}\end{cases}\\
\str: (e[\phi_1\otimes E_{i_1i_2}|\ldots |\phi_n\otimes E_{i_ni_{n+1}}]\mapsto \begin{cases}(-1)^{\ast\ast} e[\phi_1|\ldots |\phi_n] & i_1=i_{n+1}\\ 0 & \text{otherwise}\end{cases}
\end{eqnarray*}
where $\ast=(n-1)|v_{i_0}|+\sum_{s=1}^n|v_{i_s}|$ and $\ast\ast=n|v_{i_1}|+\sum_{s=2}^n|v_{i_s}|$.

\medskip

\begin{proposition}\label{newmix}$\,$\\

a) $({\rrCe}(A_w), \bem+\bw, B^e)$ and $({\rrCe}(\C[Y]), \bem+\bw, B^e)$ are mixed complexes.

b) The map $\mathrm{exp}(-\bdw):{\rC^e}(A_w)\to{\rrCe}(A_w)$ induces a morphism of mixed complexes
$$({\rC^e}(A_w), b^e, B^e)\to({\rrCe}(A_w), \bem+\bw, B^e)$$

c) $\str: ({\rrCe}(A_w), \bem+\bw, B^e)\to ({\rrCe}(\C[Y]), \bem+\bw, B^e)$ is a morphism of mixed complexes.

d) The map $\epsilon: {\rrCe}(\C[Y])\to \Omega(Y)$ determined by
\begin{eqnarray*}
\epsilon_{n+1}(\phi_0[\phi_1|\ldots |\phi_n])=\frac{1}{n!}\phi_0d\phi_1\wedge\ldots \wedge d\phi_n,\quad \epsilon_{n+1}(e[\phi_1|\ldots |\phi_n])=\frac{1}{n!}d\phi_1\wedge\ldots \wedge d\phi_n
\end{eqnarray*}
is a morphism of mixed complexes $({\rrCe}(\C[Y]), \bem+\bw, B^e)\to(\Omega(Y)\,,\, -dw\,,\, d)$.

e) The composition $\epsilon\cdot\str\cdot\mathrm{exp}(-\bdw): ({\rC^e}(A_w), b^e, B^e)\to (\Omega(Y)\,,\, -dw\,,\, d)$ is a $u$-quasi-isomorphism.
\end{proposition}
\noindent{\bf Proof} is given in section \ref{proofnewmix}.
 \hfill $\Box$

\medskip

\subsection{Comparing connections}
Consider the following operators on $\rC^e(A_w)$:
\begin{eqnarray*}
\Uwc_{n+1}=\frac12\begin{pmatrix}-\lm_{n+2}^{(0)} w^{(1)}_{n+1} & -\lm_{n+2}^{(0)} w^{(1)}_{n+1} \\0 & 0 \end{pmatrix},\quad \Vwc_{n+1}=\frac12\begin{pmatrix}0 & 0 \\ -h^e_{n+2}\sum_{i=1}^{n+1}\sum_{l=i+1}^{n+2}\tau_{n+2}^{l} w^{(i)}_{n+1} & 0 \end{pmatrix}
\end{eqnarray*}
Note that the operators extend to $\rrCe(A_w)$.

\medskip
\begin{proposition}\label{conntw} \begin{eqnarray*}
a)\quad \nabla^{\flat,w}:=\frac{d}{du}+\frac{2\Uwc}{u^2}+\frac{2\Vwc+\Gamma}{u}
\end{eqnarray*} is a $u$-connection\footnote{Including $\frac12$ in the definition of $\Uwc$ and $\Vwc$ and then multiplying them by 2 may appear somewhat strange; this is motivated by the result of Lemma \ref{0}.} on $({\rrCe}(A_w), \bem+\bw, B^e)$.

b) The map $\mathrm{exp}(-\bdw):{\rC^e}(A_w)\to{\rrCe}(A_w)$ is a morphism of mixed complexes with $u$-connections
$$({\rC^e}(A_w), b^e, B^e, {\nabla})\to({\rrCe}(A_w), \bem+\bw, B^e, {\nabla}^{\flat,w})$$
\end{proposition}
\noindent{\bf Proof} is given in section \ref{proofconntw}.
 \hfill $\Box$

\medskip

Clearly, the operators $\Uwc,\Vwc$ preserve the subspace ${\rrCe}(\C[Y])$ in $\rrCe(A_w)$. Thus, we have a $u$-connection on $({\rrCe}(\C[Y]), \bem+\bw, B^e)$
which we will denote by $\nabla^{\sharp,w}$ to distinguish it from the previous one.

\medskip
\begin{proposition}\label{connw}
The morphism $\str$ is a morphism of mixed complexes with $u$-connections
\begin{eqnarray*}
({\rrCe}(A_w), \bem+\bw, B^e, \nabla^{\flat,w})\to({\rrCe}(\C[Y]), \bem+\bw, B^e, \nabla^{\sharp,w})
\end{eqnarray*}
\end{proposition}
\noindent{\bf Proof} is given in section \ref{proofconnw}.
 \hfill $\Box$

\medskip

We will conclude by relating the $u$-connection $\nabla^{\sharp,w}$ on $({\rrCe}(\C[Y]), \bem+\bw, B^e)$ to the $u$-connection (\ref{nablaw}) on $(\Omega(Y)\,,\, -dw\,,\, d)$.

\medskip

\begin{proposition}\label{epsconn} The morphism $\epsilon$ from Proposition \ref{newmix} is a morphism of mixed complexes with $u$-connections
\begin{eqnarray*}
({\rrCe}(\C[Y]), \bem+\bw, B^e, \nabla^{\sharp,w})\to (\Omega(Y)\,,\, -dw\,,\, d,\nabla^{w})
\end{eqnarray*}
\end{proposition}
\noindent{\bf Proof} is given in section \ref{proofepsconn}.
 \hfill $\Box$

\medskip

Summarizing the results of Propositions \ref{newmix}, \ref{conntw}, \ref{connw}, \ref{epsconn} and using Propositions \ref{composition}, \ref{qishom}, we arrive at the following corollary:

\begin{corollary} The mixed complexes with $u$-connections
\begin{eqnarray*}
({\rC^e}(A_w), b^e, B^e, {\nabla}) \quad\text{and}\quad (\Omega(Y)\,,\, -dw\,,\, d,\nabla^{w})
\end{eqnarray*}
are homotopy gauge equivalent.
\end{corollary}
This corollary implies Theorem \ref{intro2} since the connections are, in fact, homotopy gauge equivalent in the stronger sense which we mentioned at the end of section \ref{uconnections}. Indeed, the $u$-quasi-isomorphism that relates the connections is an honest morphism of mixed complexes, and therefore establishes an isomorphism of the canonical $u$-lattices by the proof of Lemma \ref{standard}.

\medskip

\section{Applications of Theorem \ref{intro2}}\label{applicat}
Let $w$ be a polynomial on $Y=\C^k$ satisfying (\ref{crit}) and having an isolated singularity at the origin. For the time being, we will treat $w$ as the germ of an analytic function, as it is common in singularity theory.

One of the basic objects associated with $w$ is the Gauss-Manin system \cite{Br}. Roughly, it is a regular holonomic $D$-module on $|t|\ll1$ which extends the flat holomorphic bundle associated with the local system $\bigcup_{0<|t|\ll1}H^{k-1}(\tilde{w}^{-1}(t),\C)$ (here $\tilde{w}$ is the corresponding Milnor fibration). The Gauss-Manin system, together with a certain subspace known as the Brieskorn lattice, encodes various important characteristics and numerical invariants of the singularity such as the complex monodromy, the Steenbrink-Varchenko Hodge numbers, the spectrum, etc. (for a detailed treatment of the subject see e.g. \cite{SS}). The aim of this section is to explain a relationship between the Gauss-Manin system and the formal $D$-modules from Theorem \ref{intro2}, and thereby convince (hopefully) the reader that the aforementioned invariants can be recovered, {\it up to minor details}, starting from the category of matrix factorization.

Let us begin by recalling the definition of the Gauss-Manin system and the Brieskorn lattice.

The Gauss-Manin system can be defined as
\[
\mathcal{G}=\Omega^k(Y^{\rm an})_{(0)}[u^{-1}]/(-u^{-1}dw+d)\Omega^{k-1}(Y^{\rm an})_{(0)}[u^{-1}],
\]
where $\Omega^\bullet(Y^{\rm an})_{(0)}$ are the spaces of germs of holomorphic differential forms of appropriate degrees and the $D$-module structure is given by
\[
t(\eta\,u^{-i})=w\cdot\eta\, u^{-i}-i\eta\, u^{-i+1},\quad \partial_t(\eta\, u^{-i})=\eta\, u^{-i-1}\quad  (\eta\in \Omega^k(Y^{\rm an})_{(0)})
\]

An important property of this $D$-module is that the operator $\partial_t$ is invertible:
\[
\partial_t^{-1}(\eta\, u^{-i})=dw\wedge\eta'\, u^{-i}\quad\text{if}\quad d\eta'=\eta
\]
In other words, the obvious $\C[u^{-1}]$-module structure on $\mathcal{G}$ extends to a $\C[u,u^{-1}]$-module structure. It was observed by F.~Pham that, in fact, this extends further to a structure of $\C\{\!\{u\}\!\}[u^{-1}]$-vector space on $\mathcal{G}$ where $\C\{\!\{u\}\!\}:=\{\sum a_lu^l\in\C[[u]]\,|\,\exists R>0\,\,\text{s.\,t.}\,\sum a_l\frac{R^l}{l!}<\infty\}$. The dimension of $\cG$ is $\mu$, the Milnor number of the singularity.

The Brieskorn lattice is defined as follows:
\[
\mathcal{G}^{0}:=\Omega^k(Y^{\rm an})_{(0)}/dw\wedge d\,\Omega^{k-2}(Y^{\rm an})_{(0)}
\]
The canonical map $\Omega^k(Y^{\rm an})_{(0)}\to\mathcal{G}$ is known to induce an embedding $\mathcal{G}^{0}\hookrightarrow\mathcal{G}$ whose image  we will also denote by $\mathcal{G}^{0}$. $\mathcal{G}^{0}$ is a free $\C\{\!\{u\}\!\}$-submodule in $\mathcal{G}$ of rank $\mu$ and so $\mathcal{G}=\mathcal{G}^{0}[u^{-1}]$ (this is one reason to call $\mathcal{G}^{0}$ {\it lattice}).

Observe that $[t,u]=u^2$, so it is natural to set $\partial_u:=u^{-2}t$. Since the resulting $\C\{\!\{u\}\!\}[u^{-1}][\partial_u]$-module structure on $\cG$ encodes the original $D$-module structure (it is called the Fourier-Laplace transform of the latter), we will from now on forget the variable $t$ and work exclusively with $u$.

We don't know whether or not the analytic versions of the Gauss-Manin system and the Brieskorn lattice we just described admit any categorical interpretation. Fortunately, due to the regularity of the $D$-module $\cG$ it suffices for a number of questions (e.g. computing the aforementioned invariants of the singularity) to work with the following {\it formal} versions:
\[
\widehat{\cG}^{\,0}:=\C[[u]]\otimes_{\C\{\!\{u\}\!\}}\cG^{0},\quad \widehat{\cG}:=\widehat{\cG}^{\,0}[u^{-1}]
\]
These formal versions are precisely what we are able to reconstruct using the categorical approach. Namely, it was proved in \cite[Section 2]{Sch} that
\[
\widehat{\cG}^{\,0}\simeq H^k(\widehat{\Omega}^\bullet(Y)_{(0)}[[u]], -dw +ud),\quad \widehat{\cG}\simeq H^k(\widehat{\Omega}^\bullet(Y)_{(0)}((u)), -dw +ud)
\]
where $\widehat{\Omega}^\bullet(Y)_{(0)}$ are the spaces of formal differential forms. Moreover, the argument in \cite{Sch} shows that we would have got the same $\widehat{\cG}^{\,0}$ and $\widehat{\cG}$ if we started with the spaces ${\Omega}^\bullet(Y^{\rm alg})_{(0)}$ of differential forms with coefficients in the local algebra $\C[Y]_{(0)}$ of rational functions regular at the origin. This, together with the remark at the very end of the Introduction regarding the generality of our results, implies that the $\C[[u]]$-module $\widehat{\cG}^{\,0}$ and the $\C((u))$-linear space $\widehat{\cG}$ can be extracted from the category of matrix factorizations for the pair $(\C[Y]_{(0)}, w)$.

What about the $D$-module structure? It follows from the above definitions that the action of $\partial_u$ on $\widehat{\cG}\simeq H^k(\widehat{\Omega}^\bullet(Y)_{(0)}((u)), -dw +ud)$ is given by the formula
\begin{equation}\label{nablawtilde}
\partial_u=\frac{d}{du}+\frac{w}{u^2}
\end{equation}
The connection on the twisted de Rham cohomology induced by (\ref{nablaw}), which we are able to produce starting from the category of matrix factorizations, differs from (\ref{nablawtilde}) by the term $\frac{k}{2}\cdot\frac{1}{u}$, a kind of Tate twist (cf. the definition right after the equation (\ref{compofdmod})). Note that $k$, the number of variables $w$ depends on, {\it cannot} be extracted from the category $\mathrm{MF}(w)$ due to the celebrated Kn\"orrer periodicity: the categories $\mathrm{MF}(w)$ and $\mathrm{MF}(w+xy)$ are equivalent.

To summarize the previous discussion, Theorem \ref{intro2} allows us to reproduce, up to the above Tate twist, the formal Gauss-Manin system and the formal Brieskorn lattice starting from the categorical data. As a consequence, the category of matrix factorizations does remember about the invariants that are encoded in the Gauss-Manin system, but only modulo information involving the number $k$. Let us explain this using two examples.

Our first example is the {\it spectrum} of $w$. It is defined as follows. The $D$-module $\widehat{\mathcal{G}}$ is regular and, thus, it carries the canonical decreasing filtration ${V}^{\alpha}\widehat{\cG}$ (see section \ref{NCHodge}). This filtration induces a filtration on
\[
\widehat{\cG}^{\,0}/u\widehat{\cG}^{\,0}\simeq \widehat{\Omega}^k(Y)_{(0)}/dw\wedge\widehat{\Omega}^{k-1}(Y)_{(0)}\,(\simeq {\rm Jac}(w))
\]
by the images of ${V}^\alpha\widehat{\cG}\cap\widehat{\cG}^{\,0}$. Then the spectrum is the function
\[
d:\Q\to\Z,\quad d(\alpha):=\mathrm{dim}\,\mathrm{Gr}_V^\alpha(\widehat{\cG}^{\,0}/u\widehat{\cG}^{\,0})
\]
Its support is known to belong to the interval $(0, k)$; moreover, $d(\alpha)=d(k-\alpha)$ (cf. \cite[Example 1.8]{Sa}).

We can repeat this construction starting from the connection induced by (\ref{nablaw}) instead of (\ref{nablawtilde}). The resulting $V$-filtration will differ from the previous one by a shift: ${V}_{new}^{\alpha}\widehat{\cG}={V}_{old}^{\alpha+\frac{k}{2}}\widehat{\cG}$. Thus, we will obtain a shifted version of the spectrum $d_{new}(\alpha)=d_{old}(\alpha+\frac{k}{2})$, an even function
\[
d_{new}(-\alpha)=d_{new}(\alpha)
\]
with support in $(-\frac{k}{2}, \frac{k}{2})$.

Our second example is the {\it Steenbrink Hodge filtration} on the vanishing cohomology of $w$. The latter is nothing but $\rHP_k(\mathrm{MF}(w))$ with its non-commutative Hodge filtration (see section \ref{NCHodge}), as one can infer by inspecting the formulas presented in \cite{SS}. Since $k$ is not a categorical invariant, we recover the Steenbrink Hodge filtration only up to an undetermined shift.

\bigskip

\appendix

\section{Proofs for Section \ref{mixcompl}}
\medskip
\subsection{Proof of Lemma \ref{diffconn}}\label{proofdiffconn}
By (\ref{u-conn})
\begin{eqnarray*}
\left(\frac{df(u)}{du}+\cA'(u)f(u)-f(u)\cA(u)\right)(b+uB)\\
=\frac{df(u)}{du}(b+uB)+(b'+uB')\cA'(u)f(u)+\frac{1}{2u}(b'-uB')f(u)\\-(b'+uB')f(u)\cA(u)-\frac{1}{2u}f(u)(b-uB)\\
\end{eqnarray*}
Thus,
\begin{eqnarray*}
\left(\frac{df(u)}{du}+\cA'(u)f(u)-f(u)\cA(u)\right)(b+uB)-(b'+uB')\left(\frac{df(u)}{du}+\cA'(u)f(u)-f(u)\cA(u)\right)\\
=\frac{df(u)}{du}(b+uB)-(b'+uB')\frac{df(u)}{du}+\frac{1}{2u}(b'-uB')f(u)-\frac{1}{2u}f(u)(b-uB)
\end{eqnarray*}
To see that the latter equals 0, one can simplify the first part of the expression using the Leibniz rule as follows:
\begin{eqnarray*}
\frac{df(u)}{du}(b+uB)-(b'+uB')\frac{df(u)}{du}\\
=\frac{d}{du}\left(f(u)(b+uB)-(b'+uB')f(u)\right)-\left(f(u)\frac{d}{du}(b+uB)-\frac{d}{du}(b'+uB')f(u)\right)\\
=-f(u)B+B'f(u)
\end{eqnarray*}

\medskip

\subsection{Proof of Proposition \ref{hommaps}}\label{proofhommaps}

It suffices to show that any $u$-morphism, $u$-homotopic to 0, is a morphism of $u$-connections. Let
$$
f(u)=H(u)(b+uB)+(b'+uB')H(u)
$$
Then, computing modulo terms $u$-homotopic to 0
\begin{eqnarray*}
\frac{df(u)}{du}+\cA'(u)f(u)-f(u)\cA(u)\\
\sim H(u)B+B'H(u)+\cA'(u)H(u)(b+uB)+\cA'(u)(b'+uB')H(u)\\-H(u)(b+uB)\cA(u)-(b'+uB')H(u)\cA(u)\\
=H(u)B+B'H(u)+\cA'(u)H(u)(b+uB)+((b'+uB')\cA'(u)+\frac{1}{2u}(b'-uB'))H(u)\\-H(u)(\cA(u)(b+uB)-\frac{1}{2u}(b-uB))-(b'+uB')H(u)\cA(u)\\
%\end{eqnarray*}
%\begin{eqnarray*}
\sim H(u)B+B'H(u)+\frac{1}{2u}(b'-uB')H(u)+H(u)\frac{1}{2u}(b-uB)\\=\frac{1}{2u}(b'+uB')H(u)+H(u)\frac{1}{2u}(b+uB)
\end{eqnarray*}
($\sim$ above stands  for `$u$-homotopic').
%\begin{eqnarray*}
%\end{eqnarray*}
\medskip

\subsection{Proof of Proposition \ref{qishom}}\label{proofqishom}

Due to Proposition \ref{easylemma}, we only need to prove
\begin{lemma}\label{ofm}
If one of the $u$-morphisms $p(u)$, $\iota(u)$ is a morphism of $u$-connections then the other one is also a morphisms of $u$-connections.
\end{lemma}
\noindent{\bf Proof.} Let us assume for instance that $\iota(u)$ is a morphism of $u$-connections:
$$
\frac{d\iota(u)}{du}+\cA'(u)\iota(u)-\iota(u)\cA(u)\sim0
$$
Then
\begin{eqnarray*}
p(u)\frac{d\iota(u)}{du}p(u)+p(u)\cA'(u)\iota(u)p(u)-p(u)\iota(u)\cA(u)p(u)\sim0
\end{eqnarray*}
or, using the Leibniz rule
\begin{eqnarray}\label{equivhom}
\frac{d}{du}(p(u)\iota(u))p(u)-\frac{dp(u)}{du}\iota(u)p(u)+p(u)\cA'(u)\iota(u)p(u)-p(u)\iota(u)\cA(u)p(u)\sim0
\end{eqnarray}
Using
\begin{eqnarray}\label{homoto}
p(u)\iota(u)=id_{\rC} + (b+uB)H(u)+H(u)(b+uB)
\end{eqnarray}
\begin{eqnarray}\label{homoto1}
\iota(u)p(u)=id_{\rC'} + (b'+uB')H'(u)+H'(u)(b'+uB')
\end{eqnarray}
(\ref{equivhom}) is equivalent to
\begin{eqnarray*}
\frac{dp(u)}{du}+\cA(u)p(u)-p(u)\cA'(u)\\
\sim \frac{d}{du}\left((b+uB)H(u)+H(u)(b+uB)\right)p(u)-\frac{dp(u)}{du}\left((b'+uB')H'(u)+H'(u)(b'+uB')\right)\\
+p(u)\cA'(u)\left((b'+uB')H'(u)+H'(u)(b'+uB')\right)-\left((b+uB)H(u)+H(u)(b+uB)\right)\cA(u)p(u)
\end{eqnarray*}
Thus, it suffices to show that the latter expression is $u$-homotopic to 0. By (\ref{u-conn}) and the Leibniz rule, it is $u$-homotopic to
\begin{eqnarray*}
(BH(u)+H(u)B)p(u)-\frac{dp(u)}{du}\left((b'+uB')H'(u)+H'(u)(b'+uB')\right)\\
+\frac{1}{2u}p(u)(b'-uB')H'(u)+\frac{1}{2u}H(u)(b-uB)p(u)\\
\sim BH(u)p(u)+H(u)Bp(u)-Bp(u)H'(u)+p(u)B'H'(u)\\
+\frac{1}{2u}p(u)(b'-uB')H'(u)+\frac{1}{2u}H(u)(b-uB)p(u)\\
%\end{eqnarray*}
%\begin{eqnarray*}
=BH(u)p(u)-Bp(u)H'(u)
+\frac{1}{2u}p(u)(b'+uB')H'(u)+\frac{1}{2u}H(u)(b+uB)p(u)\\
=BH(u)p(u)-Bp(u)H'(u)
+\frac{1}{2u}(b+uB)p(u)H'(u)+\frac{1}{2u}H(u)p(u)(b'+uB')\\
\sim BH(u)p(u)+\frac{1}{2u}(b-uB)p(u)H'(u)-\frac{1}{2u}(b+uB)H(u)p(u)\\
=\frac{1}{2u}(b-uB)p(u)H'(u)-\frac{1}{2u}(b-uB)H(u)p(u)\\
=\frac{1}{2u}(b-uB)(p(u)H'(u)-H(u)p(u))
\end{eqnarray*}
The latter coincides with
$$
(b+uB)\left(\cA(u)(H(u)p(u)-p(u)H'(u))\right)+\left(\cA(u)(H(u)p(u)-p(u)H'(u))\right)(b'+uB')
$$
Indeed, this follows from (\ref{u-conn}) and the equality
\begin{eqnarray*}
(H(u)p(u)-p(u)H'(u))(b'+uB')=-(b+uB)(H(u)p(u)-p(u)H'(u))
\end{eqnarray*}
which, in its turn, is due to (\ref{homoto}) and (\ref{homoto1}).
\hfill $\Box$

\medskip

\subsection{Proof of Proposition \ref{gauge}}\label{gaugeequiv}
First of all, we need to show that (\ref{gaugeconn}) {\it is} a $u$-connection, i.e. we need to verify the property (\ref{u-conn}):
\begin{eqnarray*}
\left[\frac{d}{du}+p(u)\frac{d\iota(u)}{du}+p(u)\cA'(u)\iota(u)+\frac{1}{2u}H(u)(b-uB), b+uB\right]\\
=B+p(u)\left[\frac{d\iota(u)}{du},b+uB\right]+p(u)[\cA'(u),b+uB]\iota(u)+\left[\frac{1}{2u}H(u)(b-uB), b+uB\right]\\
=B-p(u)\left[\iota(u),\frac{d}{du}(b+uB)\right]+\frac{1}{2u}p(u)(b-uB)\iota(u)+\left[\frac{1}{2u}H(u)(b-uB), b+uB\right]\\
=B-p(u)\left[\iota(u),B\right]+\frac{1}{2u}p(u)(b+uB-2uB)\iota(u)+\left[\frac{1}{2u}H(u)(b-uB), b+uB\right]\\
=B-p(u)\iota(u)B+p(u)B\iota(u)+\frac{1}{2u}p(u)\iota(u)(b+uB)-p(u)B\iota(u)+\left[\frac{1}{2u}H(u)(b-uB), b+uB\right]\\
=\frac{1}{2u}(b+uB)-((b+uB)H(u)+H(u)(b+uB))B+\frac{1}{2u}((b+uB)H(u)+H(u)(b+uB))(b+uB)\\+\left[\frac{1}{2u}H(u)(b-uB), b+uB\right]\\
=\frac{1}{2u}(b+uB)
\end{eqnarray*}

By Lemma \ref{ofm} to complete the proof it is enough to show that $\iota(u)$ is a morphism from $\nabla$ to $\nabla'$, i.e.
\begin{eqnarray*}
\frac{d\iota(u)}{du}+\cA'(u)\iota(u)-(\iota(u)p(u)\frac{d\iota(u)}{du}+\iota(u)p(u)\cA'(u)\iota(u)+\frac{1}{2u}\iota(u)H(u)(b-uB))\sim0
\end{eqnarray*}
This can be rewritten using (\ref{homoto1}) as follows:
\begin{eqnarray*}
-\left((b'+uB')H'(u)+H'(u)(b'+uB')\right)\frac{d\iota(u)}{du}-\left((b'+uB')H'(u)+H'(u)(b'+uB')\right)\cA'(u)\iota(u)\\
-\frac{1}{2u}\iota(u)H(u)(b-uB)\sim0
\end{eqnarray*}
Using the Leibniz rule and (\ref{u-conn}), the left-hand side is $u$-homotopic to
\begin{eqnarray*}
-\frac{d}{du}\left((b'+uB')H'(u)\iota(u)+H'(u)\iota(u)(b+uB)\right)+\frac{d}{du}\left((b'+uB')H'(u)+H'(u)(b'+uB')\right)\iota(u)\\
+\frac{1}{2u}H'(u)(b'-uB')\iota(u)-\frac{1}{2u}\iota(u)H(u)(b-uB)\\
\sim
-B'H'(u)\iota(u)-H'(u)\iota(u)B+B'H'(u)\iota(u)+H'(u)B'\iota(u)\\
+\frac{1}{2u}H'(u)(b'-uB')\iota(u)-\frac{1}{2u}\iota(u)H(u)(b-uB)\\
=
-H'(u)\iota(u)B+H'(u)B'\iota(u)
+\frac{1}{2u}H'(u)(b'-uB')\iota(u)-\frac{1}{2u}\iota(u)H(u)(b-uB)
\end{eqnarray*}
\begin{eqnarray*}
=
-H'(u)\iota(u)B
+\frac{1}{2u}H'(u)(b'+uB')\iota(u)-\frac{1}{2u}\iota(u)H(u)(b-uB)\\
=
-H'(u)\iota(u)B
+\frac{1}{2u}H'(u)\iota(u)(b+uB)-\frac{1}{2u}\iota(u)H(u)(b-uB)\\
=
\frac{1}{2u}H'(u)\iota(u)(b-uB)-\frac{1}{2u}\iota(u)H(u)(b-uB)\\
=
(H'(u)\iota(u)-\iota(u)H(u))\frac{1}{2u}(b-uB)
\end{eqnarray*}
Repeating the argument at the end of the proof of Lemma \ref{ofm}, the latter is $u$-homotopic to 0.

\medskip
\section{Useful formulas}\label{auxlemmas}
In this appendix we list various ``commutation relations'' between the operators introduced in section \ref{cyclcompl}.

To avoid duplicating results, we will view $\ld$'s, $\lm$'s, $\tau$ etc. as operators on $\rC(A^e)$ where $A^e$ is the unitalization of $A$. Then most of the formulas below may be interpreted in two ways, namely, as equalities between operators whose domain is either $\rC(A)$ or $\rC^+(A)$. Similarly, we will use the symbols like $\bm$, $\bmv$ or $\bmh$ to denote the operators on $\rC(A^e)$ given by the same formulas as in the main body of the text.
Let us also denote by $\lm^{(\ast)}$ the operator on $\rC(A^e)$ defined by
\begin{eqnarray*}
\lm^{(\ast)}_{n+1}:=\lm^{(n)}_{n+1}
\end{eqnarray*}
Using this notation, the operators $\bmh$, $\bmv$ and $\bm$ are related as follows
\begin{eqnarray*}
\bm=\bmv+\bmh=\lm^{(0)}+\lm^{(\ast)}+\bmh
\end{eqnarray*}

We will also omit subscripts in symbols like $\lm^{(i)}_{n+1}$, when it does not lead to confusion.

\medskip

\begin{lemma}\label{Mitau}
\begin{eqnarray*}
a)\quad \ld^{(k)}\tau^j=\tau^j\ld^{(k+j)}
\end{eqnarray*}
\begin{eqnarray*}
b)\quad \lm_{n}^{(k)}\tau_{n}^{j}=\begin{cases} \tau_{n-1}^{j}\lm_{n}^{(k+j)} & k+j\leq n-1\\ \tau_{n-1}^{j-1}\lm_{n}^{(k+j-n)} & k+j\geq n\end{cases}
\end{eqnarray*}
\begin{eqnarray*}
c)\quad \bd\tau=\tau\bd
\end{eqnarray*}
\begin{eqnarray*}
d)\quad (\bm-\lm^{(\ast)})\tau=\tau(\bm-\lm^{(0)})
\end{eqnarray*}
\begin{eqnarray*}
e)\quad \bm(1-\tau)=(1-\tau)(\bm-\lm^{(0)})
\end{eqnarray*}
\begin{eqnarray*}
f)\quad (\bm-\lm^{(\ast)})N=N\bm
\end{eqnarray*}
\end{lemma}

\noindent{\bf Proof.} a) and b) follow from definitions; c) follows from a);  d) and e) follow from b). Finally, both hand sides of f) are equal to $N\lm^{(0)}N$.
\hfill $\Box$

\medskip

\begin{lemma}\label{lemmaDiMi}
\begin{eqnarray*}
a)\quad \ld^{(i)}\ld^{(j)}=-\ld^{(j)}\ld^{(i)}\quad (\forall i,j)
\end{eqnarray*}
\begin{eqnarray*}
b)\quad \lm^{(i)}\lm^{(j)}=-\lm^{(j)}\lm^{(i+1)}\quad (j\leq i)
\end{eqnarray*}
\begin{eqnarray*}
c)\quad \ld^{(i)}_{n}\lm^{(j)}_{n+1}=
\begin{cases}
-\lm^{(j)}_{n+1}\ld^{(i+1)}_{n+1} \quad 0\leq j< i\leq n-1\\
-\lm^{(j)}_{n+1}\ld^{(i)}_{n+1} \quad 0\leq i< j\leq n-1\\
-\lm^{(n)}_{n+1}\ld^{(i)}_{n+1} \quad 1\leq i\leq n-1, j=n\\
-\lm^{(i)}_{n+1}\ld^{(i+1)}_{n+1}-\lm^{(i)}_{n+1}\ld^{(i)}_{n+1} \quad 0\leq i=j\leq n
\end{cases}
\end{eqnarray*}
(When $i=n$ in the last line, $\ld^{(i)}_{n}$ and $\ld^{(i+1)}_{n+1}$ stand for $\ld^{(0)}_{n}$ and $\ld^{(0)}_{n+1}$, respectively.)
\end{lemma}

\noindent{\bf Proof.}
a) Using the definitions of the $\ld$'s, $\ld^{(0)}_{n+1}\ld^{(l)}_{n+1}=-\ld^{(l)}_{n+1}\ld^{(0)}_{n+1}$ when $l\neq0$. This holds true for $l=0$ as well (both hand sides vanish). The general case reduces to this special one by means of Lemma \ref{Mitau} a).

b) Again, by the definition of the $\lm$'s, $\lm^{(l)}_{n}\lm^{(0)}_{n+1}=-\lm^{(0)}_{n}\lm^{(l+1)}_{n+1}$ when $1\leq l\leq n-2$. This holds true for $l=0$ and $n-1$, as in these cases it is equivalent to the associativity of the multiplication. The general case reduces to this special one by means of Lemma \ref{Mitau} b).

Part c) is proved similarly (the last case is just the Leibniz rule).
 \hfill $\Box$

\medskip
The following lemma is a corollary of the previous one:
\begin{lemma}\label{lemmabdMi}
\begin{eqnarray*}
a)\quad \bd\ld^{(i)}=-\ld^{(i)}\bd
\end{eqnarray*}
\begin{eqnarray*}
b)\quad \bd\lm^{(i)}=-\lm^{(i)}\bd
\end{eqnarray*}
\end{lemma}

\medskip

The remaining formulas involve the operator $h^e$.

\begin{lemma}\label{lemmaDiMiH}
\begin{eqnarray*}
a)\quad \ld^{(0)}\h^e=0, \quad \ld^{(i)}\h^e=-\h^e\ld^{(i-1)} \quad (i\geq1)
\end{eqnarray*}
\begin{eqnarray*}
b)\quad \lm^{(0)}\h^e=id, \quad \lm^{(i)}_{n+2}\h^e_{n+1}=-\h^e_{n}\lm^{(i-1)}_{n+1} \quad (1\leq i\leq n),  \quad \lm^{(\ast)}\h^e=-\tau^{-1}
\end{eqnarray*}
\end{lemma}
\noindent{\bf Proof.} The equalities are straightforward except, perhaps, the very last one. We need to show that
$
\tau_{n+1}\lm^{(n+1)}_{n+2}\h^e_{n+1}=-id_{n+1}.
$
By (\ref{Mi})
\begin{eqnarray*}
\tau_{n+1}\lm^{(n+1)}_{n+2}\h^e_{n+1}(a_0[a_1|a_2|\ldots |a_n])=\tau_{n+1}\lm^{(n+1)}_{n+2}(e[a_0|a_1|a_2|\ldots |a_n])\\
=-(-1)^{|sa_n|\sum_{k=0}^{n-1}|sa_k|}\tau_{n+1}(a_n[a_0|\ldots |a_{n-1}])\\
=-(-1)^{|sa_n|\sum_{k=0}^{n-1}|sa_k|}(-1)^{|sa_n|\sum_{k=0}^{n-1}|sa_k|)}a_0[a_1|\ldots |a_{n}]
=-a_0[a_1|\ldots |a_{n}]
\end{eqnarray*}
 \hfill $\Box$

\medskip

From this lemma one easily deduces
\begin{lemma}\label{lemmabh}
\begin{eqnarray*}
a)\quad \bd h^e=-h^e\bd
\end{eqnarray*}
\begin{eqnarray*}
b)\quad \bmv h^e=1-\tau^{-1}
\end{eqnarray*}
\begin{eqnarray*}
c)\quad \bmh h^e=-h^e(\bm-\lm^{(\ast)})=-h^e(\lm^{(0)}+\bmh)
\end{eqnarray*}
\end{lemma}

\medskip

The formulas above hold true if we replace $\rC(A^e)$ by $\rC(A)$ and $h^e$ by $h$.

\medskip

\section{Proofs for Section \ref{cyclcompl}}
\subsection{Proof of Proposition \ref{quasiisocycl}}\label{homotopy}

We will start by proving
\begin{lemma} $\iota(u)$ and $p(u)$ are $u$-morphisms of mixed complexes.
\end{lemma}
\noindent{\bf Proof.}
Let us show that $\iota(u)$ is a $u$-morphism. The only part that is not immediate is the equality
\begin{eqnarray*}
B^e\begin{pmatrix}id_{\rC (A)} \\ 0  \end{pmatrix}-\begin{pmatrix}id_{\rC (A)} \\ 0  \end{pmatrix}B=\begin{pmatrix}0 \\ h^ehN  \end{pmatrix}b-b^e\begin{pmatrix}0 \\ h^ehN  \end{pmatrix}
\end{eqnarray*}
By Lemmas \ref{Mitau} c), \ref{lemmabh} a)
\begin{eqnarray*}
\begin{pmatrix}0 \\ h^ehN  \end{pmatrix}\bd-\bed\begin{pmatrix}0 \\ h^ehN  \end{pmatrix}=0
\end{eqnarray*}
so we need to show that
\begin{eqnarray*}
B^e\begin{pmatrix}id_{\rC (A)} \\ 0  \end{pmatrix}-\begin{pmatrix}id_{\rC (A)} \\ 0  \end{pmatrix}B=\begin{pmatrix}0 \\ h^ehN  \end{pmatrix}\bm-\bem\begin{pmatrix}0 \\ h^ehN  \end{pmatrix}
\end{eqnarray*}
or, equivalently
\begin{eqnarray*}
\begin{pmatrix}-(1-\tau^{-1})hN \\ h^eN  \end{pmatrix}=\begin{pmatrix}-\bmv h^ehN \\ h^ehN\bm-\bmh h^ehN  \end{pmatrix}
\end{eqnarray*}
The first components of the two vectors coincide by Lemma \ref{lemmabh} b), so it remains to prove that
\begin{eqnarray*}
h^eN=h^ehN\bm-\bmh h^ehN.
\end{eqnarray*}
The latter follows from Lemmas \ref{Mitau} f), \ref{lemmabh} c), \ref{lemmaDiMiH} b).

Let us show now that $p(u)$ is a $u$-morphism. The only non-trivial part is the equality
\begin{eqnarray*}
\bm(1-\tau^{-1})h\lm^{(0)}=\bmv+(1-\tau^{-1})h\lm^{(0)}\bmh
\end{eqnarray*}
of operators from $\rC^+(A)$ to $\rC(A)$. Since $h^e: \rC(A)\to\rC^+(A)$ is bijective, it is enough to prove that
\begin{eqnarray*}
\bm(1-\tau^{-1})h\lm^{(0)}h^e=\bmv h^e+(1-\tau^{-1})h\lm^{(0)}\bmh h^e
\end{eqnarray*}
viewed as operators on $\rC(A)$. Since $\lm^{(0)}h^e=1$ and $\bmv h^e=(1-\tau^{-1})$ (Lemmas \ref{lemmaDiMiH}, \ref{lemmabh}),
the latter is equivalent to
\begin{eqnarray*}
\bm(1-\tau^{-1})h=(1-\tau^{-1})+(1-\tau^{-1})h\lm^{(0)}\bmh h^e
\end{eqnarray*}
which, by Lemmas \ref{Mitau} e), d),  reduces to
\begin{eqnarray*}
(\lm^{(0)}+\bmh)h=1+h\lm^{(0)}\bmh h^e
\end{eqnarray*}
What remains is to use Lemmas \ref{lemmaDiMiH} b), \ref{lemmabh} c).
 \hfill $\Box$

To complete the proof of Proposition \ref{quasiisocycl}, we will prove

\begin{lemma}\label{iotap}
\begin{eqnarray}\label{piota}
p(u)\iota(u)=id_{\rC (A)} + (b+uB)H(u)+H(u)(b+uB), \quad H(u):=u(1-\tau^{-1})hhhN
\end{eqnarray}
\begin{eqnarray}\label{iotapi}
\qquad\iota(u)p(u)=id_{\rC^e(A)} + (b^e+uB^e)H^e(u)+H^e(u)(b^e+uB^e), \quad H^e(u):=\begin{pmatrix}0 & 0 \\ 0 & h^eh\lm^{(0)}  \end{pmatrix}
\end{eqnarray}
\end{lemma}
\noindent{\bf Proof.}

\noindent {\it(\ref{piota}):}  Observe that
\begin{eqnarray*}
p(u)\iota(u)=id_{\rC (A)} + u(1-\tau^{-1})hhN
\end{eqnarray*}
We will show that
\begin{eqnarray*}
bH(u)+H(u)b=u(1-\tau^{-1})hhN, \quad  BH(u)+H(u)B=0
\end{eqnarray*}
The second equality is obvious. Let us prove the first one. By Lemmas \ref{Mitau} c), \ref{lemmabh} a)
\begin{eqnarray*}
\bd H(u)+H(u)\bd=0
\end{eqnarray*}
so we only need to show that
\begin{eqnarray*}
\bm (1-\tau^{-1})hhhN+(1-\tau^{-1})hhhN\bm=(1-\tau^{-1})hhN
\end{eqnarray*}
This follows from Lemmas \ref{Mitau} e),d),f), \ref{lemmaDiMiH} b), \ref{lemmabh} c).

\medskip

\noindent {\it (\ref{iotapi}):}  Observe that
\begin{eqnarray*}
\iota(u) p(u)=id_{\rC^e(A)} + \begin{pmatrix}0 & (1-\tau^{-1})h\lm^{(0)} \\ 0 & -id_{\rC^+(A)}  \end{pmatrix}+u\begin{pmatrix}0 & 0 \\ h^ehN & 0  \end{pmatrix}
\end{eqnarray*}
We will show that
\begin{eqnarray*}
b^eH^e(u)+H^e(u)b^e=\begin{pmatrix}0 & (1-\tau^{-1})h\lm^{(0)} \\ 0 & -id_{\rC^+(A)}  \end{pmatrix}, \quad  B^eH^e(u)+H^e(u)B^e=\begin{pmatrix}0 & 0 \\ h^ehN & 0  \end{pmatrix}
\end{eqnarray*}
Again, the latter equality is obvious and we will prove only the former one. By Lemmas \ref{Mitau} c), \ref{lemmabh} a), \ref{lemmabdMi} b)
\begin{eqnarray*}
\bed H^e(u)+H^e(u)\bed=0
\end{eqnarray*}
so it remains to show that
\begin{eqnarray*}
\begin{pmatrix}\bm &  \bmv \\0 & \bmh \end{pmatrix}\begin{pmatrix}0 & 0 \\ 0 & h^eh\lm^{(0)}  \end{pmatrix}
+\begin{pmatrix}0 & 0 \\ 0 & h^e h\lm^{(0)}  \end{pmatrix}\begin{pmatrix}\bm &  \bmv \\0 & \bmh \end{pmatrix}
=\begin{pmatrix}0 & (1-\tau^{-1})h\lm^{(0)} \\ 0 & -id_{\rC^+(A)}  \end{pmatrix}
\end{eqnarray*}
or, equivalently
\begin{eqnarray*}
\bmv h^eh \lm^{(0)}=(1-\tau^{-1})h \lm^{(0)},\quad \bmh h^eh \lm^{(0)}+h^eh \lm^{(0)}\bmh=-id_{\rC^+(A)}
\end{eqnarray*}
The first equality is due to Lemma \ref{lemmabh} b). The second equality is due to
Lemmas \ref{lemmabh} c), \ref{lemmaDiMiH} b), \ref{lemmaDiMi} b).
\hfill $\Box$

\medskip

\subsection{Proof of Proposition \ref{conncycl}}\label{proofUV}

That the operators $\nabla$ and $\nabla^\circ$ are dual to each other follows from the formulas we wrote out while proving Proposition \ref{opp}. Thus, it suffices to show that only one of the operators ${\nabla}$, ${\nabla^{\circ}}$ is a $u$-connection. Let us show that, say, ${\nabla^{\circ}}$ is a $u$-connection.

\medskip

\begin{lemma}\label{UV}
\begin{eqnarray*}
a)\quad [\Ud,\bed]=[\Ud,\bem]=0
\end{eqnarray*}
\begin{eqnarray*}
b)\quad [\Vd,B^e]=[\Vd,\bed]=0
\end{eqnarray*}
\begin{eqnarray*}
c)\quad [\Vd,\bem]+[\Ud,B^e]=\frac12\bed
\end{eqnarray*}
\begin{eqnarray*}
d)\quad [\Gamma,\bed]=0,  \quad  [\Gamma,\bem]=\frac12\bem,   \quad [\Gamma,B^e]=-\frac12B^e
\end{eqnarray*}
\end{lemma}
\noindent{\bf Proof.}
\noindent {\it Part a):} The first commutator equals 0 by Lemma \ref{lemmabdMi}. Let us prove the second equality. Writing both $\bem$ and $\Ud$ in the matrix notation, the second equality is equivalent to
\begin{eqnarray}\label{bU}
\lm_{n}^{(n-1)}\ld_{n}^{(n-1)}\cdot\sum_{i=0}^{n}\lm_{n+1}^{(i)}=\sum_{i=0}^{n-1}\lm_{n}^{(i)}\cdot\lm_{n+1}^{(n)}\ld_{n+1}^{(n)}
\end{eqnarray}
where both hand sides are viewed as operators either on $\rC (A)$ or on $\rC^+ (A)$.

By Lemma \ref{lemmaDiMi} c)
\begin{eqnarray*}
\lm_{n}^{(n-1)}\ld_{n}^{(n-1)}\sum_{i=0}^{n}\lm_{n+1}^{(i)}=\lm_{n}^{(n-1)}\left(\ld_{n}^{(n-1)}\sum_{i=0}^{n-2}\lm_{n+1}^{(i)}+\ld_{n}^{(n-1)}\lm_{n+1}^{(n-1)}+\ld_{n}^{(n-1)}\lm_{n+1}^{(n)}\right)\\
=-\lm_{n}^{(n-1)}\left(\sum_{i=0}^{n-2}\lm_{n+1}^{(i)}\ld_{n+1}^{(n)}+\lm_{n+1}^{(n-1)}\ld_{n+1}^{(n-1)}+\lm_{n+1}^{(n-1)}\ld_{n+1}^{(n)}+\lm_{n+1}^{(n)}\ld_{n+1}^{(n-1)}\right)\\
=-\lm_{n}^{(n-1)}\left(\sum_{i=0}^{n-1}\lm_{n+1}^{(i)}\ld_{n+1}^{(n)}+\lm_{n+1}^{(n-1)}\ld_{n+1}^{(n-1)}+\lm_{n+1}^{(n)}\ld_{n+1}^{(n-1)}\right)\\
=-\lm_{n}^{(n-1)}\sum_{i=0}^{n-1}\lm_{n+1}^{(i)}\ld_{n+1}^{(n)}-\left(\lm_{n}^{(n-1)}\lm_{n+1}^{(n-1)}+\lm_{n}^{(n-1)}\lm_{n+1}^{(n)}\right)\ld_{n+1}^{(n-1)}
\end{eqnarray*}
By Lemma \ref{lemmaDiMi} b), the latter expression equals the right-hand side of (\ref{bU}).

\medskip

\noindent {\it Part b):} The only non-trivial statement is the second equality. It is a consequence of Lemmas \ref{Mitau} c), \ref{lemmabdMi} a), \ref{lemmabh} a).

\medskip

\noindent {\it Part c):} Observe that
\begin{eqnarray*}
2[\Vd,\bem]_{n}+2[\Ud, B^e]_n=\begin{pmatrix}\ast_{11} & 0 \\ \ast_{21} & \ast_{22} \end{pmatrix}
\end{eqnarray*}
where
\begin{eqnarray*}%\label{ast11}
\ast_{11}=-\bmv_{n+1} h^e_{n}\sum_{i=1}^{n-1}\sum_{j=1}^{i}\tau_{n}^{j}\ld_{n}^{(i)}+\lm^{(n)}_{n+1}\ld^{(n)}_{n+1} h^e_nN_n
\end{eqnarray*}
\begin{eqnarray*}
\ast_{21}=h^e_{n-1}\sum_{i=1}^{n-2}\sum_{j=1}^{i}\tau_{n-1}^{j}\ld_{n-1}^{(i)}\bm_{n}-\bmh_{n+1} h^e_{n}\sum_{i=1}^{n-1}\sum_{j=1}^{i}\tau_{n}^{j}\ld_{n}^{(i)}-h^e_{n-1}N_{n-1}\lm^{(n-1)}_n\ld^{(n-1)}_n
\end{eqnarray*}
\begin{eqnarray*}
\ast_{22}= h^e_{n-1}\sum_{i=1}^{n-2}\sum_{j=1}^{i}\tau_{n-1}^{j}\ld_{n-1}^{(i)}\bmv_{n}-h^e_{n-1}N_{n-1}\lm^{(n-1)}_n\ld^{(n-1)}_n
\end{eqnarray*}
By Lemma \ref{lemmaDiMiH},
\begin{eqnarray*}
\ast_{11}=-(1-\tau^{-1}_n)\sum_{i=1}^{n-1}\sum_{j=1}^{i}\tau_{n}^{j}\ld_{n}^{(i)}-\lm^{(n)}_{n+1}h^e_n\ld^{(n-1)}_{n} N_n\\
=-\sum_{i=1}^{n-1}(\tau_{n}^{i}-1)\ld_{n}^{(i)}+\tau_n^{-1}\ld^{(n-1)}_{n} N_n
=-\sum_{i=1}^{n-1}\tau_{n}^{i}\ld_{n}^{(i)}+\sum_{i=1}^{n-1}\ld_{n}^{(i)}+\ld^{(0)}_{n} N_n\\
=-\sum_{i=1}^{n-1}\ld_{n}^{(0)}\tau_{n}^{i}+\sum_{i=1}^{n-1}\ld_{n}^{(i)}+\ld^{(0)}_{n} N_n
=-\sum_{i=0}^{n-1}\ld_{n}^{(0)}\tau_{n}^{i}+\sum_{i=0}^{n-1}\ld_{n}^{(i)}+\ld^{(0)}_{n} N_n
=\bd_n
\end{eqnarray*}

To simplify $\ast_{22}$, let us use the fact that $h^e: \rC(A)\to\rC^+(A)$ is bijective. For $\alpha\in A\otimes sA^{\otimes (n-2)}$
\begin{eqnarray*}
h^e_{n-1}\sum_{i=1}^{n-2}\sum_{j=1}^{i}\tau_{n-1}^{j}\ld_{n-1}^{(i)}\bmv_{n}h^e_{n-1}(\alpha)-h^e_{n-1}N_{n-1}\lm^{(n-1)}_n\ld^{(n-1)}_nh^e_{n-1}(\alpha)\\
=h^e_{n-1}\sum_{i=1}^{n-2}\sum_{j=1}^{i}\tau_{n-1}^{j}\ld_{n-1}^{(i)}(1-\tau_{n-1}^{-1})(\alpha)+h^e_{n-1}N_{n-1}\lm^{(n-1)}_nh^e_{n-1}\ld^{(n-2)}_{n-1}(\alpha)\\
= h^e_{n-1}\sum_{i=1}^{n-2}\sum_{j=1}^{i}\tau_{n-1}^{j}\ld_{n-1}^{(i)}(\alpha)-h^e_{n-1}\sum_{i=1}^{n-2}\sum_{j=1}^{i}\tau_{n-1}^{j}\ld_{n-1}^{(i)}\tau^{-1}_{n-1}(\alpha)
-h^e_{n-1}N_{n-1}\tau_{n-1}^{-1}\ld^{(n-2)}_{n-1}(\alpha)
\end{eqnarray*}
\begin{eqnarray*}
=h^e_{n-1}\sum_{i=1}^{n-2}\sum_{j=1}^{i}\tau_{n-1}^{j}\ld_{n-1}^{(i)}(\alpha)-h^e_{n-1}\sum_{i=1}^{n-2}\sum_{j=1}^{i}\tau_{n-1}^{j-1}\ld_{n-1}^{(i-1)}(\alpha)-h^e_{n-1}N_{n-1}\ld^{(n-2)}_{n-1}(\alpha)\\
=h^e_{n-1}\sum_{j=1}^{n-2}\tau_{n-1}^{j-1}\ld_{n-1}^{(n-2)}(\alpha)-h^e_{n-1}\sum_{i=0}^{n-3}\ld_{n-1}^{(i)}(\alpha)-h^e_{n-1}N_{n-1}\ld^{(n-2)}_{n-1}(\alpha)\\
=h^e_{n-1}\sum_{j=0}^{n-2}\tau_{n-1}^{j-1}\ld_{n-1}^{(n-2)}(\alpha)-h^e_{n-1}\sum_{i=0}^{n-2}\ld_{n-1}^{(i)}(\alpha)-h^e_{n-1}N_{n-1}\ld^{(n-2)}_{n-1}(\alpha)
=-h^e_{n-1}\sum_{i=0}^{n-2}\ld_{n-1}^{(i)}(\alpha)
\end{eqnarray*}
By Lemma \ref{lemmabh} a) the latter equals $\bd_n(h^e_{n-1}(\alpha))$. Thus, $\ast_{22}=\bd_n$.

It remains to show that $\ast_{21}=0$. By Lemma \ref{lemmabh} c)
it is enough to prove that
\begin{eqnarray}\label{ast210}
\sum_{i=1}^{n-2}\sum_{j=1}^{i}\tau_{n-1}^{j}\ld_{n-1}^{(i)}\bm_{n}+(\bm_n-\lm_{n}^{(n-1)})\sum_{i=1}^{n-1}\sum_{j=1}^{i}\tau_{n}^{j}\ld_{n}^{(i)}-N_{n-1}\lm^{(n-1)}_n\ld^{(n-1)}_n=0
\end{eqnarray}
By Lemma \ref{lemmaDiMi} c)
\begin{eqnarray*}
\ld_{n-1}^{(i)}\bm_{n}=-\left(\sum_{k=0}^i \lm_{n}^{(k)}\right) \ld_{n}^{(i+1)}-\left(\sum_{k=i}^{n-1} \lm_{n}^{(k)}\right) \ld_{n}^{(i)}
\end{eqnarray*}
for $i\geq1$. Therefore
\begin{eqnarray*}
\sum_{i=1}^{n-2}\sum_{j=1}^{i}\tau_{n-1}^{j}\ld_{n-1}^{(i)}\bm_{n}
=-\sum_{i=1}^{n-2}\sum_{j=1}^{i}\sum_{k=0}^i\tau_{n-1}^{j}\lm_{n}^{(k)}\ld_{n}^{(i+1)}-\sum_{i=1}^{n-2}\sum_{j=1}^{i}\sum_{k=i}^{n-1}\tau_{n-1}^{j}\lm_{n}^{(k)} \ld_{n}^{(i)}\\
=-\sum_{i=2}^{n-1}\sum_{j=1}^{i-1}\sum_{k=0}^{i-1}\tau_{n-1}^{j}\lm_{n}^{(k)}\ld_{n}^{(i)}-\sum_{i=1}^{n-2}\sum_{j=1}^{i}\sum_{k=i}^{n-1}\tau_{n-1}^{j}\lm_{n}^{(k)} \ld_{n}^{(i)}
\end{eqnarray*}
Thus, the coefficient ''in front of'' $\ld_{n}^{(i)}$ on the left-hand side of (\ref{ast210}) equals
\begin{eqnarray}\label{astfinish}
-\sum_{j=1}^{i-1}\sum_{k=0}^{i-1}\tau_{n-1}^{j}\lm_{n}^{(k)}-\sum_{j=1}^{i}\sum_{k=i}^{n-1}\tau_{n-1}^{j}\lm_{n}^{(k)}+(\bm_n-\lm_{n}^{(n-1)})\sum_{j=1}^{i}\tau_{n}^{j},
\end{eqnarray}
if $2\leq i\leq n-2$,
\begin{eqnarray*}
-\sum_{k=1}^{n-1}\tau_{n-1}\lm_{n}^{(k)}+(\bm_n-\lm_{n}^{(n-1)})\tau_{n},
\end{eqnarray*}
if $i=1$, and
\begin{eqnarray*}
-\sum_{j=1}^{n-2}\sum_{k=0}^{n-2}\tau_{n-1}^{j}\lm_{n}^{(k)}+(\bm_n-\lm_{n}^{(n-1)})\sum_{j=1}^{n-1}\tau_{n}^{j} -N_{n-1}\lm_{n}^{(n-1)},
\end{eqnarray*}
if $i=n-1$. We will prove that (\ref{astfinish}) vanishes and leave the second and the third cases to the reader:
\begin{eqnarray*}
-\sum_{j=1}^{i-1}\sum_{k=0}^{i-1}\tau_{n-1}^{j}\lm_{n}^{(k)}-\sum_{j=1}^{i}\sum_{k=i}^{n-1}\tau_{n-1}^{j}\lm_{n}^{(k)}+(\bm_n-\lm_{n}^{(n-1)})\sum_{j=1}^{i}\tau_{n}^{j}\\
=-\sum_{j=1}^{i-1}\sum_{k=0}^{i-1}\tau_{n-1}^{j-k}\lm_{n}^{(0)}\tau_{n}^{k}-\sum_{j=1}^{i}\sum_{k=i}^{n-1}\tau_{n-1}^{j-k}\lm_{n}^{(0)}\tau_{n}^{k}+\sum_{k=0}^{n-2}\sum_{j=1}^{i}\tau_{n-1}^{-k}\lm_{n}^{(0)}\tau_{n}^{k+j}\\
=-\sum_{j=1}^{i-1}\sum_{k=0}^{i-1}\tau_{n-1}^{j-k}\lm_{n}^{(0)}\tau_{n}^{k}-\sum_{j=1}^{i}\sum_{k=i}^{n-1}\tau_{n-1}^{j-k}\lm_{n}^{(0)}\tau_{n}^{k}+\sum_{j=1}^{i}\sum_{k=j}^{j+n-2}\tau_{n-1}^{j-k}\lm_{n}^{(0)}\tau_{n}^{k}\\
\end{eqnarray*}
The interval for $k$ in the third expression contains the interval for $k$ in the second one, so
\begin{eqnarray*}
-\sum_{j=1}^{i-1}\sum_{k=0}^{i-1}\tau_{n-1}^{j-k}\lm_{n}^{(0)}\tau_{n}^{k}-\sum_{j=1}^{i}\sum_{k=i}^{n-1}\tau_{n-1}^{j-k}\lm_{n}^{(0)}\tau_{n}^{k}+\sum_{j=1}^{i}\sum_{k=j}^{j+n-2}\tau_{n-1}^{j-k}\lm_{n}^{(0)}\tau_{n}^{k}\\
=-\sum_{j=1}^{i-1}\sum_{k=0}^{i-1}\tau_{n-1}^{j-k}\lm_{n}^{(0)}\tau_{n}^{k}+\sum_{j=1}^{i}\sum_{k=j}^{i-1}\tau_{n-1}^{j-k}\lm_{n}^{(0)}\tau_{n}^{k}+\sum_{j=1}^{i}\sum_{k=n}^{j+n-2}\tau_{n-1}^{j-k}\lm_{n}^{(0)}\tau_{n}^{k}\\
=-\sum_{j=1}^{i-1}\sum_{k=0}^{i-1}\tau_{n-1}^{j-k}\lm_{n}^{(0)}\tau_{n}^{k}+\sum_{j=1}^{i}\sum_{k=j}^{i-1}\tau_{n-1}^{j-k}\lm_{n}^{(0)}\tau_{n}^{k}+\sum_{j=1}^{i}\sum_{k=0}^{j-2}\tau_{n-1}^{j-k-n}\lm_{n}^{(0)}\tau_{n}^{k+n}\\
=-\sum_{j=1}^{i-1}\sum_{k=0}^{i-1}\tau_{n-1}^{j-k}\lm_{n}^{(0)}\tau_{n}^{k}+\sum_{j=1}^{i}\sum_{k=j}^{i-1}\tau_{n-1}^{j-k}\lm_{n}^{(0)}\tau_{n}^{k}+\sum_{j=1}^{i}\sum_{k=0}^{j-2}\tau_{n-1}^{j-k-1}\lm_{n}^{(0)}\tau_{n}^{k}\\
=-\sum_{j=1}^{i-1}\sum_{k=0}^{i-1}\tau_{n-1}^{j-k}\lm_{n}^{(0)}\tau_{n}^{k}+\sum_{j=1}^{i}\sum_{k=j}^{i-1}\tau_{n-1}^{j-k}\lm_{n}^{(0)}\tau_{n}^{k}+\sum_{j=0}^{i-1}\sum_{k=0}^{j-1}\tau_{n-1}^{j-k}\lm_{n}^{(0)}\tau_{n}^{k}
\end{eqnarray*}
\begin{eqnarray*}
=-\sum_{j=1}^{i-1}\sum_{k=0}^{i-1}\tau_{n-1}^{j-k}\lm_{n}^{(0)}\tau_{n}^{k}+\sum_{j=1}^{i-1}\sum_{k=j}^{i-1}\tau_{n-1}^{j-k}\lm_{n}^{(0)}\tau_{n}^{k}
+\sum_{j=1}^{i-1}\sum_{k=0}^{j-1}\tau_{n-1}^{j-k}\lm_{n}^{(0)}\tau_{n}^{k}\\
=-\sum_{j=1}^{i-1}\sum_{k=0}^{j-1}\tau_{n-1}^{j-k}\lm_{n}^{(0)}\tau_{n}^{k}+\sum_{j=1}^{i-1}\sum_{k=0}^{j-1}\tau_{n-1}^{j-k}\lm_{n}^{(0)}\tau_{n}^{k}
=0
\end{eqnarray*}

\medskip

\noindent {\it Part d)} is straightforward. \hfill $\Box$

\medskip

Lemma \ref{UV} implies
\begin{eqnarray*}
[\Ud,b^e]=0,\quad [\Vd+\Gamma,b^e]+[\Ud,B^e]=\frac{1}{2}b^e,\quad [\Vd+\Gamma,B^e]=-\frac{1}{2}B^e,
\end{eqnarray*}
which proves (\ref{u-conn}) for ${\nabla^{\circ}}$.

\subsection{Proof of Proposition \ref{nabla-nabla}}\label{proofnabla-nabla} Let us show that
\begin{eqnarray*}
2u^2(\nabla-\nabla^\circ)=(b^e+uB^e)\begin{pmatrix}\ld^{(0)}-\ld^{(0)}N & \ld^{(0)} \\0 & 0 \end{pmatrix}+\begin{pmatrix}\ld^{(0)}-\ld^{(0)}N & \ld^{(0)} \\0 & 0 \end{pmatrix}(b^e+uB^e)
\end{eqnarray*}
Using the matrix expressions for $b^e$ and $B^e$, the statement can be seen to follow from the following equalities:
\begin{eqnarray*}
\bm_{n+1} \ld_{n+1}^{(0)}+\ld_{n}^{(0)}\bm_{n+1}=-\lm_{n+1}^{(0)}\ld_{n+1}^{(1)}- \lm_{n+1}^{(n)}\ld_{n+1}^{(n)},
\end{eqnarray*}
\begin{eqnarray*}
\bm_{n+1} \ld_{n+1}^{(0)}N_{n+1}+\ld_{n}^{(0)}N_{n}\bm_{n+1}=0,
\end{eqnarray*}
\begin{eqnarray*}
\ld_{n}^{(0)}N_{n}\bmv_{n+1}h^e_n=0,
\end{eqnarray*}
\begin{eqnarray*}
h^e_{n+1}N_{n+1}\ld_{n+1}^{(0)}-h^e_{n+1}N_{n+1}\ld_{n+1}^{(0)}N_{n+1}=-h^e_{n+1}\sum_{i=1}^{n}\sum_{j=1}^{n+1}\tau_{n+1}^{j}\ld_{n+1}^{(i)}
\end{eqnarray*}
Since the last two equalities are quite straightforward, we will give proofs of the first two only. The first equality follows from Lemma \ref{lemmaDiMi} c):
\begin{eqnarray*}
\ld_{n}^{(0)}\bm_{n+1}=\ld_{n}^{(0)}\lm^{(0)}_{n+1}+\ld_{n}^{(0)}\lm^{(n)}_{n+1}+\ld_{n}^{(0)}\sum_{i=1}^{n-1}\lm^{(i)}_{n+1}\\
=-\lm^{(0)}_{n+1}\ld_{n+1}^{(1)}-\lm^{(0)}_{n+1}\ld_{n+1}^{(0)}-\lm^{(n)}_{n+1}\ld_{n+1}^{(0)}-\lm^{(n)}_{n+1}\ld_{n+1}^{(n)}-\sum_{i=1}^{n-1}\lm^{(i)}_{n+1}\ld_{n+1}^{(0)}\\
=-\lm^{(0)}_{n+1}\ld_{n+1}^{(1)}-\lm^{(n)}_{n+1}\ld_{n+1}^{(n)}-\sum_{i=0}^{n}\lm^{(i)}_{n+1}\ld_{n+1}^{(0)}=-\lm^{(0)}_{n+1}\ld_{n+1}^{(1)}-\lm^{(n)}_{n+1}\ld_{n+1}^{(n)}-\bm_{n+1} \ld_{n+1}^{(0)}
\end{eqnarray*}
The proof of the second equality uses Lemmas \ref{Mitau} f), \ref{lemmaDiMi} c) and the computation we just did:
\begin{eqnarray*}
\bm_{n+1} \ld_{n+1}^{(0)}N_{n+1}+\ld_{n}^{(0)}N_{n}\bm_{n+1}=\bm_{n+1} \ld_{n+1}^{(0)}N_{n+1}+\ld_{n}^{(0)}(\bm_{n+1}-\lm^{(n)}_{n+1})N_{n+1}\\
=(\bm_{n+1} \ld_{n+1}^{(0)}+\ld_{n}^{(0)}\bm_{n+1})N_{n+1}-\ld_{n}^{(0)}\lm^{(n)}_{n+1}N_{n+1}\\
=(-\lm^{(0)}_{n+1}\ld_{n+1}^{(1)}-\lm^{(n)}_{n+1}\ld_{n+1}^{(n)})N_{n+1}-\ld_{n}^{(0)}\lm^{(n)}_{n+1}N_{n+1}\\=(-\lm^{(0)}_{n+1}\ld_{n+1}^{(1)}-\lm^{(n)}_{n+1}\ld_{n+1}^{(n)}+\lm^{(n)}_{n+1}\ld_{n+1}^{(0)}+\lm^{(n)}_{n+1}\ld_{n+1}^{(n)})N_{n+1}\\
=(-\lm^{(0)}_{n+1}\ld_{n+1}^{(1)}+\lm^{(n)}_{n+1}\ld_{n+1}^{(0)})N_{n+1}=(-\lm^{(0)}_{n+1}\tau_{n+1}^{-1}\ld_{n+1}^{(0)}\tau_{n+1}+\lm^{(0)}_{n+1}\tau_{n+1}^{-1}\ld_{n+1}^{(0)})N_{n+1}=0
\end{eqnarray*}

\medskip

\subsection{Proof of Corollary \ref{alterconncycl}}\label{proofalterconncycl}
Let us show that
\begin{eqnarray*}
2u^2(\widetilde{\nabla}-\nabla)=(b^e+uB^e)H(u)+H(u)(b^e+uB^e)
\end{eqnarray*}
where $H(u)=H_0+uH_1: {\rC^e}(A)\to {\rC^e}(A)$ and
$$
H_0=\begin{pmatrix}\lm^{(0)} & \lm^{(0)} \\0 & 0 \end{pmatrix},\quad H_1=\begin{pmatrix} 0 & 0 \\ h^eN' & 0 \end{pmatrix}, \quad N'_{n+1}=\sum_{j=1}^{n+1}j\tau_{n+1}^{j}
$$

We will use the following property of $N'$:
$$
(1-\tau^{-1})N'=-N-2\Gamma+1
$$

Let us compute $(\bed+\bem+uB^e)H_0+H_0(\bed+\bem+uB^e)$ first. By Lemma \ref{lemmabdMi} b), $\bed H_0+H_0\bed=0$.
By Lemma \ref{lemmaDiMi} b)
\begin{eqnarray*}
\bem H_0+H_0\bem=\begin{pmatrix}\bm\lm^{(0)}+\lm^{(0)}\bm & \bm\lm^{(0)}+\lm^{(0)}\bm \\0 & 0 \end{pmatrix}
=2\Umc
\end{eqnarray*}
Finally,
\begin{eqnarray*}
B^eH_0+H_0B^e=\begin{pmatrix} N & 0 \\h^eN\lm^{(0)} & h^eN\lm^{(0)} \end{pmatrix}=:T\\
\end{eqnarray*}
Thus,
\begin{eqnarray}\label{I}
(\bed+\bem+uB^e)H_0+H_0(\bed+\bem+uB^e)=2\Umc+uT
\end{eqnarray}

Next, let us compute $(\bed+\bem+uB^e)H_1+H_1(\bed+\bem+uB^e)$. Note that $\bed H_1+H_1\bed=0$ by Lemmas \ref{Mitau} c) and \ref{lemmabh} a). Also, $B^eH_1+H_1B^e=0$ for obvious reasons. Furthermore, by Lemma \ref{lemmabh} b), c)
\begin{eqnarray*}
\bem H_1+H_1\bem=\begin{pmatrix} \bmv h^eN' & 0 \\ \bmh h^eN' + h^eN' \bm & h^eN'\bmv \end{pmatrix}\\
=\begin{pmatrix} (1-\tau^{-1})N' & 0 \\ -h^e(\lm^{(0)}+\bmh)N' + h^eN' \bm & h^eN'\bmv \end{pmatrix}\\
=\begin{pmatrix} -N-2\Gamma+id & 0 \\ -h^e(\lm^{(0)}+\bmh)N' + h^eN' \bm & h^eN'\bmv \end{pmatrix}
\end{eqnarray*}
\begin{eqnarray}\label{II}
=\begin{pmatrix} -N-2\Gamma+id & 0 \\ -h^e\lm^{(0)}N'-h^e\bmh N' + h^eN' \bm & -h^eN\lm^{(0)}-2\Gamma \end{pmatrix}
\end{eqnarray}
To explain the last transition above, we'll take $\alpha\in\rC(A)$ and notice that
$$
h^eN'\bmv h^e(\alpha)=h^eN'(1-\tau^{-1})(\alpha)=-h^e(N+2\Gamma-1)(\alpha)=-h^eN\lm^{(0)}(h^e(\alpha))-2\Gamma(h^e(\alpha))
$$
Let us simplify the remaining entry of the matrix (\ref{II}):
\begin{eqnarray*}
N'_{n} \bm_{n+1}-\lm^{(0)}_{n+1}N'_{n+1}-\bmh_{n+1} N'_{n+1}=\sum_{i=0}^{n}\sum_{j=1}^{n}j\tau_{n}^{j} \lm^{(i)}_{n+1}-\sum_{j=1}^{n+1}\sum_{i=0}^{n-1}j \lm^{(i)}_{n+1}\tau_{n+1}^{j}\\
= \sum_{i=0}^{n}\sum_{j=1}^{n}j\tau_{n}^{j} \lm^{(i)}_{n+1}-\sum_{j=1}^{n+1}\sum_{i=0}^{n-1}j \tau_{n}^{-i}\lm^{(0)}_{n+1}\tau_{n+1}^{i+j}\\
= \sum_{i=0}^{n}\sum_{j=1}^{n}j\tau_{n}^{j} \lm^{(i)}_{n+1}-\sum_{j=1}^{n}\sum_{i=0}^{n-j}j \tau_{n}^{j}\lm^{(i+j)}_{n+1}-\sum_{j=2}^{n+1}\sum_{i=n-j+1}^{n-1}j \tau_{n}^{j-1}\lm^{(i+j-n-1)}_{n+1}
\end{eqnarray*}
By reindexing and resumming one gets
\begin{eqnarray*}
\sum_{i=0}^{n}\sum_{j=1}^{n}j\tau_{n}^{j} \lm^{(i)}_{n+1}-\sum_{i=1}^{n}\sum_{j=1}^{i}j \tau_{n}^{j}\lm^{(i)}_{n+1}-\sum_{i=0}^{n-1}\sum_{j=i+1}^{n}(j+1) \tau_{n}^{j}\lm^{(i)}_{n+1}\\=-\sum_{i=0}^{n-1}\sum_{j=i+1}^{n}\tau_{n}^{j}\lm^{(i)}_{n+1}=-N_{n}\lm^{(0)}_{n+1}-\sum_{i=1}^{n-1}\sum_{j=i+1}^{n}\tau_{n}^{j}\lm^{(i)}_{n+1}
\end{eqnarray*}
Inserting this formula into (\ref{II}) gives
\begin{eqnarray}\label{IV}
\bem H_1+H_1\bem=-T-2\Gamma+2\Vmc
\end{eqnarray}

Thus, by (\ref{I}) and (\ref{IV})
\begin{eqnarray*}
(\bed+\bem+uB^e)\frac{H(u)}{2u^2}+\frac{H(u)}{2u^2}(\bed+\bem+uB^e)=\frac{\Umc}{u^2}+\frac{\Vmc-\Gamma}{u}
\end{eqnarray*}
which finishes the proof.

\medskip

\subsection{Proof of Corollary \ref{nonext}}\label{proofnonext}

The proof is based on formula (\ref{gaugeconn}) applied to (\ref{ext}), with $\iota(u)$, $p(u)$, and $H(u)$ given by (\ref{iota}),(\ref{p}), and (\ref{piota}), respectively. Let us compute all the summands in (\ref{gaugeconn}).

First, observe that
\begin{eqnarray*}
p(u)\frac{d\iota(u)}{du}=\begin{pmatrix}1 & (1-\tau^{-1})h\lm^{(0)}\end{pmatrix}\begin{pmatrix}0 \\ h^ehN  \end{pmatrix}=(1-\tau^{-1})hhN\sim0
\end{eqnarray*}
since
$$
(1-\tau^{-1})hhN=(b+uB)\frac{H(u)}{u}+\frac{H(u)}{u}(b+uB)
$$
where $H(u)=u(1-\tau^{-1})hhhN$ (see the proof of Part a) of Lemma \ref{iotap}).
Furthermore, by Lemmas \ref{lemmaDiMiH} a),b)
\begin{eqnarray*}
p(u)\Udl\iota(u)=\frac12\begin{pmatrix}1 & (1-\tau^{-1})h\lm^{(0)}  \end{pmatrix}\begin{pmatrix}-\lm^{(0)}\ld^{(1)} & -\lm^{(0)}\ld^{(1)} \\0 & 0 \end{pmatrix}\begin{pmatrix}1 \\ uh^ehN  \end{pmatrix}=-\frac12\lm^{(0)}\ld^{(1)}=\nUd
\end{eqnarray*}
and similarly $p(u)\Vdl\iota(u)=\nVd$. It remains to compute $p(u)\Gamma\iota(u)$:
\begin{eqnarray*}
p(u)\Gamma\iota(u)=\Gamma+ u(1-\tau^{-1})h\lm^{(0)}\Gamma h^ehN= \Gamma+ u(1-\tau^{-1})h\lm^{(0)}(h^ehN\Gamma-h^ehN)\\
= \Gamma+ u(1-\tau^{-1})hhN\Gamma-u(1-\tau^{-1})hhN
\sim \Gamma+ ((b+uB) H(u)+H(u)(b+uB))\Gamma
\end{eqnarray*}
By Part d) of Lemma \ref{UV}, the latter is $u$-homotopic to
\begin{eqnarray*}
\Gamma-\frac12 H(u)(\bm-uB)
\end{eqnarray*}
To summarize,
\begin{eqnarray*}
p(u)\frac{d\iota(u)}{du}+p(u)\frac{\Udl}{u^2}\iota(u)+p(u)\frac{\Vdl}{u}\iota(u)+p(u)\frac{\Gamma}{u}\iota(u)+\frac{1}{2u}H(u)(b-uB)\\
\sim \frac{\nUd}{u^2}+\frac{\nVd}{u}+\frac{\Gamma}{u}-\frac1{2u} H(u)(\bm-uB)+\frac{1}{2u}H(u)(b-uB)\\=\frac{\nUd}{u^2}+\frac{\nVd}{u}+\frac{\Gamma}{u}+\frac{1}{2}(1-\tau^{-1})hhhN\bd
\end{eqnarray*}

\medskip

\subsection{Proof of Proposition \ref{zgr}}\label{proofzgr}
Let us start by introducing the following operators on $\rC^e(A)$:
\begin{eqnarray*}
\gamma^{(i)}_{n+1}(a_0[a_1|\ldots |a_n])=\deg(a_i)a_0[a_1|\ldots|a_n]
\end{eqnarray*}
They extend to $\rC(A^e)$ and the extended operators satisfy the relations
$$
\gamma^{(i)}_{n+1}=\tau_{n+1}^{-i}\gamma^{(0)}_{n+1}\tau_{n+1}^{i}
$$
Obviously,
$$
\Gamma'_{n+1}=\sum_{i=0}^n\gamma^{(i)}_{n+1}+2\Gamma_{n+1}
$$

\medskip

\begin{lemma}\label{gammaDiMi}
\begin{eqnarray*}
a)\quad \gamma^{(i)}\ld^{(j)}=\begin{cases} \ld^{(i)}\gamma^{(i)}+\ld^{(i)} & i=j\\ \ld^{(j)}\gamma^{(i)} & i\neq j \end{cases}
\end{eqnarray*}
\begin{eqnarray*}
b)\quad \gamma^{(i)}_n\lm^{(j)}_{n+1}=\begin{cases} \lm^{(i)}_{n+1}\gamma^{(i)}_{n+1}+\lm^{(i)}_{n+1}\gamma^{(i+1)}_{n+1} & i=j\\
\lm^{(j)}_{n+1}\gamma^{(i)}_{n+1} & i<j\\ \lm^{(j)}_{n+1}\gamma^{(i+1)}_{n+1} & i>j
\end{cases}
\end{eqnarray*}
\begin{eqnarray*}
c)\quad \gamma^{(i)}h^e=\begin{cases} 0 & i=0\\ h^e\gamma^{(i-1)} & i\geq 1 \end{cases}
\end{eqnarray*}
\end{lemma}
\noindent{\bf Proof.} Part c) is obvious. The rest is very similar to the proof of Lemma \ref{lemmaDiMi}: the case $i=0$ is easy; Lemma \ref{Mitau} a),b) reduces the general case to this special one. \hfill $\Box$

\medskip
Let us show that
\begin{eqnarray*}
2u^2(\nabla_{gr}-\nabla)=(b^e+uB^e)\hH(u)+\hH(u)(b^e+uB^e)
\end{eqnarray*}
where $\hH(u)=\hH_0+u\hH_1: {\rC^e}(A)\to {\rC^e}(A)$ and
$$
(\hH_0)_{n+1}=\begin{pmatrix}\lm^{(0)}_{n+1}\gamma^{(1)}_{n+1} & \lm^{(0)}_{n+1}\gamma^{(1)}_{n+1} \\0 & 0 \end{pmatrix},\quad (\hH_1)_{n+1}=\begin{pmatrix} 0 & 0 \\ h^e_{n+1}\sum_{i=1}^{n}\sum_{j=i+1}^{n+1}\tau_{n+1}^{j}\gamma_{n+1}^{(i)} & 0 \end{pmatrix}
$$

Let us compute $(\bed+\bem+uB^e)\hH_0+\hH_0(\bed+\bem+uB^e)$ first.
By Lemmas \ref{lemmabdMi} b), \ref{gammaDiMi} a)
\begin{eqnarray*}
\bed \hH_0+\hH_0\bed=-2\Udl
\end{eqnarray*}
By Lemma \ref{gammaDiMi} b)
\begin{eqnarray*}
\bem \hH_0+\hH_0\bem=\begin{pmatrix}\ast & \ast\\0 & 0 \end{pmatrix}
\end{eqnarray*}
where $$\ast=\bm\lm^{(0)}\gamma^{(1)}+\lm^{(0)}(\bm\gamma^{(1)}+\lm^{(0)}\gamma^{(2)}+\lm^{(1)}\gamma^{(2)}-\lm^{(0)}\gamma^{(1)})$$
Using Lemma \ref{lemmaDiMi} b)
\begin{eqnarray*}
\bm\lm^{(0)}\gamma^{(1)}+\lm^{(0)}(\bm\gamma^{(1)}+\lm^{(0)}\gamma^{(2)}+\lm^{(1)}\gamma^{(2)}-\lm^{(0)}\gamma^{(1)})\\
=(\bm\lm^{(0)}+\lm^{(0)}\bm)\gamma^{(1)}+(\lm^{(0)}\lm^{(0)}+\lm^{(0)}\lm^{(1)})\gamma^{(2)}-\lm^{(0)}\lm^{(0)}\gamma^{(1)}\\
=\lm^{(0)}\lm^{(0)}\gamma^{(1)}-\lm^{(0)}\lm^{(0)}\gamma^{(1)}=0
\end{eqnarray*}
i.e.
$
\bem \hH_0+\hH_0\bem=0
$.
Also, by Lemmas \ref{gammaDiMi} c),  \ref{lemmaDiMiH} b)
\begin{eqnarray*}\label{4a}
B^e\hH_0+\hH_0B^e=\begin{pmatrix} \gamma^{(0)}N & 0 \\h^eN\lm^{(0)}\gamma^{(1)} & h^eN\lm^{(0)}\gamma^{(1)} \end{pmatrix}=:T'\\
\end{eqnarray*}
Thus,
\begin{eqnarray}\label{5a}
(\bed+\bem+uB^e)\hH_0+\hH_0(\bed+\bem+uB^e)=-2\Udl+uT'
\end{eqnarray}

Let us compute $(\bed+\bem+uB^e)\hH_1+\hH_1(\bed+\bem+uB^e)$ now. First, by Lemmas \ref{Mitau} c), \ref{lemmabh} a), \ref{gammaDiMi} a)
\begin{eqnarray}\label{aa}
\bed \hH_1+\hH_1\bed=-2\Vdl
\end{eqnarray}
Also, $B^e\hH_1+\hH_1B^e=0$. It remains to simplify $\bem \hH_1+\hH_1\bem$.

\begin{eqnarray*}
\bem_{n+1} (\hH_1)_n+(\hH_1)_{n-1}\bem_{n}=\begin{pmatrix} \ast_{11} & 0 \\ \ast_{21} & \ast_{22} \end{pmatrix}\\
\end{eqnarray*}
where
$$
\ast_{11}=\bmv_{n+1}h^e_{n}\sum_{i=1}^{n-1}\sum_{j=i+1}^{n}\tau_{n}^{j}\gamma_{n}^{(i)},
$$
$$
\ast_{21}=\bmh_{n+1}h^e_{n}\sum_{i=1}^{n-1}\sum_{j=i+1}^{n}\tau_{n}^{j}\gamma_{n}^{(i)}+h^e_{n-1}\sum_{i=1}^{n-2}\sum_{j=i+1}^{n-1}\tau_{n-1}^{j}\gamma_{n-1}^{(i)}\bm_{n},
$$
$$
\ast_{22}=h^e_{n-1}\sum_{i=1}^{n-2}\sum_{j=i+1}^{n-1}\tau_{n-1}^{j}\gamma_{n-1}^{(i)}\bmv_{n}
$$
%%%%%%%%%%%%%%%%%%%%%%%%%%%%%%%%%%%%%%%%%%%%%%%%%%%
%%%%%%%%%%%%%%%%%%%%%%%%%%%%%%%%%%%%%%%%%%%%%%%%%%
By Lemma \ref{lemmabh},
\begin{eqnarray*}
\ast_{11}=(1-\tau^{-1}_n)\sum_{i=1}^{n-1}\sum_{j=i+1}^{n}\tau_{n}^{j}\gamma_{n}^{(i)}
=-\sum_{i=1}^{n-1}(\tau_{n}^{i}-1)\gamma_{n}^{(i)}
=-\sum_{i=1}^{n-1}\tau_{n}^{i}\gamma_{n}^{(i)}+\sum_{i=1}^{n-1}\gamma_{n}^{(i)}\\
=-\sum_{i=0}^{n-1}\tau_{n}^{i}\gamma_{n}^{(i)}+\sum_{i=0}^{n-1}\gamma_{n}^{(i)}
=-\sum_{i=0}^{n-1}\gamma_{n}^{(0)}\tau_{n}^{i}+\Gamma'_n-2\Gamma_n
\end{eqnarray*}
\begin{eqnarray}\label{3a}
=-\gamma_{n}^{(0)}N_{n}+\Gamma'_n-2\Gamma_n
\end{eqnarray}

To simplify $\ast_{22}$, take $\alpha\in A\otimes sA^{\otimes (n-2)}$ and apply $\ast_{22}$ to $h^e_{n-1}(\alpha)$:
\begin{eqnarray*}
h^e_{n-1}\sum_{i=1}^{n-2}\sum_{j=i+1}^{n-1}\tau_{n-1}^{j}\gamma_{n-1}^{(i)}\bmv_{n}h^e_{n-1}(\alpha)
=h^e_{n-1}\sum_{i=1}^{n-2}\sum_{j=i+1}^{n-1}\tau_{n-1}^{j}\gamma_{n-1}^{(i)}(1-\tau_{n-1}^{-1})(\alpha)\\
=h^e_{n-1}\sum_{i=1}^{n-2}\sum_{j=i+1}^{n-1}\tau_{n-1}^{j}\gamma_{n-1}^{(i)}(\alpha)-h^e_{n-1}\sum_{i=1}^{n-2}\sum_{j=i+1}^{n-1}\tau_{n-1}^{j-1}\gamma_{n-1}^{(i-1)}(\alpha)\\
=h^e_{n-1}\sum_{i=1}^{n-2}\sum_{j=i+1}^{n-1}\tau_{n-1}^{j}\gamma_{n-1}^{(i)}(\alpha)-h^e_{n-1}\sum_{i=0}^{n-3}\sum_{j=i+1}^{n-2}\tau_{n-1}^{j}\gamma_{n-1}^{(i)}(\alpha)\\
=h^e_{n-1}\gamma_{n-1}^{(n-2)}(\alpha)+h^e_{n-1}\sum_{i=1}^{n-3}\sum_{j=i+1}^{n-1}\tau_{n-1}^{j}\gamma_{n-1}^{(i)}(\alpha)\\
-h^e_{n-1}\sum_{i=1}^{n-3}\sum_{j=i+1}^{n-2}\tau_{n-1}^{j}\gamma_{n-1}^{(i)}(\alpha)-h^e_{n-1}\sum_{j=1}^{n-2}\tau_{n-1}^{j}\gamma_{n-1}^{(0)}(\alpha)\\
=h^e_{n-1}\gamma_{n-1}^{(n-2)}(\alpha)+h^e_{n-1}\sum_{i=1}^{n-3}\gamma_{n-1}^{(i)}(\alpha)
-h^e_{n-1}\sum_{j=1}^{n-2}\tau_{n-1}^{j}\gamma_{n-1}^{(0)}(\alpha)\\
=h^e_{n-1}\sum_{i=0}^{n-2}\gamma_{n-1}^{(i)}(\alpha)-h^e_{n-1}N_{n-1}\gamma_{n-1}^{(0)}(\alpha)
=\sum_{i=0}^{n-1}\gamma_{n}^{(i)}h^e_{n-1}(\alpha)-h^e_{n-1}N_{n-1}\gamma_{n-1}^{(0)}\lm_{n}^{(0)}h^e_{n-1}(\alpha)\\
=(\Gamma'_n-2\Gamma_n)h^e_{n-1}(\alpha)-h^e_{n-1}N_{n-1}\lm_{n}^{(0)}\gamma_{n-1}^{(1)}h^e_{n-1}(\alpha)
\end{eqnarray*}
Thus,
\begin{eqnarray}\label{6a}
\ast_{22}=\Gamma'_n-2\Gamma_n-h^e_{n-1}N_{n-1}\lm_{n}^{(0)}\gamma_{n-1}^{(1)}
\end{eqnarray}

Furthermore, by Lemmas \ref{lemmabh} c), \ref{Mitau} b), \ref{gammaDiMi} b)
\begin{eqnarray*}
\ast_{21}=-h^e_{n-1}\sum_{i=1}^{n-1}\sum_{j=i+1}^{n}\sum_{k=0}^{n-2}\lm^{(k)}_n\tau_{n}^{j}\gamma_{n}^{(i)}+h^e_{n-1}\sum_{i=1}^{n-2}\sum_{j=i+1}^{n-1}\sum_{k=0}^{n-1}\tau_{n-1}^{j}\gamma_{n-1}^{(i)}\lm^{(k)}_n\\
=-h^e_{n-1}\sum_{i=1}^{n-2}\sum_{j=i+1}^{n-1}\sum_{k=0}^{n-j-1}\tau_{n-1}^{j}\lm_{n}^{(k+j)}\gamma_{n}^{(i)}
-h^e_{n-1}\sum_{i=1}^{n-1}\sum_{j=i+1}^{n}\sum_{k=n-j}^{n-2}\tau_{n-1}^{j-1}\lm_{n}^{(k+j-n)}\gamma_{n}^{(i)}\\
+h^e_{n-1}\sum_{i=1}^{n-2}\sum_{j=i+1}^{n-1}\tau_{n-1}^{j}(\lm^{(i)}_{n}\gamma^{(i)}_{n}+\lm^{(i)}_{n}\gamma^{(i+1)}_{n})
+h^e_{n-1}\sum_{i=1}^{n-2}\sum_{j=i+1}^{n-1}\sum_{k=i+1}^{n-1}\tau_{n-1}^{j}\lm^{(k)}_{n}\gamma^{(i)}_{n}\\
+h^e_{n-1}\sum_{i=1}^{n-2}\sum_{j=i+1}^{n-1}\sum_{k=0}^{i-1}\tau_{n-1}^{j}\lm^{(k)}_{n}\gamma^{(i+1)}_{n}\\
\end{eqnarray*}
or, after reindexing
\begin{eqnarray*}
-h^e_{n-1}\sum_{i=1}^{n-2}\sum_{j=i+1}^{n-1}\sum_{k=j}^{n-1}\tau_{n-1}^{j}\lm_{n}^{(k)}\gamma_{n}^{(i)}
-h^e_{n-1}\sum_{i=1}^{n-1}\sum_{j=i}^{n-1}\sum_{k=0}^{j-1}\tau_{n-1}^{j}\lm_{n}^{(k)}\gamma_{n}^{(i)}\\
+h^e_{n-1}\sum_{i=1}^{n-2}\sum_{j=i+1}^{n-1}\tau_{n-1}^{j}\lm^{(i)}_{n}\gamma^{(i)}_{n}
+h^e_{n-1}\sum_{i=2}^{n-1}\sum_{j=i}^{n-1}\tau_{n-1}^{j}\lm^{(i-1)}_{n}\gamma^{(i)}_{n}\\
+h^e_{n-1}\sum_{i=1}^{n-2}\sum_{j=i+1}^{n-1}\sum_{k=i+1}^{n-1}\tau_{n-1}^{j}\lm^{(k)}_{n}\gamma^{(i)}_{n}
+h^e_{n-1}\sum_{i=2}^{n-1}\sum_{j=i}^{n-1}\sum_{k=0}^{i-2}\tau_{n-1}^{j}\lm^{(k)}_{n}\gamma^{(i)}_{n}\\
=-h^e_{n-1}\sum_{i=1}^{n-2}\sum_{j=i+1}^{n-1}\sum_{k=j}^{n-1}\tau_{n-1}^{j}\lm_{n}^{(k)}\gamma_{n}^{(i)}
-h^e_{n-1}\sum_{i=1}^{n-1}\sum_{j=i}^{n-1}\sum_{k=0}^{j-1}\tau_{n-1}^{j}\lm_{n}^{(k)}\gamma_{n}^{(i)}\\
+h^e_{n-1}\sum_{i=1}^{n-2}\sum_{j=i+1}^{n-1}\sum_{k=i}^{n-1}\tau_{n-1}^{j}\lm^{(k)}_{n}\gamma^{(i)}_{n}
+h^e_{n-1}\sum_{i=2}^{n-1}\sum_{j=i}^{n-1}\sum_{k=0}^{i-1}\tau_{n-1}^{j}\lm^{(k)}_{n}\gamma^{(i)}_{n}\\
=-h^e_{n-1}\sum_{j=1}^{n-1}\tau_{n-1}^{j}\lm^{(0)}_{n}\gamma^{(1)}_{n}=-h^e_{n-1}N_{n-1}\lm^{(0)}_{n}\gamma^{(1)}_{n}
\end{eqnarray*}
%%%%%%%%%%%%%%%%%%%%%%%%%%%%
This, together with (\ref{5a}), (\ref{aa}), (\ref{3a}), and (\ref{6a}) shows that
\begin{eqnarray*}
(\bed+\bem+uB^e)\frac{\hH(u)}{2u^2}+\frac{\hH(u)}{2u^2}(\bed+\bem+uB^e)=\frac{\Gamma'}{2u}-\frac{\Udl}{u^2}-\frac{\Vdl+\Gamma}{u}
\end{eqnarray*}

%%%%%%%%%%%%%%%%%%%%%%%%%%%%%%%%%%%%%%%%%%%%%%%%%%%

\medskip
\section{Proofs for Section \ref{locsing}}
\subsection{Auxiliary lemma}\label{auxlem}

Let us derive some basic properties of operators similar to $\bw$ and $\bdw$.

Let $(A,d)$ be an (abstract) dg algebra and $a\in A$. Consider the following operator on $\rC^e(A)$:
$$
\ID_{n+1}=\sum_{i=1}^{n+1}a^{(i)}_{n+1}
$$
where
\begin{eqnarray}\label{ID}
a^{(i)}_{n+1}(a_0[a_1|\ldots |a_n])=(-1)^{|sa|\sum_{j=0}^{i-1}|sa_j|}a_0[a_1|\ldots|a_{i-1}|a|a_{i}| \ldots|a_n]
\end{eqnarray}
Note that the parity of $\ID$ is opposite to that of $a$. Clearly, $\ID$ preserves the subspaces $\rC(A)$ and $\rC^+(A)$, and therefore we may also think of it as represented by the matrix
$$
\begin{pmatrix} b(a) & 0 \\ 0 & b(a) \end{pmatrix},
$$
where the two $b(a)$'s are the restrictions of $\ID$ onto $\rC(A)$ and $\rC^+(A)$, respectively. It is also obvious that $\ID$ extends to ${\rrCe}(A)$.

We will need yet another operator $\dD$ on $\rC^e(A)$ defined as follows
$$
\dD_{n+1}:=\sum_{i=0}^{n+1}\dD^{(i)}_{n+1}
$$
where
$$
\dD^{(i)}_{n+1}(a_0[a_1|\ldots |a_n])=(-1)^{|a|\sum_{j=0}^{i-1}|sa_j|}a_0[a_1|\ldots|[a,a_i]|\ldots|a_n]
$$
($[\,,\,]$ stands for the super-commutator).
\medskip

\begin{lemma}\label{intertw}
\begin{eqnarray*}
a)\quad [\bed, \ID]=\IdD
\end{eqnarray*}
\begin{eqnarray*}
b)\quad [\bem, \ID]=(-1)^{|a|}\dD
\end{eqnarray*}
\begin{eqnarray*}
c)\quad [B^e, \ID]=0
\end{eqnarray*}
\begin{eqnarray*}
d)\quad [\JD, \ID]=0,\quad \forall a, a'
\end{eqnarray*}
\end{lemma}
\noindent{\bf Proof.} As in section \ref{auxlemmas}, it will be convenient to extend all the operators to $\rC(A^e)$ and establish the relations on the whole of $\rC(A^e)$.

Observe that the extensions of $a^{(i)}_{n+1}$ and $\dD^{(i)}_{n+1}$ to $\rC(A^e)$ satisfy the properties
$$
a^{(i)}_{n+1}=\tau_{n+2}^{-i}a^{(0)}_{n+1}\tau_{n+1}^{i},\quad \dD^{(i)}_{n+1}=\tau_{n+1}^{-i}\dD^{(0)}_{n+1}\tau_{n+1}^{i}
$$
where
$
a^{(0)}_{n+1}(a_0[a_1|\ldots |a_n])=a[a_0|a_1|\ldots |a_n].
$

a) By (\ref{Di})
\begin{eqnarray*}\ld^{(i)}_{n+2}a^{(0)}_{n+1}=
\begin{cases}
da^{(0)}_{n+1},\quad i=0\\
(-1)^{|a|+1}a^{(0)}_{n+1}\ld^{(i-1)}_{n+1},\quad i\neq0\\
\end{cases}
\end{eqnarray*}
Therefore, by Lemma \ref{Mitau} a)
\begin{eqnarray*}\ld^{(i)}_{n+2}a^{(j)}_{n+1}=\ld^{(i)}_{n+2}\tau_{n+2}^{-j}a^{(0)}_{n+1}\tau_{n+1}^{j}
=\tau_{n+2}^{-j}\ld^{(i-j)}_{n+2}a^{(0)}_{n+1}\tau_{n+1}^{j}\\
=\begin{cases}
da^{(j)}_{n+1},\quad i=j\\
(-1)^{|a|+1}\tau_{n+2}^{-j}a^{(0)}_{n+1}\ld^{(i-j-1)}_{n+1}\tau_{n+1}^{j},\quad i\neq j\\
\end{cases}\\
=\begin{cases}
da^{(j)}_{n+1},\quad i=j\\
(-1)^{|a|+1}a^{(j)}_{n+1}\ld^{(i-1)}_{n+1},\quad i>j\\
(-1)^{|a|+1}a^{(j)}_{n+1}\ld^{(i)}_{n+1},\quad i<j\\
\end{cases}
\end{eqnarray*}
Then
\begin{eqnarray*}
\bed_{n+2}\ID_{n+1}=\sum_{i=0}^{n+1}\sum_{j=1}^{n+1}\ld^{(i)}_{n+2}a^{(j)}_{n+1}\\
=\sum_{i=1}^{n+1}\ld^{(i)}_{n+2}a^{(i)}_{n+1}+\sum_{i=0}^{n}\sum_{j=i+1}^{n+1}\ld^{(i)}_{n+2}a^{(j)}_{n+1}+\sum_{i=2}^{n+1}\sum_{j=1}^{i-1}\ld^{(i)}_{n+2}a^{(j)}_{n+1}\\
=\IdD_{n+1}+(-1)^{|a|+1}\sum_{i=0}^{n}\sum_{j=i+1}^{n+1}a^{(j)}_{n+1}\ld^{(i)}_{n+1}+(-1)^{|a|+1}\sum_{i=2}^{n+1}\sum_{j=1}^{i-1}
a^{(j)}_{n+1}\ld^{(i-1)}_{n+1}\\
=\IdD_{n+1}+(-1)^{|a|+1}\sum_{i=0}^{n}\sum_{j=i+1}^{n+1}a^{(j)}_{n+1}\ld^{(i)}_{n+1}+(-1)^{|a|+1}\sum_{i=1}^{n}\sum_{j=1}^{i}
a^{(j)}_{n+1}\ld^{(i)}_{n+1}\\
=\IdD_{n+1}+(-1)^{|a|+1}\ID_{n+1}\bed_{n+1}
\end{eqnarray*}

\medskip

b) Set

$$
\lD^{(0)}_{n+1}(a_0[a_1|\ldots |a_n])=aa_0[a_1|\ldots |a_n], \quad \lD^{(i)}_{n+1}:=\tau_{n+1}^{-i}\lD^{(0)}_{n+1}\tau_{n+1}^{i}
$$
$$
\rD^{(0)}_{n+1}(a_0[a_1|\ldots |a_n])=(-1)^{|a||a_0|}a_0a[a_1|\ldots |a_n], \quad \rD^{(i)}_{n+1}:=\tau_{n+1}^{-i}\rD^{(0)}_{n+1}\tau_{n+1}^{i}
$$
Then, clearly,
$
\dD^{(i)}_{n+1}=\lD^{(i)}_{n+1}-\rD^{(i)}_{n+1}.
$
The following is easy to check using (\ref{Mi})
\begin{eqnarray*}\lm^{(i)}_{n+2}a^{(0)}_{n+1}=
\begin{cases}
(-1)^{|a|}\lD^{(0)}_{n+1},\quad i=0\\
(-1)^{|a|+1}a^{(0)}_{n}\lm^{(i-1)}_{n+1},\quad 1\leq i\leq n\\
(-1)^{|a|+1}\rD^{(0)}_{n+1}\tau_{n+1}^{-1},\quad i=n+1
\end{cases}
\end{eqnarray*}
In combination with Lemma \ref{Mitau} b) this gives
\begin{eqnarray*}\lm^{(i)}_{n+2}a^{(j)}_{n+1}=\lm^{(i)}_{n+2}\tau_{n+2}^{-j}a^{(0)}_{n+1}\tau_{n+1}^{j}=\lm^{(i)}_{n+2}\tau_{n+2}^{n+2-j}a^{(0)}_{n+1}\tau_{n+1}^{j}\\
=\begin{cases} \tau_{n+1}^{1-j}\lm_{n+2}^{(n+2+i-j)}a^{(0)}_{n+1}\tau_{n+1}^{j} & i-j<0\\ \tau_{n+1}^{-j}\lm_{n+2}^{(i-j)}a^{(0)}_{n+1}\tau_{n+1}^{j} & i-j\geq 0\end{cases}
\end{eqnarray*}
\begin{eqnarray*}
=\begin{cases}
(-1)^{|a|+1}\tau_{n+1}^{1-j}a^{(0)}_{n}\lm_{n+1}^{(n+1+i-j)}\tau_{n+1}^{j} & -n-1\leq i-j\leq -2\\
(-1)^{|a|+1}\tau_{n+1}^{1-j}\rD^{(0)}_{n+1}\tau_{n+1}^{-1}\tau_{n+1}^{j} & i-j=-1\\
(-1)^{|a|}\tau_{n+1}^{-j}\lD^{(0)}_{n+1}\tau_{n+1}^{j} & i-j=0\\
(-1)^{|a|+1}\tau_{n+1}^{-j}a^{(0)}_{n}\lm_{n+1}^{(i-j-1)}\tau_{n+1}^{j} & 1\leq i-j\leq n\\
(-1)^{|a|+1}\tau_{n+1}^{-j}\rD^{(0)}_{n+1}\tau_{n+1}^{-1}\tau_{n+1}^{j} & i-j=n+1
\end{cases}\\
=\begin{cases}
(-1)^{|a|+1}a^{(j-1)}_{n}\lm_{n+1}^{(i)} & -n-1\leq i-j\leq -2\\
(-1)^{|a|+1}\rD^{(j-1)}_{n+1} & i-j=-1\\
(-1)^{|a|}\lD^{(j)}_{n+1} & i-j=0\\
(-1)^{|a|+1}a^{(j)}_{n}\lm_{n+1}^{(i-1)} & 1\leq i-j\leq n\\
(-1)^{|a|+1}\rD^{(j)}_{n+1}\tau_{n+1}^{-1} & i-j=n+1
\end{cases}
\end{eqnarray*}
Then
\begin{eqnarray*}
\bem_{n+2}\ID_{n+1}=\sum_{i=0}^{n+1}\sum_{j=1}^{n+1}\lm^{(i)}_{n+2}a^{(j)}_{n+1}\\
=\sum_{j=1}^{n+1}\lm^{(j)}_{n+2}a^{(j)}_{n+1}+\sum_{j=1}^{n+1}\lm^{(j-1)}_{n+2}a^{(j)}_{n+1}
+\sum_{i=0}^{n+1}\sum_{j=i+2}^{n+1}\lm^{(i)}_{n+2}a^{(j)}_{n+1}
+\sum_{i=0}^{n+1}\sum_{j=1}^{i+1}\lm^{(i)}_{n+2}a^{(j)}_{n+1}\\
=(-1)^{|a|}\sum_{j=1}^{n+1}\lD^{(j)}_{n+1}
+(-1)^{|a|+1}\sum_{j=1}^{n+1}\rD^{(j-1)}_{n+1}\\
+(-1)^{|a|+1}\sum_{i=0}^{n+1}\sum_{j=i+2}^{n+1}a^{(j-1)}_{n}\lm_{n+1}^{(i)}
+(-1)^{|a|+1}\sum_{i=0}^{n+1}\sum_{j=1}^{i-1}a^{(j)}_{n}\lm_{n+1}^{(i-1)}\\
=(-1)^{|a|}\dD_{n+1}
+(-1)^{|a|+1}\sum_{i=0}^{n+1}\sum_{j=i+2}^{n+1}a^{(j-1)}_{n}\lm_{n+1}^{(i)}
+(-1)^{|a|+1}\sum_{i=0}^{n+1}\sum_{j=1}^{i-1}a^{(j)}_{n}\lm_{n+1}^{(i-1)}\\
=(-1)^{|a|}\dD_{n+1}
+(-1)^{|a|+1}\sum_{i=0}^{n-1}\sum_{j=i+1}^{n}a^{(j)}_{n}\lm_{n+1}^{(i)}
+(-1)^{|a|+1}\sum_{i=1}^{n}\sum_{j=1}^{i}a^{(j)}_{n}\lm_{n+1}^{(i)}\\
=(-1)^{|a|}\dD_{n+1}+(-1)^{|a|+1}\ID_{n}\bem_{n+1}
\end{eqnarray*}

\medskip

c) First, note that $B^e$ is the restriction of $h^eN$ onto $\rC^e(A)$. Then observe that
\begin{eqnarray*}
N_{n+1}\ID_{n}=N_{n+1}\sum_{i=1}^{n}\tau_{n+1}^{-i}a^{(0)}_{n}\tau_{n}^{i}=\sum_{i=1}^{n}N_{n+1}a^{(0)}_{n}\tau_{n}^{i}
=N_{n+1}a^{(0)}_{n}\sum_{i=1}^{n}\tau_{n}^{i}=N_{n+1}a^{(0)}_{n}N_{n}
\end{eqnarray*}
and similarly
\begin{eqnarray*}
(a^{(0)}_{n}+\ID_{n})N_n=\sum_{i=0}^{n}\tau_{n+1}^{-i}a^{(0)}_{n}\tau_{n}^{i}N_n=\sum_{i=0}^{n}\tau_{n+1}^{-i}a^{(0)}_{n}N_n
=N_{n+1}a^{(0)}_{n}N_{n}
\end{eqnarray*}
Furthermore, one checks easily that
\begin{eqnarray}\label{hba}
h^e_{n+1}a^{(i)}_{n}=(-1)^{|a|+1}a^{(i+1)}_{n+1}h^e_{n}
\end{eqnarray}
Thus,
\begin{eqnarray*}
h^e_{n+1}N_{n+1}\ID_{n}=h^e_{n+1}(a^{(0)}_{n}+\ID_{n})N_n=h^e_{n+1}\sum_{i=0}^{n}a^{(i)}_{n}N_n
=(-1)^{|a|+1}\sum_{i=1}^{n+1}a^{(i)}_{n+1}h^e_{n}N_n\\=(-1)^{|a|+1}\ID_{n+1}h^e_{n}N_n
\end{eqnarray*}

\medskip
d) By (\ref{ID})
\begin{eqnarray*}
{a'}^{(j)}_{n+2}a^{(i)}_{n+1}(a_0[a_1|\ldots |a_n])\\
=\begin{cases}
(-1)^{|sa|\sum_{k=0}^{i-1}|sa_k|}(-1)^{|sa'|\sum_{l=0}^{j-1}|sa_l|}a_0[a_1|\ldots|a_{j-1}|a'|\ldots|a|a_{i}| \ldots|a_n] & j\leq i\\
(-1)^{|sa|\sum_{k=0}^{i-1}|sa_k|}(-1)^{|sa'|\sum_{l=0}^{j-2}|sa_l|}(-1)^{|sa'||sa|}a_0[a_1|\ldots|a_{i-1}|a|\ldots|a'|a_{j-1}|\ldots|a_n] & j>i
\end{cases}
\end{eqnarray*}
which means
\begin{eqnarray}\label{baba}
{a'}^{(j)}_{n+2}a^{(i)}_{n+1}=(-1)^{|sa'||sa|}
\begin{cases}
a^{(i+1)}_{n+2}{a'}^{(j)}_{n+1}, \quad j\leq i\\
a^{(i)}_{n+2}{a'}^{(j-1)}_{n+1}, \quad j>i
\end{cases}
\end{eqnarray}
The latter implies the statement.

\medskip
\subsection{Proof of Proposition \ref{newmix}}\label{proofnewmix}

\noindent{\it Part a)} follows from Lemma \ref{intertw} immediately since in our case  $[D(w),w]=0$, $\mathrm{ad}(w)=0$, and $\bw^2=0$ (the latter is a special case of Part d) of the Lemma).

\medskip

\noindent{\it Part b):} 
One needs to show that
\begin{eqnarray*}
b^e\mathrm{exp}(\bdw)=\mathrm{exp}(\bdw)(\bem +\bw), \quad B^e\mathrm{exp}(\bdw)=\mathrm{exp}(\bdw) B^e
\end{eqnarray*}
The proof is straightforward application of the following formulas
\begin{eqnarray*}
\quad [\bed, \bdw]=2\bw, \quad [\bem, \bdw]=-\bed, \quad [B^e, \bdw]=0, \quad [\bw, \bdw]=0
\end{eqnarray*}
which, in turn, follow from $[D(w),D(w)]=2w$, $\mathrm{ad}(D(w))=\bed$ and Lemma \ref{intertw}.

\medskip

\noindent{\it Part c):} That $\str$ is a morphism of complexes $({\rrCe}(\C[Y]\otimes\DD), \bem)\to ({\rrCe}(\C[Y]),\bem)$ is a special case of a more general fact (see, for instance, \cite{Shk}) but it is also easy to prove in our special case. It is enough to show that $\str\cdot\lm^{(0)}=\lm^{(0)}\cdot\str$ and $\str\cdot\tau=\tau\cdot\str$:
\begin{eqnarray*}
\str\cdot\lm^{(0)}((\phi_0\otimes E_{i_0i_1})[\phi_1\otimes E_{i_1i_2}|\ldots |\phi_n\otimes E_{i_ni_0}])\\
=(-1)^{|v_{i_0}|+|v_{i_1}|}\str((\phi_0\phi_1\otimes E_{i_0i_2})[\phi_2\otimes E_{i_2i_3}|\ldots |\phi_n\otimes E_{i_ni_0}])\\
=(-1)^{|v_{i_0}|+|v_{i_1}|}(-1)^{(n-2)|v_{i_0}|+\sum_{s=2}^n|v_{i_s}|}\phi_0\phi_1[\phi_2|\ldots |\phi_n]\\
=(-1)^{(n-1)|v_{i_0}|+\sum_{s=1}^n|v_{i_s}|}\phi_0\phi_1[\phi_2|\ldots |\phi_n]\\
=\lm^{(0)}\cdot\str((\phi_0\otimes E_{i_0i_1})[\phi_1\otimes E_{i_1i_2}|\ldots |\phi_n\otimes E_{i_ni_0}]),
\end{eqnarray*}
\begin{eqnarray*}
\str\cdot\tau_{n+1}((\phi_0\otimes E_{i_0i_1})[\phi_1\otimes E_{i_1i_2}|\ldots |\phi_n\otimes E_{i_ni_0}])\\
=(-1)^{(|s\phi_0|+|v_{i_0}|+|v_{i_1}|)(\sum_{t=1}^{n-1}(|s\phi_t|+|v_{i_t}|+|v_{i_{t+1}}|)+|s\phi_n|+|v_{i_n}|+|v_{i_{0}}|)}\str((\phi_1\otimes E_{i_1i_2})[\ldots|\phi_0\otimes E_{i_0i_1}])\\
=(-1)^{(1+|v_{i_0}|+|v_{i_1}|)(n+|v_{i_0}|+|v_{i_1}|)}\str((\phi_1\otimes E_{i_1i_2})[\ldots|\phi_0\otimes E_{i_0i_1}])\\
=(-1)^{n(1+|v_{i_0}|+|v_{i_1}|)}\str((\phi_1\otimes E_{i_1i_2})[\ldots|\phi_0\otimes E_{i_0i_1}])\\
=(-1)^{n(1+|v_{i_0}|+|v_{i_1}|)}(-1)^{(n-1)|v_{i_1}|+|v_{i_0}|+\sum_{s=2}^n|v_{i_s}|}\phi_1[\ldots |\phi_0]\\
=(-1)^{n+(n+1)|v_{i_0}|+\sum_{s=1}^n|v_{i_s}|}\phi_1[\ldots |\phi_0]=(-1)^{n+(n-1)|v_{i_0}|+\sum_{s=1}^n|v_{i_s}|}\phi_1[\ldots |\phi_0]\\
=(-1)^{(n-1)|v_{i_0}|+\sum_{s=1}^n|v_{i_s}|}\tau_{n+1}(\phi_0[\phi_1|\ldots |\phi_n])\\
=\tau_{n+1}\cdot\str((\phi_0\otimes E_{i_0i_1})[\phi_1\otimes E_{i_1i_2}|\ldots |\phi_n\otimes E_{i_ni_0}])
\end{eqnarray*}
It remains to verify that $\str$ commutes with $\bw$ and $B^e$:
\begin{eqnarray*}
\str\cdot w^{(t)}_{n+1}((\phi_0\otimes E_{i_0i_1})[\phi_1\otimes E_{i_1i_2}|\ldots |\phi_n\otimes E_{i_ni_0}])\\
=(-1)^{t+|v_{i_0}|+|v_{i_t}|}\str((\phi_0\otimes E_{i_0i_1})[\phi_1\otimes E_{i_1i_2}|\ldots|\phi_{t-1}\otimes E_{i_{t-1}i_t}|w\otimes1|\ldots |\phi_n\otimes E_{i_ni_0}])\\
=(-1)^{t+|v_{i_0}|+|v_{i_t}|}\str((\phi_0\otimes E_{i_0i_1})[\phi_1\otimes E_{i_1i_2}|\ldots|\phi_{t-1}\otimes E_{i_{t-1}i_t}|w\otimes E_{i_{t}i_t}|\ldots |\phi_n\otimes E_{i_ni_0}])\\
=(-1)^{t+|v_{i_0}|+|v_{i_t}|}(-1)^{n|v_{i_0}|+|v_{i_t}|+\sum_{s=1}^n|v_{i_s}|}\phi_0[\phi_1|\ldots|\phi_{t-1}|w|\ldots |\phi_n])\\
=(-1)^{t}(-1)^{(n-1)|v_{i_0}|+\sum_{s=1}^n|v_{i_s}|}\phi_0[\phi_1|\ldots|\phi_{t-1}|w|\ldots |\phi_n])\\
= w^{(t)}_{n+1}\cdot\str((\phi_0\otimes E_{i_0i_1})[\phi_1\otimes E_{i_1i_2}|\ldots |\phi_n\otimes E_{i_ni_0}])
\end{eqnarray*}
(in the above computation, we implicitly use the fact that $|E_{ij}|=|v_{i}|+|v_{j}|$),
\begin{eqnarray*}
\str\cdot h^e((\phi_0\otimes E_{i_0i_1})[\phi_1\otimes E_{i_1i_2}|\ldots |\phi_n\otimes E_{i_ni_0}])\\
=\str(e[\phi_0\otimes E_{i_0i_1}[\phi_1\otimes E_{i_1i_2}|\ldots |\phi_n\otimes E_{i_ni_0}])\\
=(-1)^{(n+1)|v_{i_0}|+\sum_{s=1}^n|v_{i_s}|}e[\phi_0|\ldots |\phi_n]=h^e\cdot\str((\phi_0\otimes E_{i_0i_1})[\phi_1\otimes E_{i_1i_2}|\ldots |\phi_n\otimes E_{i_ni_0}])
\end{eqnarray*}

\medskip

\noindent{\it Part d):} The statement is an easy variation on the classical Hochschild-Kostant-Rosenberg map (cf. \cite{CalTu}).

\medskip

\noindent{\it Part e):} The quasi-isomorphism of complexes constructed in \cite{Seg} is the composition of $\epsilon\cdot\str\cdot\mathrm{exp}(-\bdw)$ with the embedding
$$
({\rC}(A_w), b)\to({\rC^e}(A_w), b^e)
$$
Since the latter embedding is a quasi-isomorphism (see Proposition \ref{quasiisocycl}), we conclude that $\epsilon\cdot\str\cdot\mathrm{exp}(-\bdw): ({\rC^e}(A_w), b^e)\to (\Omega(Y), -dw)$ is a quasi-isomorphism. It remains to apply Lemma \ref{standard}.
\hfill $\Box$

\medskip

\subsection{Proof of Proposition \ref{conntw}}\label{proofconntw}

It will be convenient for us to work with $\widetilde{\nabla}$ instead of $\nabla$ (see Corollary \ref{alterconncycl}). The proof of Proposition \ref{conntw} will be a combination of two Lemmas.

Let us extend $\mathrm{exp}(-\bdw)$ to a map from  ${\rrCe}(A_w)$ to itself.
\medskip
\begin{lemma}\label{0}
$$
\mathrm{exp}(-\bdw)(\Udl+\Uml)\mathrm{exp}(\bdw)=\Uml+\Uwc,
$$
$$
\mathrm{exp}(-\bdw)(\Vdl+\Vml)\mathrm{exp}(\bdw)=\Vml+\Vwc,
$$
In particular\footnote{Note the similarity between (\ref{alterext}) and (\ref{00}); the connections seem to be special cases of a more general formula which should work for arbitrary curved $A_\infty$ algebras.}
\begin{equation}\label{00}
\frac{d}{du}+\frac{\Uml+\Uwc}{u^2}+\frac{\Vml+\Vwc}{u}
\end{equation}
is a $u$-connection on $({\rrCe}(A_w),\bem+\bw, B^e)$ homotopy gauge equivalent to $\widetilde{\nabla}$ and $\nabla$.
\end{lemma}
\noindent{\bf Proof.} We will be using the following formulas obtained earlier (see the proof of Lemma \ref{intertw}):
\begin{eqnarray}\label{ldbdw}
\ld^{(i)}_{n+2}D(w)^{(j)}_{n+1}
=\begin{cases}
2 w^{(j)}_{n+1},\quad i=j\\
D(w)^{(j)}_{n+1}\ld^{(i-1)}_{n+1},\quad i>j\\
D(w)^{(j)}_{n+1}\ld^{(i)}_{n+1},\quad i<j\\
\end{cases}
\end{eqnarray}
\begin{eqnarray}\label{lmbdw}
\lm^{(i)}_{n+2}D(w)^{(j)}_{n+1}
=\begin{cases}
D(w)^{(j-1)}_{n}\lm_{n+1}^{(i)} & -n-1\leq i-j\leq -2\\
\rdw^{(j-1)}_{n+1} & i-j=-1\\
-\ldw^{(j)}_{n+1} & i-j=0\\
D(w)^{(j)}_{n}\lm_{n+1}^{(i-1)} & 1\leq i-j\leq n\\
\rdw^{(j)}_{n+1}\tau_{n+1}^{-1} & i-j=n+1
\end{cases}
\end{eqnarray}
\begin{eqnarray}\label{lmbw}
\lm^{(i)}_{n+2} w^{(j)}_{n+1}
=\begin{cases}
- w^{(j-1)}_{n}\lm_{n+1}^{(i)} & -n-1\leq i-j\leq -2\\
-\rw^{(j-1)}_{n+1} & i-j=-1\\
\lw^{(j)}_{n+1} & i-j=0\\
- w^{(j)}_{n}\lm_{n+1}^{(i-1)} & 1\leq i-j\leq n\\
-\rw^{(j)}_{n+1}\tau_{n+1}^{-1} & i-j=n+1
\end{cases}
\end{eqnarray}

In what follows, we represent $\bdw$ by the matrix
$$
\begin{pmatrix}\bdwne & 0 \\ 0 & \bdwne \end{pmatrix}
$$
as we agreed in section \ref{auxlem}.

\medskip

Let us simplify $\mathrm{exp}(-\bdw)\Udl\mathrm{exp}(\bdw)$. Clearly,
\begin{eqnarray*}
2[\Udl,\bdw]_n\\=\begin{pmatrix}-\lm_{n+1}^{(0)}\ld_{n+1}^{(1)}\bdwne_n+\bdwne_{n-1}\lm_{n}^{(0)}\ld_{n}^{(1)} & -\lm_{n+1}^{(0)}\ld_{n+1}^{(1)}\bdwne_n+\bdwne_{n-1}\lm_{n}^{(0)}\ld_{n}^{(1)} \\0 & 0 \end{pmatrix}
\end{eqnarray*}
By (\ref{ldbdw}), (\ref{lmbdw})
\begin{eqnarray*}
-\lm_{n+1}^{(0)}\ld_{n+1}^{(1)}\bdwne_n=-\sum_{j=1}^n\lm_{n+1}^{(0)}\ld_{n+1}^{(1)}D(w)^{(j)}_n=-\lm_{n+1}^{(0)}\ld_{n+1}^{(1)}D(w)^{(1)}_n-\sum_{j=2}^n\lm_{n+1}^{(0)}\ld_{n+1}^{(1)}D(w)^{(j)}_n\\
=-2\lm_{n+1}^{(0)} w^{(1)}_{n}-\sum_{j=2}^n\lm_{n+1}^{(0)}D(w)^{(j)}_n\ld_{n}^{(1)}=-2\lm_{n+1}^{(0)} w^{(1)}_{n}-\bdwne_{n-1}\lm_{n}^{(0)}\ld_{n}^{(1)}
\end{eqnarray*}
and therefore
\begin{eqnarray}\label{udlbdw}
[\Udl,\bdw]=2\Uwc
\end{eqnarray}
The operators $\bdw$ and $-\lm_{n+1}^{(0)} w_{n}^{(1)}=\rw^{(0)}_{n}$ are easily seen to commute, so we conclude that
\begin{eqnarray}\label{expudl}
\mathrm{exp}(-\bdw)\Udl\mathrm{exp}(\bdw)=\Udl+2\Uwc
\end{eqnarray}

Furthermore, using $-\lm_{n+1}^{(0)}\lm_{n+2}^{(1)}=\lm_{n+1}^{(0)}\lm_{n+2}^{(0)}$
\begin{eqnarray*}
2[\Uml,\bdw]_{n+1}\\=\begin{pmatrix} \lm_{n+1}^{(0)}\lm_{n+2}^{(0)}\bdwne_{n+1}-\bdwne_{n-1}\lm_{n}^{(0)}\lm_{n+1}^{(0)} & \lm_{n+1}^{(0)}\lm_{n+2}^{(0)}\bdwne_{n+1}-\bdwne_{n-1}\lm_{n}^{(0)}\lm_{n+1}^{(0)} \\0 & 0 \end{pmatrix}
\end{eqnarray*}
By (\ref{lmbdw})
\begin{eqnarray*}
\lm_{n+1}^{(0)}\lm_{n+2}^{(0)}\bdwne_{n+1}=\sum_{j=1}^{n+1}\lm_{n+1}^{(0)}\lm_{n+2}^{(0)}D(w)^{(j)}_{n+1}
=\lm_{n+1}^{(0)}\lm_{n+2}^{(0)}D(w)^{(1)}_{n+1}+\sum_{j=2}^{n+1}\lm_{n+1}^{(0)}\lm_{n+2}^{(0)}D(w)^{(j)}_{n+1}\\
=\lm_{n+1}^{(0)}\rdw^{(0)}_{n+1}+\sum_{j=2}^{n+1}\lm_{n+1}^{(0)}D(w)^{(j-1)}_{n}\lm_{n+1}^{(0)}
=\lm_{n+1}^{(0)}\rdw^{(0)}_{n+1}+\sum_{j=1}^{n}\lm_{n+1}^{(0)}D(w)^{(j)}_{n}\lm_{n+1}^{(0)}\\
=\lm_{n+1}^{(0)}\rdw^{(0)}_{n+1}+\lm_{n+1}^{(0)}D(w)^{(1)}_{n}\lm_{n+1}^{(0)}+\sum_{j=2}^{n}\lm_{n+1}^{(0)}D(w)^{(j)}_{n}\lm_{n+1}^{(0)}\\
=\lm_{n+1}^{(0)}\rdw^{(0)}_{n+1}+\rdw^{(0)}_{n}\lm_{n+1}^{(0)}+\bdwne_{n-1}\lm_{n}^{(0)}\lm_{n+1}^{(0)}
\end{eqnarray*}
Observe that
\begin{eqnarray*}
\lm_{n+1}^{(0)}\rdw^{(0)}_{n+1}(a_0[a_1|\ldots |a_n])+\rdw^{(0)}_{n}\lm_{n+1}^{(0)}(a_0[a_1|\ldots |a_n])=-a_0[D(w),a_1][a_2|\ldots |a_n]\\
=\lm_{n+1}^{(0)}\ld_{n+1}^{(1)}(a_0[a_1|\ldots |a_n])
\end{eqnarray*}
Hence
\begin{eqnarray*}
[\Uml,\bdw]=-\Udl
\end{eqnarray*}
Therefore, by (\ref{udlbdw})
\begin{eqnarray*}
[[\Uml,\bdw],\bdw]=-[\Udl,\bdw]=-2\Uwc
\end{eqnarray*}
and we conclude that
\begin{eqnarray*}
\mathrm{exp}(-\bdw)\Uml\mathrm{exp}(\bdw)=\Umc-\Udl-\Uwc
\end{eqnarray*}
By (\ref{expudl}) and the latter formula
\begin{eqnarray*}
\mathrm{exp}(-\bdw)(\Udl+\Uml)\mathrm{exp}(\bdw)=\Umc+\Uwc
\end{eqnarray*}

\medskip
Let us compute $\mathrm{exp}(-\bdw)\Vdl\mathrm{exp}(\bdw)$ now. We have
\begin{eqnarray*}
2[\Vdl,\bdw]_n\\=\begin{pmatrix}0 & 0 \\ -h^e_{n+1}\sum_{i=1}^{n}\sum_{l=i+1}^{n+1}\tau_{n+1}^{l}\ld_{n+1}^{(i)}\bdwne_n+\bdwne_{n+1}h^e_{n}\sum_{i=1}^{n-1}\sum_{l=i+1}^{n}\tau_{n}^{l}\ld_{n}^{(i)} & 0 \end{pmatrix}
\end{eqnarray*}
By (\ref{ldbdw})
\begin{eqnarray*}
-h^e_{n+1}\sum_{i=1}^{n}\sum_{l=i+1}^{n+1}\tau_{n+1}^{l}\ld_{n+1}^{(i)}\bdwne_n=-h^e_{n+1}\sum_{i=1}^{n}\sum_{l=i+1}^{n+1}\sum_{j=1}^{n}\tau_{n+1}^{l}\ld_{n+1}^{(i)}D(w)_n^{(j)}\\
=-h^e_{n+1}\sum_{i=1}^{n-1}\sum_{l=i+1}^{n+1}\sum_{j=i+1}^{n}\tau_{n+1}^{l}D(w)^{(j)}_{n}\ld^{(i)}_{n}-h^e_{n+1}\sum_{i=2}^{n}\sum_{l=i+1}^{n+1}\sum_{j=1}^{i-1}\tau_{n+1}^{l}D(w)^{(j)}_{n}\ld^{(i-1)}_{n}\\-2h^e_{n+1}\sum_{i=1}^{n}\sum_{l=i+1}^{n+1}\tau_{n+1}^{l} w^{(i)}_{n}
\end{eqnarray*}
Using the formula $D(w)^{(j)}_{n}=\tau_{n+1}^{-j}D(w)^{(0)}_{n}\tau_{n}^{j}$, the latter equals
\begin{eqnarray*}
-h^e_{n+1}\sum_{i=1}^{n-1}\sum_{l=i+1}^{n+1}\sum_{j=i+1}^{n}\tau_{n+1}^{l-j}D(w)^{(0)}_{n}\tau_{n}^{j}\ld^{(i)}_{n}-h^e_{n+1}\sum_{i=2}^{n}\sum_{l=i+1}^{n+1}\sum_{j=1}^{i-1}\tau_{n+1}^{l-j}D(w)^{(0)}_{n}\tau_{n}^{j}\ld^{(i-1)}_{n}\\-2h^e_{n+1}\sum_{i=1}^{n}\sum_{l=i+1}^{n+1}\tau_{n+1}^{l} w^{(i)}_{n}\\
=-h^e_{n+1}\sum_{i=1}^{n-1}\sum_{l=i+1}^{n+1}\sum_{j=i+1}^{l-1}D(w)^{(n+1+j-l)}_{n}\tau_{n}^{l-1}\ld^{(i)}_{n}-h^e_{n+1}\sum_{i=1}^{n-1}\sum_{l=i+1}^{n+1}\sum_{j=l}^{n}D(w)^{(j-l)}_{n}\tau_{n}^{l}\ld^{(i)}_{n}\\
-h^e_{n+1}\sum_{i=2}^{n}\sum_{l=i+1}^{n+1}\sum_{j=1}^{i-1}D(w)^{(n+1+j-l)}_{n}\tau_{n}^{l-1}\ld^{(i-1)}_{n}-2h^e_{n+1}\sum_{i=1}^{n}\sum_{l=i+1}^{n+1}\tau_{n+1}^{l} w^{(i)}_{n}
\end{eqnarray*}
After reindexing we get
\begin{eqnarray*}
-h^e_{n+1}\sum_{i=1}^{n-1}\sum_{l=i+1}^{n}\sum_{j=n+i+1-l}^{n}D(w)^{(j)}_{n}\tau_{n}^{l}\ld^{(i)}_{n}-h^e_{n+1}\sum_{i=1}^{n-1}\sum_{l=i+1}^{n}\sum_{j=0}^{n-l}D(w)^{(j)}_{n}\tau_{n}^{l}\ld^{(i)}_{n}
\end{eqnarray*}
\begin{eqnarray*}
-h^e_{n+1}\sum_{i=1}^{n-1}\sum_{l=i+1}^{n}\sum_{j=n+1-l}^{n+i-l}D(w)^{(j)}_{n}\tau_{n}^{l}\ld^{(i)}_{n}-2h^e_{n+1}\sum_{i=1}^{n}\sum_{l=i+1}^{n+1}\tau_{n+1}^{l} w^{(i)}_{n}\\
=-h^e_{n+1}\sum_{i=1}^{n-1}\sum_{l=i+1}^{n}\sum_{j=0}^{n}D(w)^{(j)}_{n}\tau_{n}^{l}\ld^{(i)}_{n}-2h^e_{n+1}\sum_{i=1}^{n}\sum_{l=i+1}^{n+1}\tau_{n+1}^{l} w^{(i)}_{n}
\end{eqnarray*}
By (\ref{hba}) the latter equals
\begin{eqnarray*}
-\sum_{j=0}^{n}D(w)^{(j+1)}_{n+1}h^e_{n}\sum_{i=1}^{n-1}\sum_{l=i+1}^{n}\tau_{n}^{l}\ld^{(i)}_{n}-2h^e_{n+1}\sum_{i=1}^{n}\sum_{l=i+1}^{n+1}\tau_{n+1}^{l} w^{(i)}_{n}\\
=-\bdwne_{n+1}h^e_{n}\sum_{i=1}^{n-1}\sum_{l=i+1}^{n}\tau_{n}^{l}\ld^{(i)}_{n}-2h^e_{n+1}\sum_{i=1}^{n}\sum_{l=i+1}^{n+1}\tau_{n+1}^{l} w^{(i)}_{n}
\end{eqnarray*}
Thus,
\begin{eqnarray}\label{vdlbdw}
[\Vdl,\bdw]=2\Vwc
\end{eqnarray}

By essentially repeating the computation we just did word for word, but using (\ref{baba}) instead of (\ref{ldbdw}), one can show that
\begin{eqnarray*}
[[\Vdl,\bdw],\bdw]=0
\end{eqnarray*}
Therefore,
\begin{eqnarray}\label{expvdl}
\mathrm{exp}(-\bdw)\Vdl\mathrm{exp}(\bdw)=\Vdl+2\Vwc
\end{eqnarray}

\medskip

It remains to compute $\mathrm{exp}(-\bdw)\Vml\mathrm{exp}(\bdw)$. The computation is very similar to the previous one.
\begin{eqnarray*}
2[\Vml,\bdw]_n\\=\begin{pmatrix}0 & 0 \\ -h^e_n\sum_{i=1}^{n-1}\sum_{l=i+1}^{n}\tau_{n}^{l}\lm_{n+1}^{(i)}\bdwne_n+\bdwne_{n+1}h^e_n\sum_{i=1}^{n-1}\sum_{l=i+1}^{n}\tau_{n}^{l}\lm_{n+1}^{(i)} & 0 \end{pmatrix}
\end{eqnarray*}
By (\ref{lmbdw})
\begin{eqnarray*}
-h^e_n\sum_{i=1}^{n-1}\sum_{l=i+1}^{n}\tau_{n}^{l}\lm_{n+1}^{(i)}\bdwne_n=-h^e_n\sum_{i=1}^{n-1}\sum_{l=i+1}^{n}\sum_{j=1}^{n}\tau_{n}^{l}\lm_{n+1}^{(i)}D(w)_n^{(j)}\\
=-h^e_n\sum_{i=1}^{n-2}\sum_{l=i+1}^{n}\sum_{j=i+2}^{n}\tau_{n}^{l}D(w)^{(j-1)}_{n-1}\lm_{n}^{(i)}
-h^e_n\sum_{i=2}^{n-1}\sum_{l=i+1}^{n}\sum_{j=1}^{i-1}\tau_{n}^{l}D(w)^{(j)}_{n-1}\lm_{n}^{(i-1)}\\
+h^e_n\sum_{i=1}^{n-1}\sum_{l=i+1}^{n}\tau_{n}^{l}\ldw^{(i)}_{n}
-h^e_n\sum_{i=1}^{n-1}\sum_{l=i+1}^{n}\tau_{n}^{l}\rdw^{(i)}_{n}
\end{eqnarray*}
\begin{eqnarray*}
=-h^e_n\sum_{i=1}^{n-2}\sum_{l=i+1}^{n}\sum_{j=i+1}^{n-1}\tau_{n}^{l}D(w)^{(j)}_{n-1}\lm_{n}^{(i)}
-h^e_n\sum_{i=1}^{n-2}\sum_{l=i+2}^{n}\sum_{j=1}^{i}\tau_{n}^{l}D(w)^{(j)}_{n-1}\lm_{n}^{(i)}\\
+h^e_n\sum_{i=1}^{n-1}\sum_{l=i+1}^{n}\tau_{n}^{l}\ld^{(i)}_{n}
\end{eqnarray*}
Using $D(w)^{(j)}_{n-1}=\tau_{n}^{-j}D(w)^{(0)}_{n-1}\tau_{n-1}^{j}$, we get
\begin{eqnarray*}
-h^e_n\sum_{i=1}^{n-2}\sum_{l=i+1}^{n}\sum_{j=i+1}^{n-1}\tau_{n}^{l-j}D(w)^{(0)}_{n-1}\tau_{n-1}^{j}\lm_{n}^{(i)}
-h^e_n\sum_{i=1}^{n-2}\sum_{l=i+2}^{n}\sum_{j=1}^{i}\tau_{n}^{l-j}D(w)^{(0)}_{n-1}\tau_{n-1}^{j}\lm_{n}^{(i)}\\
+h^e_n\sum_{i=1}^{n-1}\sum_{l=i+1}^{n}\tau_{n}^{l}\ld^{(i)}_{n}\\
=-h^e_n\sum_{i=1}^{n-2}\sum_{l=i+1}^{n}\sum_{j=i+1}^{l-1}D(w)^{(n+j-l)}_{n-1}\tau_{n-1}^{l-1}\lm_{n}^{(i)}
-h^e_n\sum_{i=1}^{n-2}\sum_{l=i+1}^{n}\sum_{j=l}^{n-1}D(w)^{(j-l)}_{n-1}\tau_{n-1}^{l}\lm_{n}^{(i)}\\
-h^e_n\sum_{i=1}^{n-2}\sum_{l=i+2}^{n}\sum_{j=1}^{i}D(w)^{(n+j-l)}_{n-1}\tau_{n-1}^{l-1}\lm_{n}^{(i)}
+h^e_n\sum_{i=1}^{n-1}\sum_{l=i+1}^{n}\tau_{n}^{l}\ld^{(i)}_{n}
\end{eqnarray*}
Just as in the previous computation, re-indexing shows that the latter expression equals
\begin{eqnarray*}
-h^e_n\sum_{j=0}^{n-1}D(w)^{(j)}_{n-1}\sum_{i=1}^{n-2}\sum_{l=i+1}^{n-1}\tau_{n-1}^{l}\lm_{n}^{(i)}
+h^e_n\sum_{i=1}^{n-1}\sum_{l=i+1}^{n}\tau_{n}^{l}\ld^{(i)}_{n}
\end{eqnarray*}
and, by (\ref{hba}), we finally get
\begin{eqnarray*}
-\bdwne_{n}h^e_{n-1}\sum_{i=1}^{n-2}\sum_{l=i+1}^{n-1}\tau_{n-1}^{l}\lm_{n}^{(i)}
+h^e_n\sum_{i=1}^{n-1}\sum_{l=i+1}^{n}\tau_{n}^{l}\ld^{(i)}_{n}
\end{eqnarray*}
Thus,
\begin{eqnarray*}
[\Vml,\bdw]=-\Vdl
\end{eqnarray*}
and by (\ref{vdlbdw})
\begin{eqnarray*}
[[\Vml,\bdw],\bdw]=-[\Vdl,\bdw]=-2\Vwc
\end{eqnarray*}
Therefore,
\begin{eqnarray*}
\mathrm{exp}(-\bdw)\Vml\mathrm{exp}(\bdw)=\Vmc-\Vdl-\Vwc
\end{eqnarray*}
This, along with (\ref{expvdl}), proves that
\begin{eqnarray*}
\mathrm{exp}(-\bdw)(\Vdl+\Vml)\mathrm{exp}(\bdw)=\Vmc+\Vwc
\end{eqnarray*}

Lemma \ref{0} is proved.
\hfill $\Box$

To finish the proof of Proposition \ref{conntw}, one needs to add the following $u$-morphism to the $u$-connection (\ref{00}):

\begin{lemma} The $u$-morphism
$$
\frac1{u^2}(-\Uml+\Uwc)+\frac1{u}(-\Vml+\Vwc+\Gamma)
$$
from $({\rrCe}(A_w),\bem+\bw, B^e)$ to itself is $u$-homotopic to 0.
\end{lemma}
\noindent{\bf Proof} is similar to that of Corollary \ref{alterconncycl}. We will show that
\begin{equation}\label{III}
2(\Uml-\Uwc)+2u(\Vml-\Vwc-\Gamma)
\end{equation}
is $u$-homotopic to 0.

We will be using the notation introduced in the proof of Corollary \ref{alterconncycl}.

Let us compute $(\bem+\bw+uB^e)H_0+H_0(\bem+\bw+uB^e)$ first. We already know that
\begin{eqnarray*}
\bem H_0+H_0\bem=2\Umc
\end{eqnarray*}
and
\begin{eqnarray*}
B^eH_0+H_0B^e=T
\end{eqnarray*}
By (\ref{lmbw})
\begin{eqnarray*}
\bw H_0+H_0\bw=\begin{pmatrix}\bwne\lm^{(0)}+\lm^{(0)}\bwne & \bwne\lm^{(0)}+\lm^{(0)}\bwne \\0 & 0 \end{pmatrix}
=-2\Uwc
\end{eqnarray*}
Thus,
\begin{eqnarray}\label{VV}
(\bem+\bw+uB^e)H_0+H_0(\bem+\bw+uB^e)=(2\Umc-2\Uwc)+uT
\end{eqnarray}

Let us compute $(\bem+\bw+uB^e)H_1+H_1(\bem+\bw+uB^e)$. We already know that
$$
B^eH_1+H_1B^e=0
$$
and
\begin{eqnarray}\label{VVV}
\bem H_1+H_1\bem=-T-2\Gamma+2\Vmc
\end{eqnarray}
It remains to simplify $\bw H_1+H_1\bw$. The computation is very similar to that of $\bem H_1+H_1\bem$, so let us simply present the result:
\begin{eqnarray*}
\bw H_1+H_1\bw=-2\Vwc
\end{eqnarray*}
This, together with (\ref{VV}) and (\ref{VVV}) finishes the proof.
\hfill $\Box$

\medskip

\subsection{Proof of Proposition \ref{connw}}\label{proofconnw}

We need to show that $$\str\cdot\Uwc=\Uwc\cdot\str,\quad\str\cdot\Vwc=\Vwc\cdot\str$$
This is a consequence of the proof of Part c) in section \ref{proofnewmix}.

\medskip

\subsection{Proof of Proposition \ref{epsconn}}\label{proofepsconn}

We need to show that the morphism of complexes
\begin{eqnarray*}
\left(\frac{w}{u^2}+\frac{\Gamma}{u}\right)\epsilon-\epsilon\left(\frac{2\Uwc}{u^2}+\frac{2\Vwc+\Gamma}{u}\right)
\end{eqnarray*}
from $({\rrCe}(\C[Y])((u)), \bem+\bw+uB^e)$ to $(\Omega(Y)((u))\,,\, -dw+ud)$ is homotopic to 0. Both ${w}\epsilon-2\epsilon\Uwc$ and $\Gamma\epsilon-\epsilon\Gamma$ are easily seen to vanish, so we only need to show that $\epsilon\Vwc$ is homotopic to 0.

Observe that
\begin{eqnarray*}
\epsilon_{n+2} h^e_{n+1}\tau_{n+1}^{j}(\phi_0[\phi_1|\ldots |\phi_n])=(-1)^{jn}\epsilon_{n+2} h^e_{n+1}(\phi_j[\phi_{j+1}|\ldots |\phi_{j-1}])\\
=(-1)^{jn}\frac{1}{(n+1)!}d\phi_j\wedge d\phi_{j+1}\wedge\ldots \wedge d\phi_{j-1}=\frac{1}{(n+1)!}d\phi_0\wedge d\phi_{1}\wedge\ldots \wedge d\phi_{n}\\
=\frac1{n+1}d(\epsilon_{n+1}(\phi_0[\phi_1|\ldots |\phi_n]))
\end{eqnarray*}
and
\begin{eqnarray*}
d(\epsilon_{n+1}( w^{(i)}_{n}(\phi_0[\phi_1|\ldots |\phi_{n-1}])))=(-1)^id(\epsilon_{n+1}(\phi_0[\phi_1|\ldots|\phi_{i-1}|w|\ldots |\phi_{n-1}]))\\
=(-1)^i\frac1{n!}d(\phi_0d\phi_1\wedge\ldots\wedge d\phi_{i-1}\wedge dw\wedge\ldots \wedge d\phi_{n-1})=\frac1{n!}dw\wedge d\phi_0\wedge d\phi_1\wedge\ldots\wedge  d\phi_{n-1}\\
=\frac1{n}dw\wedge d(\epsilon_{n}(\phi_0[\phi_1|\ldots |\phi_{n-1}]))
\end{eqnarray*}
Therefore,
\begin{eqnarray*}
2\epsilon_{n+2}\Vwc_n=-\sum_{i=1}^{n}\sum_{j=i+1}^{n+1}\epsilon_{n+2}\cdot h^e_{n+1}\cdot \tau_{n+1}^{j}\cdot  w^{(i)}_{n}
=-\sum_{i=1}^{n}\frac{n+1-i}{n+1}d\cdot \epsilon_{n+1}\cdot  w^{(i)}_{n}\\
=-\sum_{i=1}^{n}\frac{n+1-i}{n(n+1)}(dw\cdot d)\cdot \epsilon_{n}=\frac12 (dw\cdot d)\cdot \epsilon_{n}
\end{eqnarray*}
It remains to show that
$$
(dw\cdot d):(\Omega(Y)((u))\,,\, -dw+ud)\to (\Omega(Y)((u))\,,\, -dw+ud)
$$
is homotopic to 0, which is obvious:
$$
dw\cdot d=(-dw+ud)H+ H(-dw+ud), \quad H:= \frac1{u}w\cdot d
$$
\hfill $\Box$


\bigskip

\begin{thebibliography}{99}

\bibitem{Bar1} S.~Barannikov, Quantum periods. I. Semi-infinite variations of Hodge structures. Internat. Math. Res. Notices 2001, no. 23, 1243–-1264.
\bibitem{Bar2} S.~Barannikov, Non-commutative periods and mirror symmetry in higher dimensions. Comm. Math. Phys. 228 (2002), no. 2, 281–-325.
\bibitem{Br} E.~Brieskorn, Die Monodromie der isolierten Singularitäten von Hyperflächen. Manuscripta Math. 2 1970 103–-161.
\bibitem{CalTu} A.~Caldararu, J.~Tu, Curved A-infinity algebras and Landau-Ginzburg models. arXiv:1007.2679
\bibitem{DTT} V.~A.~Dolgushev, D.~E.~Tamarkin, B.~L.~Tsygan, Noncommutative calculus and the Gauss-Manin connection. Higher structures in geometry and physics, 139–158, Progr. Math., 287, Birkhäuser/Springer, New York, 2011.
\bibitem{Dyck} T.~Dyckerhoff, Compact generators in categories of matrix factorizations. Duke Math. J. 159 (2011), no. 2, 223–-274.
\bibitem{Get} E.~Getzler, Cartan homotopy formulas and the Gauss-Manin connection in cyclic homology. (English summary) Quantum deformations of algebras and their representations (Ramat-Gan, 1991/1992; Rehovot, 1991/1992), 65–78, Israel Math. Conf. Proc., 7, Bar-Ilan Univ., Ramat Gan, 1993.
\bibitem{HS} C.~Hertling, C.~Sabbah, Examples of non-commutative Hodge structures. J. Inst. Math. Jussieu 10 (2011), no. 3, 635–-674.
\bibitem{HSe} C.~Hertling, C.~Sevenheck, Twistor structures, tt*-geometry and singularity theory. From Hodge theory to integrability and TQFT tt*-geometry, 49–-73, Proc. Sympos. Pure Math., 78, Amer. Math. Soc., Providence, RI, 2008.
\bibitem{Iri} H.~Iritani, An integral structure in quantum cohomology and mirror symmetry for toric orbifolds. Adv. Math. 222 (2009), no. 3, 1016–-1079.
\bibitem{Ka} D. Kaledin, Motivic structures in non-commutative geometry. Proceedings of the International Congress of Mathematicians. Volume II, 461–-496, Hindustan Book Agency, New Delhi, 2010.
\bibitem{KKP} L.~Katzarkov, M.~Kontsevich, T.~Pantev, Hodge theoretic aspects of mirror symmetry.
From Hodge theory to integrability and TQFT $tt^*$-geometry, 87–-174, Proc. Sympos. Pure Math., 78, Amer. Math. Soc., Providence, RI, 2008.
\bibitem{Ke} B.~Keller, On the cyclic homology of ringed spaces and schemes. Doc. Math. {\bf 3} (1998), 231–-259.
\bibitem{Kel} B. Keller, On differential graded categories.  International Congress of Mathematicians. Vol. II,  151--190, Eur. Math. Soc., Z\"{u}rich, 2006.
\bibitem{Kon} M.~Kontsevich, XI Solomon Lefschetz Memorial Lecture series: Hodge structures in non-commutative geometry. Notes by Ernesto Lupercio. Contemp. Math., 462, Non-commutative geometry in mathematics and physics, 1–-21, Amer. Math. Soc., Providence, RI, 2008.
\bibitem{KS} M.~Kontsevich, Y.~Soibelman, Notes on $A_\infty$-algebras, $A_\infty$-categories and non-commutative geometry. Homological mirror symmetry, 153–-219, Lecture Notes in Phys., 757, Springer, Berlin, 2009.
\bibitem{Lod} J.-L. Loday, Cyclic homology. Grundlehren der Mathematischen Wissenschaften, 301. Springer-Verlag, Berlin, 1992.
\bibitem{PP} A.~Polishchuk, L.~Positselski, Hochschild (co)homology of the second kind I. arXiv:1010.0982
\bibitem{Ri} G.~S.~Rinehart, Differential forms on general commutative algebras. Trans. Amer. Math. Soc. 108 (1963), 195–-222.
\bibitem{S} C.~Sabbah, D$\acute{\rm e}$formations isomonodromiques et vari$\acute{\rm e}$t$\acute{\rm e}$s de Frobenius, Savoirs Actuels,
CNRS $\acute{\rm E}$ditions $\&$ EDP Sciences, Paris, 2002.
\bibitem{Sa} C.~Sabbah, Fourier-Laplace transform of a variation of polarized complex Hodge structure, II. New developments in algebraic geometry, integrable systems and mirror symmetry (RIMS, Kyoto, 2008), 289-–347, Adv. Stud. Pure Math., 59, Math. Soc. Japan, Tokyo, 2010.
\bibitem{Sab} C.~Sabbah, Non-commutative Hodge structures. arXiv:1107.5890
\bibitem{ST} K.~Saito, A.~Takahashi, From primitive forms to Frobenius manifolds. From Hodge theory to integrability and TQFT tt*-geometry, 31–-48, Proc. Sympos. Pure Math., 78, Amer. Math. Soc., Providence, RI, 2008.
\bibitem{SS} J.~Scherk, J.~H.~M.~Steenbrink, On the mixed Hodge structure on the cohomology of the Milnor fibre, Math. Ann. 271 (1985), no. 4, 641–-665.
\bibitem{Sch} M.~Schulze, A normal form algorithm for the Brieskorn lattice. J. Symbolic Comput. 38 (2004), no. 4, 1207–-1225.
\bibitem{Seg} E.~Segal, The closed state space of affine Landau-Ginzburg B-models.  arXiv:0904.1339v2
\bibitem{Shk} D.~Shklyarov, Hirzebruch-Riemann-Roch theorem for differential graded algebras. arXiv:0710.1937. 
\bibitem{Ts} B.~Tsygan, On the Gauss-Manin connection in cyclic homology. Methods Funct. Anal. Topology 13 (2007), no. 1, 83-–94.
\bibitem{We} C.~Weibel, The Hodge filtration and cyclic homology. K-Theory {\bf 12} (1997), no. 2, 145-–164.
\end{thebibliography}
\end{document}